\documentclass[review]{elsarticle}
\usepackage[colorlinks,linkcolor=blue,citecolor=cyan,anchorcolor=blue]{hyperref}
\usepackage{lineno}
\modulolinenumbers[5]

\usepackage{datetime}
\usdate 

\journal{Naval Research Logistics}

\usepackage[letterpaper]{geometry}
\usepackage{amsmath}
\usepackage{amsfonts, amsthm, bbm, dsfont} 

\newtheorem{innercustomthm}{Lemma}
\newenvironment{customthm}[1]
  {\renewcommand\theinnercustomthm{#1}\innercustomthm}
  {\endinnercustomthm}

\usepackage{color,xcolor,graphicx,subfig}
\usepackage{enumerate}
\usepackage{wrapfig}
\usepackage{tikz}
\usepackage{fixltx2e}
\usepackage{listings}
\usepackage{mdframed}
\usepackage{placeins}
\newcommand{\gbc}{C_0}
\newcommand{\R}{\mathbb{R}}

\newcommand{\E}{\mathbb{E}}
\newcommand{\lipone}{\text{\rm Lip(1)}}
\newcommand{\Prob}{\mathbb{P}}
\providecommand{\abs}[1]{\left\lvert#1\right\rvert}
\providecommand{\norm}[1]{\lVert#1\rVert}
\newtheorem{lemma}{Lemma}

\newtheorem{theorem}{Theorem}

\numberwithin{equation}{section}

\newtheorem{remark}{Remark}

\newcommand{\blue}[1]{\textcolor{black}{#1}}

\newtheorem{lemma-sec}{Lemma}[section]
\newtheorem{theorem-sec}{Theorem}[section]

\usepackage{url}
\makeatletter
\g@addto@macro{\UrlBreaks}{\UrlOrds}
\makeatother
\usepackage{verbatim}
\usepackage{subfig}

\usepackage{chngcntr}

\usepackage{subfiles}

\usepackage{etoolbox}
\makeatletter
    \patchcmd{\tnotemark}{\ding{73}}{1}{}{\@latex@error{Failed to path \string\tnotemark\space for \string\ding{73}}}
    \patchcmd{\tnotetext}{\ding{73}}{1}{}{\@latex@error{Failed to path \string\tnotetext\space for \string\ding{73}}}
\makeatother

\begin{document}

\begin{frontmatter}

\title{Steady-state Diffusion Approximations for Discrete-time Queue \\
in Hospital Inpatient Flow Management\tnoteref{grant}}

\author{Jiekun Feng} 
\address{Department of Statistical Science, Cornell University}
\ead{jf646@cornell.edu}

\author{Pengyi Shi}
\address{Krannert School of Management, Purdue University}
\ead{shi178@purdue.edu}

\tnotetext[grant]{Supported in part by NSF Grant CMMI-1335724 and CMMI-1537795.}

%

\begin{abstract}
In this paper, we analyze a discrete-time queue 
that is motivated from studying hospital inpatient flow management, 
where the customer count process captures the midnight inpatient census.
The stationary distribution of the customer count  
has no explicit form and is difficult to compute in certain parameter regimes. 
Using the Stein's method framework, 
we identify a continuous random variable to approximate the steady-state customer count. 
The continuous random variable corresponds to the stationary distribution 
of a diffusion process with \emph{state-dependent} diffusion coefficients. 
We characterize the error bounds of this approximation 
under a variety of system load conditions -- from lightly loaded to heavily loaded. 
We also identify the critical role that the service rate plays 
in the convergence rate of the error bounds. 
We perform extensive numerical experiments to support the theoretical findings 
and to demonstrate the approximation quality.  
In particular, we show that our approximation performs better than those based on constant diffusion coefficients
when the number of servers is small, 
\blue{which is relevant to decision making in a single hospital ward.} 


\quad

\end{abstract}


\end{frontmatter}

\linenumbers

\section{Introduction}
\label{sec:intro}

In this paper, we analyze a $GI/Geo/N$ discrete-time queue, or \emph{discrete queue} in short. 
This discrete queue has $N$ identical servers and a buffer that can hold infinitely many customers.
Customer arrivals and departures occur at discrete time epochs $k=0, 1, 2, \dots$.
At each epoch $k$, a total number of $D_k$ customers depart from the queue first,
and then a total number of $A_k$ customers arrive. 
If there are enough servers, we admit all waiting customers (if any) and new arrivals into service;
otherwise, we admit as many customers as possible, 
following the first-come-first-served queueing discipline, 
until all servers are occupied and hold the remaining customers in the buffer.
For the arrival process, we assume that
$\{A_k, k=0,1\dots\}$ forms a sequence of independent and identically distributed (i.i.d.) random variables.
For the departure process, 
we assume that each customer in service at the beginning of epoch $k$ 
(excluding new arrivals) 
has a constant probability $\mu \in (0,1)$ of departing in epoch $k$.
It is equivalent to assuming that the ``service time'' in this discrete queue 
follows a \emph{geometric} distribution with the success parameter being $\mu$, 
which is why we use ``Geo'' in the notation of the queue;   
see~\cite{daishi2015c} for a more rigorous proof on this equivalence using a coupling argument.

This discrete queue is motivated from 
studying the inpatient midnight census in hospitals~\cite{daishi2015c}, 
where the servers correspond to inpatient beds, 
and customer arrivals and departures correspond 
to patient bed-requests and discharges in a day;
see more motivation in Section~\ref{sec:mot}.
We focus on the customer count process and analyze its steady-state performance. 
Let $X_k$ denote the total number of customers in system at the beginning of epoch $k$,
including both the customers in service and those waiting in the buffer.
Under our arrival and departure assumptions, 
the customer count process $X=\{X_k: k = 0,1,\dotsc\}$
forms a discrete-time Markov chain (DTMC) and is characterized by the following relationship:
\begin{equation}
    X_{k+1} = X_k + A_k - D_k, \quad k=0,1, \dots.
\label{eq:evol}
\end{equation}
Here, the total number of departures $D_k$ 
follows a binomial distribution with parameters $(Z_k,\mu)$,
where $Z_k \equiv X_k \wedge N$ is the number of busy servers at the beginning of epoch $k$
with $\wedge$ denoting the minimum between two real numbers. 
In the rest of this paper, we focus on the \emph{Poisson} arrival case,
that is, $A_k$ follows a common Poisson distribution with mean $\Lambda$.
We specify the treatments for general arrival distributions
in an online supplement~\cite{FengShi2016a} 
and establish the corresponding error bounds there to keep this paper focused.

When
\begin{equation}
\Lambda < N\mu, \quad \textrm{or equivalently,} \quad R \equiv \frac{\Lambda}{\mu} < N,
\label{eq:stable-condition}
\end{equation}
the DTMC $X$ has a unique stationary distribution $\pi$~\cite{daishi2015c}.
Here, $R$ is the \emph{offered load} of the discrete queue.
We use $X_\infty$ to denote the steady-state customer count,
and correspondingly, its distribution is $\pi$.
We also define a \emph{scaled} version of $X_\infty$ as
\begin{equation}
\tilde X_\infty = (X_\infty - R)/\sqrt R .
\label{eq:scaled_count}
\end{equation}


\subsection{Results summary} 
\label{subsec:results_summary}
We identify a continuous random variable (r.v.) $Y_\infty$ 
to approximate the scaled steady-state customer count $\tilde X_\infty$. 
We establish bounds on the approximation errors under various system conditions. 
Specifically, $Y_\infty$ is defined by the density 
\begin{align} 
p(x) &= \frac{\kappa}{a(x)}\exp\left(\int^x_0
\frac{2b(y)}{a(y)}\mathrm{dy}\right),
\quad x\in \mathbb{R},
\label{eq:nu}
\end{align}
where $\kappa > 0$ is a normalizing constant such that $\int^\infty_{-\infty} p(x) dx =1$,
\begin{align}
b(x) &\equiv\mu\left[(x+\zeta)^--\zeta^-\right] =
\begin{cases}
-\mu x, \quad x \leq -\zeta,\\
\mu \zeta, \quad x \geq -\zeta,
\end{cases}
\label{eq:bx}
\end{align}
\begin{align}
a(x) &\equiv 
\begin{cases}
\mu(1 + \Lambda), \quad x \leq -\sqrt{R}, \\
\mu\big(2-\mu + \delta (1-\mu)x + \mu x^2\big), \quad x \in [-\sqrt{R}, -\zeta], \\
\mu \big( 2-\mu + \delta(1-\mu) \abs{\zeta} + \mu \zeta^2 \big), \quad x \geq -\zeta,
\end{cases}
\label{eq:aform}
\end{align}
and 
\begin{align}
    \zeta 
    &\equiv (R - N)/\sqrt R<0 .
    \label{eq:zeta}
\end{align} 
\blue{The density of $Y_\infty$ corresponds to the stationary distribution 
of a diffusion process with a drift $b(x)$ and a diffusion coefficient $a(x)$.
Both $a(x)$ and $b(x)$ are continuous, and come from the Stein's method framework to be detailed in Section~\ref{sec:pf}.}

Now we present our main results. 
Consider any given $N \geq 1$, $\Lambda > 0$, and $\mu \in (0,1)$ such that $1\leq R < N$ and 
\begin{equation} 
\mu = \gamma R^{-s}, \quad N - R = \beta R^q 
\label{eq:char}
\end{equation}
for some constants $\gamma,\ \beta >0$ and $s, \ q \geq 0$, 
where $s$ characterizes the rate of $\mu$ converging to 0, and $q$ characterizes the system load condition. 
Note that $q=0, \ 1/2, \ 1$ corresponds to the non-degenerate slowdown (NDS)~\cite{Atar2012}, quality-and-efficiency-driven (QED),
and quality-driven (QD) regimes~\cite{GansKoolMand2003}, respectively; 
\blue{also see~\cite{GarnMandReim2002} for a description on the QD and QED regime.} 
We establish the following error bounds on the Wasserstein distance between $\tilde X_\infty$ and $Y_\infty$, 
defined as
\begin{equation}
d_{W}\left(\tilde X_\infty, \ Y_\infty \right)   \equiv
\sup_{h\in \lipone}\abs{\E h(\tilde X_\infty) - \E h(Y_\infty)} 
\label{eq:dtmc_1}
\end{equation}
with $\lipone = \{h:\mathbb R\to\mathbb R,\,\abs{h(x) - h(y)}\leq \abs{x-y} \text{ for all }x,y\in \mathbb{R}\}$.

\begin{theorem}
\label{thm:dtmc_1}
For any $s\in [1/2, 1)$ and $q\in[1/2, 1]$ such that~\eqref{eq:char} is satisfied,
\begin{equation}
d_{W}\left(\tilde X_\infty, \ Y_\infty \right)  \leq C_1(\gamma, \beta) R^{-s/2}. 
\label{eq:dtmc_1}
\end{equation}
\end{theorem}
\begin{theorem}
\label{thm:dtmc_2}
For any $s\geq 1$ and $q\in[0, 1]$ such that~\eqref{eq:char} is satisfied 
and $R \geq 2\gamma$,
\begin{equation}
d_{W}\left(\tilde X_\infty, \ Y_\infty \right)   \leq C_2(\gamma, \beta) R^{-1/2}. 
\label{eq:dtmc_2}
\end{equation} 
\end{theorem}
Here, $C_1(\gamma, \beta)$ and $C_2(\gamma, \beta)$
are two constants that only depend on $\gamma$ and $\beta$. 
We give more interpretations of these two theorems below 
in Section~\ref{sec:main-contri}. Note that convergence in the Wasserstein distance 
implies convergence in distribution~\cite{GibbSu2002}.

\subsection{Main contributions}
\label{sec:main-contri}

Comparing to the recent series of papers on Stein's method for steady-state approximation,
our paper makes the following contributions.

\begin{itemize}

\item The density of $Y_\infty$ corresponds to the steady-state distribution 
of a diffusion process with a \emph{state-dependent} diffusion coefficient $a(x)$
and a piece-wise linear drift $b(x)$. 
Note that in~\cite{daishi2015c}, the most relevant paper, 
the authors studied the same DTMC $X$ and promoted the use of a continuous r.v. 
with similar state-dependent diffusion coefficients to approximate $\tilde X_\infty$.  
However, they were not able to establish the error bounds for such approximation 
-- Theorem~3 there only established the error bounds between $\tilde X_\infty$ 
and a continuous r.v. with a \emph{constant} diffusion coefficient, 
which we denote as $Y^0_\infty$ in this paper.
Comparing to $Y^0_\infty$, $Y_\infty$ is equally easy to evaluate numerically 
using the explicit form in~\eqref{eq:nu}, 
while it produces better approximation for $\tilde X_\infty$, 
especially when the system size is small.
For example, when $N=18$, 
$Y_\infty$ can reduce the relative approximation error in the expected queue length 
by as much as 10\% than using $Y^0_\infty$.  
This paper fills the gap in establishing the error bounds 
between $\tilde X_\infty$ and the state-dependent r.v. $Y_\infty$.  
As we will further illustrate below, this is not a trivial extension 
-- it is challenging to establish bounds involving the density of $Y_\infty$ 
since its form is more complicated than that of $Y^0_\infty$.

\item We characterize the convergence rate of the error bounds
under different system load conditions (reflected by $q$) 
and rates of $\mu$ converging to $0$ (reflected by $s$).  
Theorem~\ref{thm:dtmc_1} says that when $\mu$ converges to $0$ in a rate 
that is between $1/\sqrt{R}$ and $1/R$, 
and when the system runs in the QD, QED, or any regime in between,  
the error bound converges to $0$ in a rate that is between $1/\sqrt[4]{R}$ and $1/\sqrt{R}$.
Theorem~\ref{thm:dtmc_2} says that when $\mu$ converges to $0$ in a rate that
is faster than or equal to $1/R$, and when the system runs in any regimes
from QD to NDS, the error bound converges to $0$ at a \emph{constant} rate $1/\sqrt R$.
Comparing to Theorem 3 in~\cite{daishi2015c} 
that only established the error bounds in the QED regime,  
our results here cover a much wider range of $q$, 
which justify using $Y_\infty$ to approximate $X_\infty$  
under a variety of system load conditions. 
This can have more practical impact since, for example, 
the utilization of different hospital wards can vary greatly from lightly loaded (60\%) 
to very heavily loaded ($>95\%$); see for example, Table~5 in~\cite{ShiDaietAl2014b}, 
and Tables 89 and 91 in~\cite{NCHS2016}.

\item Theorems~\ref{thm:dtmc_1} and~\ref{thm:dtmc_2} also reveal the critical role 
that the service rate $\mu$ plays in the convergence of the error bound in the discrete queue.
That is, to ensure that the error bound goes to 0,
we require that the number of servers $N$ goes to $\infty$
and the service rate $\mu$ goes to $0$ \emph{at the same time},
\blue{which is supported by our numerical results.}  
This is a main difference from the continuous-time queueing systems 
studied in the literature, including~\cite{Brav2017},  
where the authors develop a state-dependent diffusion model for Erlang-C queue.  
In the continuous-time queues, just $N, R \rightarrow \infty$ ensures the convergence of the error bound.  
In addition, in the discrete queue the rate of $\mu$ converging to 0 needs to be \emph{fast} enough
to ensure convergence in different operating regimes, 
e.g., in the NDS regime, we find that the error bound does \emph{not} converge
when $s = 1/2$, but does converge when $s=1$; see numerical results in Section~\ref{sec:num}.
Note that this asymptotic regime of $\mu$ going to $0$ is also of practical relevance for hospital setting 
since a typical inpatient stay (service time) is 4-5 days~\cite{daishi2015c}, 
i.e., the service rate is small.

\item \blue{For the proof, we overcome a major challenge in establishing the error bounds 
due to the intrinsic difference between continuous- and discrete-time queues. 
Unlike our continuous counterpart -- the Erlang-C queue --  
where the customer count process is a birth-death process 
with only one arrival or departure occurring at each transition, 
in the discrete queue a number of arrivals and discharges could occur at each transition. 
The latter makes the proof of the moment bounds complicated, 
especially when the system is heavily loaded (NDS regime);
see details in Section~\ref{sec:proof-last-error}. 
In addition, since the diffusion coefficient $a(x)$ is state-dependent,  
the density $p(x)$ becomes complicated. 
We apply a new technique developed in~\cite{Brav2017} to establish the gradient bounds. 
This application is not trivial since our diffusion coefficient has a quadratic form, 
while the one in~\cite{Brav2017} has a linear form; 
see details in Section~\ref{sec:proof-gradient}. 
}

\end{itemize}

Two more things worth noting here. 
First, that there is no explicit formula to calculate $\pi$
despite the simplicity of the dynamic equation~\eqref{eq:evol}.
Using the standard Markov chain technique to solve $\pi$ \blue{imposes challenges for practitioners}. 
It can be also time consuming 
since generating the transition matrix requires calculation of the convolution between $A_k$ and $D_k$,
where the distribution of the latter depends on the number of customers $X_k$.
Even for realistic hospital settings, it could take several hours, sometimes days, to get $\pi$;
see computational results summarized in Table~\ref{tab:comp_time_1}. 
The conventional generating-function method does not work well either, 
because the distribution of $D_k$ depends on $X_k$ instead of following a common distribution
(which differentiates Equation~\eqref{eq:evol} from the dynamics of the $M/GI/1$ queue); 
see more discussions in Section~\ref{sec:lit}.
In contrast, the density function of $Y_\infty$ has an explicit form and takes almost no time to evaluate. 
\blue{Moreover, we demonstrate through this paper that the Stein's method framework provides not only a powerful tool to characterize the error bounds, 
but also an engineering tool to \emph{identify} a good approximation for the steady-state distribution. 
The latter is particularly helpful for systems without known or explicit steady-state distributions, such as our discrete queue.  
}

\blue{Second, note that we fix the time unit as one day,  
and let the service rate $\mu$ go to 0 in the limit regime. 
In the queueing literature, however, the service rate/time is usually fixed and normalized to 1.  
We would like to emphasize that this normalization can also be done in our setting 
by choosing the time unit to be $1/\mu$. 
For example, consider two discrete queues with $\mu_1=1/7$ days and $\mu_2=1/30$ days, 
corresponding to a mean service time of one week and one month, respectively. 
To normalize the service rate to $\tilde \mu_1 = \tilde \mu_2 = 1$, 
we set the time unit as one week in the first system, and as one month in the second system. 
Now we need to interpret the customer count process $\{X_k\}$ carefully under the new time unit: 
$X_k$ represents the customer count at the beginning of the $k$th time interval, 
where each interval lasts for $\mu_1=1/7$ unit in the first system 
and $\mu_2=1/30$ unit in the second system. 
In other words, if we view $\{X_k\}$ as a sampling process from a underlying continuous process
with samples taken at the beginning of each interval $k$, 
we sample the first system 7 times in one time unit (a week), and the second system 30 times in one time unit (a month). 
The service rate $\mu$ is essentially equivalent to the width of the sampling window. 
Thus, if we fix $\mu=1$ as in the literature, we require the sampling window to go to 0, 
or say, the sampling frequency $\rightarrow \infty$ in the limiting regime. 
It also becomes intuitive why we require this for the convergence of error bounds, 
as we sample more frequently,  
the discrete queue becomes closer to the (continuous-time) diffusion process. 
In this paper, for ease of exposition and a better connection with the hospital background, 
we choose to fix the time unit as 1 day 
but let the service rate $\mu$, with respect to this time unit, decrease to 0. 
}


The rest of the paper is organized as follows.
In Section~\ref{sec:mot},
we discuss the motivation of the discrete queue
from hospital inpatient flow management
and its broader application in telecommunication.
In Section~\ref{sec:lit} we review relevant papers in the literature.
In Section~\ref{sec:pf} we prove our main theorems. 
In Section~\ref{sec:pfLem}, we detail the proof
for two important lemmas to establish the error bounds. 
We describe the numerical results in Section~\ref{sec:num}
and conclude the paper in Section~\ref{sec:conclusion}.

\subsection{Motivation of the discrete queue}
\label{sec:mot}

The discrete queue is motivated by studying hospital inpatient flows~\cite{Shietal2014}.
The inpatient beds are modeled as the servers, 
and patients who need to be admitted to an inpatient bed
are modeled as customers, for example,
patients who have received treatment in the emergency department (ED)
and wait to be hospitalized -- commonly known as the ED boarding patients.
The customer count $X_k$ corresponds to the \emph{midnight} census at day $k$,
i.e., the number of patients who are occupying an inpatient bed
or waiting to be admitted at the midnight of day $k$.
Naturally, $A_k$ and $D_k$ correspond
to the total number of patient arrivals and discharges within day $k$, respectively,
and the midnight census at the next day, $X_{k+1}$, evolves as in~\eqref{eq:evol}.
Empirical studies suggest that the bed-request process of the ED boarding patients
can be modeled by a periodic Poisson process with the period being one day~\cite{Shietal2014, Armoetal2015}.
Thus, it is reasonable to assume that the daily arrival $A_k$ follows a Poisson distribution,
with $\Lambda$ corresponding to the daily arrival rate.
However, the length-of-stay (LOS) distribution is usually \emph{not} geometric. 
Nevertheless, the system performance is not very sensitive to the LOS distribution 
when the utilization is not extremely high; 
see Section 4.7 of~\cite{Shi2013}.
Therefore, we focus on the geometric setting in this paper for tractability.

The midnight census is a key performance metric
monitored by many hospitals~\cite{SimoYankDunt2010}.
Moreover, getting the stationary distribution for the midnight census is a crucial step
in predicting the time-of-day patient census, 
since the census at a certain hour $t$ equals the sum of the midnight census 
and the difference between the number of arrivals and discharges from the midnight to hour $t$; 
see the two-time-scale framework developed in~\cite{daishi2015c}. 
Also see there for more details on the importance of studying the midnight census 
as well as further justifications on model assumptions.


Besides the applications in the healthcare setting,
discrete queueing systems have been motivated from a variety of applications
in the fields of telecommunication and computer systems,
in which the time is usually divided into fixed-length time slots.
For example, multi-server discrete queues with geometric service times are studied
in the context of circuit-switched multiple-access communications channel~\cite{RubiZhan1991}
and systems with randomly interrupted servers~\cite{LaevBrun1995}
(which are equivalent to systems
with non-interrupted servers and geometric service time).
Thus, the analysis and results we gain in this paper can potentially benefit a larger community.

\subsection{Literature review}
\label{sec:lit}

Stein's method is a well known method for establishing error bounds
in various fields and applications; see, for example, the survey papers~\cite{Chat2014,Ross2011}.
The proof in our paper is mainly based on the framework developed in~\cite{BravDai2015} 
on applying Stein's method to steady-state diffusion approximations in $M/Ph/N + M$ queue;
also see a tutorial on applying this framework to Erlang-A and Erlang-C queues
in~\cite{BravDaiFeng2015} and the references there for this line of work. 
\blue{Another relevant work is done by Gurvich~\cite{Gurv2014}, who independently develops a method to prove a steady-state convergence that is similar to Theorem 1 in~\cite{BravDai2015}. In particular, his method also uses diffusion models, instead of a diffusion limit, to approximate continuous time Markov chains in steady state. Furthermore, his method includes the derivation of gradient bounds for Poisson equations as well as the moment bounds for CTMCs; these bounds are closely related with the Stein's method framework and has largely inspired the work in~\cite{BravDai2015}.} 
In addition to error bounds on the customer count distributions,~\cite{HuanGurv2016} apply the Stein's method framework
to the $M/GI/1 + GI$ queue and establish error bounds 
for a variety of performance measures including the waiting time and abandonment.

\blue{Among this line of research, the most relevant work is Chapter 3 of doctoral thesis~\cite{Brav2017}
(originally appeared as a working paper~\cite{BravDai2016}). 
The author develops a new diffusion model with state-dependent diffusion coefficients for Erlang-C queue;
this refined diffusion model allows establishing an error bound with a higher order convergence rate. 
Indeed, this work is a main driver for us to look into the state-dependent diffusion approximation in the discrete queue. 
As discussed in Section~\ref{sec:main-contri}, the limit regime we identify is different from that in~\cite{Brav2017},
and we overcome a major challenge in proving the moment bounds that uniquely arises in our discrete setting. 
Note that to establish the gradient bounds in Section~\ref{sec:proof-gradient}, we build on a new technique developed in~\cite{Brav2017}. 
Applying this technique to our setting, though, is not trivial since our diffusion coefficient $a(x)$ is different from that in~\cite{Brav2017}.  
Specifically, the diffusion coefficient in~\cite{Brav2017} takes the form  
\[ 
\mu\big ( 2 + \delta( x \wedge \abs{\zeta} ) \big),
\]
which differs from $a(x)$ in Equation~\eqref{eq:aform} by a term of $\mu^2 ( x\wedge \abs{\zeta} )^2$. 
This difference is caused by the different generators of Erlang-C queue and our discrete queue; 
see Section~\ref{sec:pf} in our paper and Section~3.2.2 in~\cite{Brav2017} for the derivation. 
Because of the quadratic term in our diffusion coefficient $a(x)$, we need to prove Lemmas~\ref{lem:ints_for_grad_bdds_} --~\ref{lem:a_} 
to prove Lemma~\ref{lem:gradboundsCW}, instead of directly applying the results in~\cite{Brav2017}. 
In particular, the quadratic term raises additional complexity 
when we bound~\eqref{eq:abdd_2_} and~\eqref{eq:abdd_3_} in Lemma~\ref{lem:a_},~\eqref{eq:ey__} in Lemma~\ref{lem:ey_}, and~\eqref{eq:int-2_} in Lemma~\ref{lem:ints_for_grad_bdds_}.}


Regarding the modeling side, various discrete queues have been studied
in the area of telecommunication and computer systems,
with different assumptions on the arrival and service time distribution,
server numbers (single or multiple), and buffer capacity;
see a detailed summary in Section~1.2.7 of~\cite{BrunKim1993}
for some of the early works and~\cite{GaoWittBrun2004} for more recent development.
The most relevant work is~\cite{GaoWittBrun2004},
where the authors study a same $GI/Geo/N$ queue.
In that paper, the authors employ the conventional generating-function method
to perform steady-state analysis. 
However, numerically implementing the generating-function method to find $\pi$
involves finding $N-1$ roots inside the unit disk from an $N$th order nonlinear equation,
and it becomes computationally difficult when $N$ is large.
Indeed, in their numerical experiments, 
the largest $N$ tried by the authors is 16.
In this paper, we develop an efficient and accurate way to approximate $\pi$,
and more importantly, we are able to provide error bounds on such approximation.~\cite{JansLeeuZwar2008} study the $M/D/s$ queue 
where the dynamic equation for the number of waiting customers has a similar form to (1.1), 
but intrinsically is different. 
Their objective is similar to ours, 
that is, finding more accurate approximations for performance metrics of interest, 
with a focus on the QED regime.  
The method used in~\cite{JansLeeuZwar2008} involves infinite series expansion, 
which requires detailed derivations tied to the queueing models, performance metrics,  
and operating regimes; see~\cite{BravDaiFeng2015} for more discussions 
on comparing this method with the Stein's method framework. 


%


\section{Proof of Theorems~\ref{thm:dtmc_1} and~\ref{thm:dtmc_2}}
\label{sec:pf}

To prove Theorems~\ref{thm:dtmc_1} and~\ref{thm:dtmc_2},
we employ the Stein's method framework for steady-state approximation.
The major components of this framework are 
Poisson equation, generator coupling, gradient and moment bounds;
see~\cite{BravDai2015} for a systematical description.
Throughout the rest of the paper, we define
\[
\delta \equiv 1/\sqrt R 
\]
for notational convenience.

\begin{proof}

Define $G_Y$ as
\begin{equation}
    G_Y f(x) = b(x) f'(x) + \frac{1}{2} a(x) f''(x), 
    \quad x\in\mathbb{R},\,f\in C^2(\mathbb{R}),
    \label{eq:gy}
\end{equation}
which is the generator of a diffusion process
with diffusion coefficient $a(x)$ and drift $b(x)$.
As noted before, the stationary distribution of this diffusion process
has a density given by~\eqref{eq:nu}.

Now, let $f=f_h$ be a solution to the \emph{Poisson equation}
 \begin{equation}
 G_Y f(x) = \E \left[h(Y_\infty)\right]-h(x),\quad x\in \mathbb R.
 \label{eq:poisson}
 \end{equation}
Lemma~\ref{lem:gradboundsCW} below shows that
$f$ is twice continuously differentiable, with an absolutely continuous second derivative.

Next, we do a \emph{generator coupling} via~\eqref{eq:poisson}.
The generator of the scaled DTMC $\tilde X$ is
\begin{equation}
    G_{\tilde X}f(x) = \E_n\left[f(x+\delta(A_0-D_0))-f(x)\right] , \quad
    x=\delta(n-R),\,n = 0,1,\dotsc, \,f\in C^2(\mathbb{R}) ,
    \label{eq:gx}
\end{equation}
where $\E_n$ is the expectation under $\Prob_n$,
the conditional probability distribution given that the starting customer count
equals $n$, $A_0 \sim$ Poisson with mean $\Lambda$,
and $D_0 \sim$ binomial with $(n\wedge N, \mu)$
and is independent of $A_0$.
It is proven in~\cite{daishi2015c} that
\begin{equation}
    \E \left[G_{\tilde{X}}f(\tilde X_\infty)\right] = 0.
    \label{eq:gx_prop}
\end{equation}
Taking expectation with respect to $\tilde X_\infty$ on both sides of \eqref{eq:poisson},
we have that
\begin{eqnarray}
    \E\left[h(Y_\infty)\right] - \E\left[h(\tilde X_\infty)\right]
    &=& \E \left[G_Y f(\tilde X_\infty)\right]\notag\\
    &=& \E \left[G_Y f(\tilde X_\infty) - G_{\tilde{X}}f(\tilde X_\infty)\right] ,
    \label{eq:gcouple}
\end{eqnarray}
where the second equality comes from~\eqref{eq:gx_prop}.
To bound the right side of \eqref{eq:gcouple}, we perform the Taylor expansion for $G_{\tilde{X}}f(x)$,
for any given $x = \delta(n - x_\infty)$ and $n=0,1,\dotsc$,
\begin{align}
    G_{\tilde{X}}f(x)&=f'(x) \delta \E_n(A_0 -D_0) +\frac{1}{2}f''(x)\delta^2 \E_n[(A_0-D_0)^2]+
    \frac{1}{2}\delta^2 \E_n[ \left(f''(\eta)-f''(x)\right)(A_0-D_0)^2]\notag\\
    &= G_Y f(x) 
+ \frac{1}{2}\delta^2 \E_n[\int^\eta_x f'''(y)\mathrm dy (A_0-D_0)^2],
\label{eq:te}
\end{align}
with
\begin{equation*}
    \abs{\eta - x}\leq \delta \abs{A_0 - D_0}.
\end{equation*}
Note that~\eqref{eq:te} follows from the absolute continuity of $f''(x)$ and the following facts:
\begin{eqnarray}
\delta \E_n(A_0 -D_0)&=&\delta [\Lambda - (n\wedge N)\mu]=
\delta (\Lambda - N\mu) + \delta (n-N)^- \mu\notag\\
&=& \mu\zeta + (x+\zeta)^- \mu
= b(x) ,
\label{eq:firstorder}
\\
\delta^2 \E_n [(A_0 -D_0)^2] &=&\delta^2 \textrm{Var}_n (A_0-D_0) +\delta^2 (\E_n(A_0-D_0))^2 \nonumber \\
&=& \delta^2 \Lambda +\delta^2 (n\wedge N)\mu(1-\mu) +b^2(x) = a(x). 
\label{eq:secondorder}
\end{eqnarray}
Note that in~\eqref{eq:aform}, 
the part of $a(x)=\mu(1+\Lambda)$ on $(-\infty, -1/\delta)$ corresponds to $n<0$. 
It is added to make the diffusion coefficient $a(x)$ exist everywhere on $\R$ and be continuous at the point $x = -1/\delta$,
so that the corresponding diffusion process is well-defined. 
Here, to derive the second equality of~\eqref{eq:secondorder},
we have used the fact that $A_0$ is Poisson so that its variance equals the mean $\Lambda$.
The treatment for the non-Poisson arrival distribution  
and the corresponding error bounds 
are detailed in the online supplement~\cite{FengShi2016a}. 
\blue{Also note that it is through the Taylor expansion of $G_{\tilde X}f(x)$ that we find the approximating random variable $Y_\infty$. 
Specifically, we use the first two terms of the Taylor expansion to find the diffusion coefficient $a(x)$ 
and the drift $b(x)$ for a diffusion process, 
where $Y_\infty$ is the steady-state random variable of this diffusion process. 
}

Combining~\eqref{eq:gcouple} and \eqref{eq:te} gives
\begin{eqnarray}
\abs{\E h(\tilde X_\infty) - \E h(Y_\infty)}
&=&\Big|\E \left[G_{\tilde{X}}f(\tilde X_\infty) - G_Y f(\tilde X_\infty)\right]\Big|\notag\\
&\leq& \frac{1}{2}\delta^2 \E\left[\epsilon(X_\infty)\right],
\label{eq:ate}
\end{eqnarray}
where 
\begin{equation}
\epsilon(X_\infty) = \E_{X_\infty}
\Big[  \Big|\int^{\tilde X_\infty+\delta(A_0-D_0)}_{\tilde X_\infty} \abs{f'''(y)} \mathrm d y\Big|(A_0-D_0)^2 \Big], 
\label{eq:def-epsilon}
\end{equation}
and we use $\E_{X_\infty}$ to denote 
the expectation conditioning on the starting customer count $X_\infty$ in the rest of the paper.
Then, Theorems~\ref{thm:dtmc_1} and~\ref{thm:dtmc_2} 
follow from part (a) and part (b) of Lemma~\ref{lem:last_error_dtmc} below, respectively.
\end{proof}

\begin{lemma}
\label{lem:last_error_dtmc} 
(a) For an $s\in [1/2, 1]$ and any $q\in[1/2, 1]$ such that~\eqref{eq:char} is satisfied,
\begin{align}
\frac{1}{2}\delta^2 \E\left[\epsilon(X_\infty)\right]\leq C_1(\gamma,\beta) R^{-s/2} . 
\label{eq:last_error_dtmc_QDQED}
\end{align}

(b) For any $s\geq 1$ and any $q\in[0, 1]$ such that~\eqref{eq:char} is satisfied 
and $R\geq 2\gamma$, 
\begin{align}
\frac{1}{2}\delta^2 \E\left[\epsilon(X_\infty)\right]\leq 
C_2(\gamma,\beta) R^{-1/2}. 
\label{eq:last_error_dtmc_QDNDS}
\end{align}
\end{lemma} 
Here, $C_1(\gamma,\beta)$ and $C_2(\gamma,\beta)$ 
are two constants depending only on $\gamma,\,\beta$.
Note that the condition $R\geq 2\gamma$ in part (b) implies $\mu\leq 1/2$, 
which is a realistic assumption for hospital inpatient setting 
since the typical service time is longer than 2 days. 
\blue{However, this condition can be relaxed -- 
for a given $\mu \in (0, \mu_0)$ for some constant $\mu_0$, 
we can prove a similar result as that in~\eqref{eq:dtmc_2} 
with the constant $C_2$ depending on $\gamma$, $\beta$, and $\mu_0$. 
We focus on the special case of $R\geq 2\gamma$ 
to keep the statement of the results clean.}

To prove Lemma~\ref{lem:last_error_dtmc}, 
we also need the following lemma on the gradient bounds, 
whose proof is detailed in Section~\ref{sec:pfLem}. 
\begin{lemma}
\label{lem:gradboundsCW}
Fix an $h \in \lipone$ with $h(0) = 0$.
There exists a solution to the Poisson equation~\eqref{eq:poisson}, $f_h$, 
that is twice continuously differentiable with an absolutely continuous second derivative,
for any $s\geq 1/2$ such that~\eqref{eq:char} is satisfied,
\begin{align}
\abs{f_h'''(x)} &\leq \begin{cases}
\frac{\gbc}{\mu} (1 + 1/\abs{\zeta}),
\quad x\leq -\zeta , \\
\frac{4}{\mu},
\quad x \geq  -\zeta, 
\end{cases}
\label{eq:WCder3}
\end{align}
where $f_h'''(x)$ is interpreted as the left derivative at the point $x = -1/\delta$ and $x=-\zeta$, 
and $\gbc =\gbc\left(\gamma \right)$ is a constant depending only on $\gamma$.
\end{lemma}

\section{Proof of Lemmas~\ref{lem:last_error_dtmc} 
and~\ref{lem:gradboundsCW}}
\label{sec:pfLem}

As mentioned in the introduction, 
a major challenge we overcome in proving Theorems 1 and 2  
is to bound the error terms in Lemma~\ref{lem:last_error_dtmc}. 
This challenge comes from the fact that 
the discrete queue could have multiple arrivals and departures between two consecutive epochs,
in other words, the jump size between $X_k$ and $X_{k+1}$ is unbounded. 
In addition, since the density $p(x)$ has a complicated form, 
we adapt a new technique developed in~\cite{Brav2017} to establish the gradient bounds in our setting. 
In the interest of space, we give the main steps in proving Lemmas~\ref{lem:last_error_dtmc} and~\ref{lem:gradboundsCW}
in Sections~\ref{sec:proof-last-error} and~\ref{sec:proof-gradient}, respectively, 
and leave some of the lengthy proofs to Section~1 of the online supplement~\cite{FengShi2016a}.

\subsection{Proof of Lemma~\ref{lem:last_error_dtmc}}
\label{sec:proof-last-error}

Let $f=f_h$ be a solution to the Poisson equation defined in~(\ref{eq:poisson}).
To prove this lemma, we first consider $\epsilon(X_\infty)$ 
when $X_\infty\leq N$, or equivalently, $\tilde X_\infty \leq - \zeta$,
\begin{eqnarray}
\epsilon(X_\infty)
&=&  \E_{X_\infty} \big[   \Big| \int^{\tilde X_\infty+\delta(A_0-D_0)}_{\tilde X_\infty} \abs{f'''(y)} \mathrm dy
\Big|(A_0-D_0)^2 \mathds{1}_{\{X_\infty + A_0 -D_0 \leq N\}}\big]     \nonumber\\
&&{} + \E_{X_\infty} \big[ \Big| \int^{\tilde X_\infty+\delta(A_0-D_0)}_{\tilde X_\infty} \abs{f'''(y)} \mathrm dy
\Big|(A_0-D_0)^2 \mathds{1}_{\{X_\infty + A_0 -D_0 > N\}}\big] \nonumber\\
&\leq& \E_{X_\infty} \big[  \Big|\int^{\tilde X_\infty+\delta(A_0-D_0)}_{\tilde X_\infty} 
\frac{\gbc}{\mu}(1+1/\abs\zeta)
 \mathrm d y\Big|(A_0-D_0)^2 \mathds{1}_{\{X_\infty + A_0 -D_0 \leq N\}}\big]\nonumber\\
&+& \E_{X_\infty}\big[  \Big|\int^{-\zeta}_{\tilde X_\infty} 
\frac{\gbc}{\mu}(1+1/\abs\zeta)
\mathrm dy\Big|(A_0-D_0)^2 \mathds{1}_{\{X_\infty + A_0 -D_0 > N\}}\big]\nonumber\\
&+& \E_{X_\infty} \big[  \Big|\int_{-\zeta}^{\tilde X_\infty+\delta(A_0-D_0)} 
\frac{4}{\mu}
\mathrm  dy\Big|(A_0-D_0)^2 \mathds{1}_{\{X_\infty + A_0 -D_0 > N\}}\big]\nonumber\\
&=&\delta 
\frac{\gbc}{\mu}(1+1/\abs\zeta)
 \E_{X_\infty} \Big[ \left|(X_\infty+A_0-D_0)\wedge N -  X_\infty\right|(A_0-D_0)^2 \Big]
\label{eq:ate1}\\
&+&\delta 
\frac{4}{\mu} 
\E_{X_\infty} \Big[ (A_0 -D_0 -(N - X_\infty))(A_0-D_0)^2 \mathds{1}_{\{X_\infty + A_0 -D_0 > N\}}\Big] .
\label{eq:ate2}
\end{eqnarray}
To get the first inequality above,  
we use the gradient bound for $f'''$ stated in Lemma~\ref{lem:gradboundsCW}. 
The reason we need to separately consider the two cases $X_\infty + A_0 -D_0 > N$
and $X_\infty + A_0 -D_0 \leq N$ is that the gradient bounds for $f'''$ are different
on $(-\infty, -\zeta]$ and $(-\zeta, \infty)$.
It is easy to see that the term inside the expectations in~\eqref{eq:ate1} and~\eqref{eq:ate2} 
can be bounded above by $\abs{A_0-D_0}^3$.
Therefore, \eqref{eq:ate1} and \eqref{eq:ate2} imply that
\begin{align}
& \frac{1}{2}\delta^2 \E \left[\epsilon(X_\infty) \mathds{1}_{\{X_\infty \leq N\}} \right]
\notag\\
&\leq \frac{1}{2}\delta^2 \E\left\{
\delta \frac{\gbc}{\mu}(1+1/\abs\zeta) \,\E_{X_\infty}\abs{A_0 - D_0}^3\,
\mathds{1}_{\{X_\infty \leq N\}} \right\}
 + \frac{1}{2}\delta^2\delta \frac{4}{\mu}\,\E\left[\E_{X_\infty}\abs{A_0 - D_0}^3\right]
\notag\\
&\leq
\delta^3 \frac{2}{\mu} \E\left[\E_{X_\infty}\abs{A_0 - D_0}^3\right]
+ \frac{\gbc}{2} \delta^3 \frac{1}{\mu}(1+1/\abs\zeta)
\E\left[\E_{X_\infty}\abs{A_0 - D_0}^3\mathds{1}_{\{X_\infty \leq N\}}\right].
\label{eq:ate5}
\end{align}

Similarly, when $X_\infty>N$, or equivalently, $\tilde X_\infty > - \zeta$, we get
\begin{align}
\epsilon(X_\infty)
&\leq \delta \frac{4}{\mu} \E_{X_\infty} \Big[ \left| X_\infty - (X_\infty+A_0-D_0)\vee N\right|(A_0-D_0)^2 \Big]
\label{eq:ate3}\\
&+\delta \frac{\gbc}{\mu} (1+1/\abs\zeta)\E_{X_\infty} \Big[ (D_0-A_0 -(X_\infty-N))(A_0-D_0)^2 \mathds{1}_{\{X_\infty + A_0 -D_0 \leq N\}}\Big].
\label{eq:ate4}
\end{align}
Note that~\eqref{eq:ate3} is subject to the same upper bound
as that for~\eqref{eq:ate2}, which is the first term in~\eqref{eq:ate5}. Thus,
\eqref{eq:ate3} and \eqref{eq:ate4} imply that
\begin{align}
& \frac{1}{2}\delta^2 \E\left[ \epsilon(X_\infty)\mathds{1}_{\{X_\infty > N\}}\right]
\notag\\
&\leq
\delta^3 \frac{2}{\mu} \E\left[\E_{X_\infty}\abs{A_0 - D_0}^3\right]
+
\frac{\gbc}{2} \delta^3 \frac{1}{\mu}(1+1/\abs\zeta)
\E \left\{   \E_{X_\infty}\left[\abs{A_0-D_0}^3
\mathds{1}_{\{X_\infty + A_0 - D_0 \leq N\}}\right]\mathds{1}_{\{X_\infty > N\}}\right\} 
\notag \\ 
&= \delta^3 \frac{2}{\mu} \E\left[\E_{X_\infty}\abs{A_0 - D_0}^3\right]
+
\frac{\gbc}{2} \delta^3 \frac{1}{\mu}(1+1/\abs\zeta)
\E \left\{  \abs{A-D}^3
\E_{X_\infty}\left[\mathds{1}_{\{X_\infty + A_0 - D_0 \leq N\}}\right]\mathds{1}_{\{X_\infty > N\}}\right\}
\label{eq:ate7} \\
& \leq
\delta^3 \frac{2}{\mu} \E\left[\E_{X_\infty}\abs{A_0 - D_0}^3\right]
+
\frac{\gbc}{2} \delta^3 \frac{1}{\mu}(1+1/\abs\zeta)
\E \left[   \abs{A-D}^3\right]\Prob (X_\infty + A_0 - D_0 \leq N)
\notag\\
& = \delta^3 \frac{2}{\mu} \E\left[\E_{X_\infty}\abs{A_0 - D_0}^3\right]
+
\frac{\gbc}{2} \delta^3 \frac{1}{\mu}(1+1/\abs\zeta)
\E \left[   \abs{A-D}^3\right]\Prob (X_\infty \leq N) .
\label{eq:ate8}
\end{align}
To get~\eqref{eq:ate7}, we use the fact that when $X_\infty > N$,
the number of discharges is $Bino(N,\mu)$,
which is independent of $X_\infty$ and $X_\infty + A_0-D_0$.
We adopt the two new notations for the arrival and discharge quantities,
$A\sim Poiss(\Lambda)$ and $D\sim Bino(N,\mu)$,
to emphasize their independence on $X_\infty$.
Consequently, we can decouple $|A - D|^3$ from $\E_{X_\infty}$.
For the last equality~\eqref{eq:ate8}, we use the fact that the system is in the steady-state
so that the next period customer count $X_\infty + A_0-D_0$ has the same distribution
as the current count $X_\infty$.

Next, we specify how to bound~\eqref{eq:ate5} and~\eqref{eq:ate8} 
under different conditions of $s$ and $q$, 
that is, $s\in [1/2,1],\,q\in [1/2,1]$ for part (a) of Lemma~\ref{lem:last_error_dtmc}, 
and $s\geq 1,\,q\in [0,1]$ for part (b) of Lemma~\ref{lem:last_error_dtmc}. 
\blue{We address two different challenges associated with proving these two parts, 
one is on bounding $\E_{X_\infty}\left[\abs{A_0-D_0}^3\right]$ for part (a), 
and the other is on bounding the idle probability $\Prob(X_\infty \leq N)$ for part (b).}

\subsubsection{Part (a) of Lemma~\ref{lem:last_error_dtmc}}

\blue{To bound~\eqref{eq:ate5} and~\eqref{eq:ate8} when $s\in [1/2,1],\,q\in [1/2,1]$,   
we need a special treatment to bound $\E_{X_\infty}\left[\abs{A_0-D_0}^3\right]$
since it is in the order of $\mu N$, and under the condition $s\in [1/2, 1)$,
$\mu N$ increases to $\infty$ when $N\rightarrow\infty$.}
Specifically, let $X_\infty = n$. For any given $n$,  
\begin{align*}
\E_n\abs{A_0-D_0}^3 &=\E_n[\abs{(A_0-\Lambda)+(z(n)\mu - D_0)+(\Lambda - z(n)\mu)}^3]\\
&\leq 9\E_n \abs{A_0-\Lambda}^3+9\E_n \abs{D_0-z(n)\mu}^3+9\E_n \abs{\Lambda -z(n)\mu}^3\\
&\leq 27\max\left( \Lambda ,\, \Lambda^{3/2} \right)+27\max\left( z(n)\mu ,\, (z(n)\mu )^{3/2} \right)   +9\abs{\Lambda -z(n)\mu}^3,
\end{align*}
where $z(n) = n\wedge N$. 
To get the first inequality, we use the $c_r$-inequality, 
while to get the second inequality, we use Lemma~\ref{lem:rvcentralmoments_} stated at the end of this section.

Taking expectation with respect to $X_\infty$, we obtain 
\begin{align}
\E\left[\E_{X_\infty}\abs{A_0-D_0}^3\right] &< 54 \max\left(N\mu,\,(N\mu)^{3/2}\right)+9\E\abs{z(X_\infty)\mu-\Lambda}^3,
\label{eq:3rd_order_mom}
\end{align}
since $\Lambda < N\mu$ from the stability condition.  
Recall that $\mu = \gamma R^{-s} = \gamma \delta^{2s}$ and $N = R + \beta R^q$.  
When $s \in [1/2,1]$ and $q\in [1/2,1]$, we have 
\begin{align}
54\max\left(N\mu,\,(N\mu)^{3/2}\right) \leq 54\gamma(1+\beta)[1\vee \sqrt{\gamma(1+\beta)}] \delta^{3(s-1)},
\label{eq:challenge0}
\end{align}
{where we use the characterizations of $N$ and $\mu$ in~\eqref{eq:char} to get the bound on $N\mu$}
\begin{align*}
    N\mu = \gamma (N\delta^2) \delta^{2s-2}\leq \gamma (\beta+1)\delta^{2s-2}. 
\end{align*}

To bound the second term on the right side of~\eqref{eq:3rd_order_mom}, we have 
\begin{align}
\E\abs{z(X_\infty)\mu-\Lambda}^3&=\E \abs{N\mu-\Lambda -\mu(X_\infty-N)^-}^3
=\delta^{-3}\E\big|b(\tilde X_\infty)\big|^3\nonumber\\
&=\delta^{-3} \mu^3 \E[\big|\tilde X_\infty^3\mathds{1}_{\{\tilde X_\infty<-\zeta\}}\big|]
+ \delta^{-3}\mu^3 \abs\zeta^3\Prob(\tilde X_\infty\geq -\zeta),
\label{eq:challenge_}
\end{align} 
where bounding $\E[\big|\tilde X_\infty^3\mathds{1}_{\{\tilde X_\infty<-\zeta\}}\big|]$ 
and $\abs\zeta \Prob(\tilde X_\infty\geq -\zeta)$ relies on the two following lemmas.

\begin{lemma-sec}[Partial third moment bound]
\label{lem:3rd_order_mom_QDQED_}
Consider the DTMC $X$ and the scaled version $\tilde X$. 
For all $N \geq 1$, $\Lambda > 0$, and $\mu \in (0,1)$ that satisfy $1 \leq R < N$
and~\eqref{eq:char} for some $s \in [1/2,1]$ and $q\in [1/2,1]$,
\begin{align}
&\E \Big[\big|\tilde X_\infty\big|^3\mathds{1}_{\{\tilde X_\infty \leq -\zeta\}} \Big] \leq C_3(\gamma,\beta)\delta^{\min\{3(s-1+q/2),0\}    }\delta^{-3/4}. 
\label{eq:3rd_order_mom_QDQED_}
\end{align} 
\end{lemma-sec}

\begin{lemma-sec}[Moment bounds]
\label{lem:mb_dtmc_}
For all $\Lambda, N$, and $\mu$ satisfying $N \geq 1$, $0 < \Lambda < N\mu $, and $0<\mu<1$,
\begin{align}
&\E \Big[\big|\tilde X_\infty \mathds{1}_{\{\tilde X_\infty \geq -\zeta\}}\big| \Big] \leq (\delta^2+1)\frac{1}{\abs{\zeta}}  + \delta,
\label{eq:xplus_}
\end{align}
where, for a set $F$, $\mathds{1}_{F}$ denotes the indicator function of $F$.
\end{lemma-sec}

The proofs of the above two lemmas are given in Sections 1.2 and 1.3 of the online supplement~\cite{FengShi2016a}, 
where we also specify the form of constant $C_3(\gamma, \beta)$. 
Applying the moment bounds~\eqref{eq:3rd_order_mom_QDQED_},~\eqref{eq:xplus_},
and the characterization of $\abs{\zeta}$ in terms of $\delta$ as given by
\[
\abs{\zeta} = \beta R^{q-1/2}
=\beta\delta^{1-2q}, 
\]
we have 
\begin{align}
\E\abs{z(X_\infty)\mu-\Lambda}^3 
&=\delta^{-3} \mu^3 \E[\big|\tilde X_\infty^3\mathds{1}_{\{\tilde X_\infty<-\zeta\}}\big|]
+ \delta^{-3}\mu^3 \abs\zeta^3\Prob(\tilde X_\infty\geq -\zeta)
\nonumber\\
&\leq
\delta^{-3} \mu^3 C_3(\gamma,\beta)\delta^{\min\{3(s-1+q/2),0\}    }\delta^{-3/4}+\delta^{-3}\mu^3 \abs\zeta^2\left[
(\delta^2+2)\frac{1}{\abs{\zeta}}  + \delta\right]\nonumber \\
&=  C_3(\gamma,\beta) \gamma^3 \delta^{\min\{3(s-1+q/2),0\}}\delta^{6s-15/4} +\gamma^3 \beta^2(1+3/\beta) \delta^{6s-2q -2}.
\label{eq:challenge}
\end{align}

Now, we are ready to complete the proof. 
Combining~\eqref{eq:ate5} and~\eqref{eq:ate8}, we get
\begin{align}
\frac{1}{2}\delta^2 \E\left[\epsilon(X_\infty)\right] &
\leq 
4 \delta^3 \frac{1}{\mu}\E\left(\E_{X_\infty}|A_0-D_0|^3\right)+
\gbc\delta^3 \frac{1}{\mu} (1+1/\abs\zeta)
\E\left(\E_{X_\infty}|A_0-D_0|^3\right)
\notag\\
&\leq
 \left(4 + \gbc\right) \delta^3 \frac{1}{\mu} (1+1/\abs\zeta)
\E\left(\E_{X_\infty}|A_0-D_0|^3\right)
. 
\label{eq:part-a}
\end{align}
Plugging~\eqref{eq:challenge0} and \eqref{eq:challenge} in~\eqref{eq:3rd_order_mom},
and using the characterization of $\abs{\zeta}$ in terms of $\delta$, we get from~\eqref{eq:part-a} that 
\begin{align}
\frac{1}{2}\delta^2 \E\left[\epsilon(X_\infty)\right]&\leq 
(  4 + \gbc)
\delta^2 \delta \frac{1}{\mu} (1+1/\abs\zeta)
\E\left(\E_{X_\infty}|A_0-D_0|^3\right)\notag \\
&\leq (  4 + \gbc)\frac{1}{\gamma} (1+1/\beta) \delta^{3-2s}
\Big\{  54\gamma(1+\beta)[1\vee \sqrt{\gamma(1+\beta)}] \delta^{3(s-1)}\notag\\
&+ 9C_3(\gamma,\beta) \gamma^3 \delta^{\min\{3(s-1+q/2),0\}}\delta^{6s-15/4}
+9\gamma^3 \beta^2(1+3/\beta) \delta^{6s-2q -2}\Big \} \notag\\
&\leq (  4 + \gbc)\frac{1}{\gamma} (1+1/\beta)\Big\{  54\gamma(1+\beta)[1\vee \sqrt{\gamma(1+\beta)}] \delta^{s}\notag\\
&+ 9C_3(\gamma,\beta) \gamma^3 \delta^{\min\{3(s-1+q/2),0\}}\delta^{4s-3/4}
+9\gamma^3 \beta^2(1+3/\beta) \delta^{4s-2q +1}\Big \}.
\label{eq:bddQDQED}
\end{align}

Since $4s-2q +1 \geq s,\, 4s-3/4 \geq s$, and
\begin{align*}
3(s-1+q/2)+4s-3/4\geq 3(s-3/4)+4s-3/4=7s-3\geq s,
\end{align*}
\eqref{eq:bddQDQED} completes the proof for part (a) of Lemma~\ref{lem:last_error_dtmc}. 
Finally, we state Lemma~\ref{lem:rvcentralmoments_}, 
the proof of which is given in Section~1 of the online supplement. \hfill\qedsymbol

\begin{lemma-sec}[Random variable absolute central moments]
\label{lem:rvcentralmoments_}
Let $A\sim $Poisson$(\lambda)$, and $D\sim$Binomial$(M,r)$. Then
\begin{align*}
\E[\abs{A-\lambda}^3] &\leq 3 \left[ \lambda \mathds{1}_{\{\lambda <1\}} +\lambda^{3/2}\mathds{1}_{\{\lambda \geq 1\}}\right]
=3\max\left( \lambda ,\, \lambda^{3/2} \right),
\\
\E[\abs{D-Mr}^3]& \leq 3 \left[ Mr \mathds{1}_{\{Mr <1\}} +(Mr)^{3/2}\mathds{1}_{\{Mr \geq 1\}}\right]
=3\max\left( Mr ,\, (Mr)^{3/2} \right).
\end{align*}
\end{lemma-sec}

\subsubsection{Part (b) of Lemma~\ref{lem:last_error_dtmc}}

For part (b), the above issue on bounding $\E_{X_\infty}\left[\abs{A_0-D_0}^3\right]$ no longer exists,
because it is in the order of $\mu N$ 
and can be bounded above by a constant under the condition $s\geq 1$. 
Let $X_\infty = n$. From Lemma~\ref{lem:rvmoments_} stated at the end of this section,  
for any $n \in \mathbb N$, we have 
\begin{equation}
\E_{n}\abs{A_0-D_0}^3 \leq 40 \left[1\vee (N\mu)^2\right]N\mu.
\label{eq:partb-3rd-moment}
\end{equation}
Since the upper bound on the right side does not depend on $n$, 
it is an upper bound for all $n$, and thus, an upper bound for $\E_{X_\infty}\abs{A_0-D_0}^3$.
Thus, 
\begin{align}
\frac{1}{2}\delta^2 \E\left[\epsilon(X_\infty)\right] 
&\leq
4 \delta^3 \frac{1}{\mu}\left\{40 \left[1\vee (N\mu)^2\right]N\mu\right\}
+ \gbc \delta^3 \frac{1}{\mu}\left\{40 \left[1\vee (N\mu)^2\right]N\mu\right\}\left[(1+1/\abs\zeta)\Prob(X_\infty\leq N)\right]
\notag\\
&= 
4 \delta(N\delta^2 )\left\{40 \left[1\vee (N\mu)^2\right]\right\}
\left\{1 +
\frac{\gbc}{4}
\left[(1+1/\abs\zeta)\Prob(X_\infty\leq N)\right]\right\},  
\label{eq:lem-1-part-b-intermed}
\end{align}
which comes from applying~\eqref{eq:partb-3rd-moment} to~\eqref{eq:ate5} and~\eqref{eq:ate8}. 

\blue{However, the inclusion of $q\in [0,1/2)$ poses an additional challenge,
where the term $1/\abs\zeta$ appearing in~\eqref{eq:ate5} and~\eqref{eq:ate8} 
can be larger than a constant in terms of the order of magnitude; in the extreme case of $q = 0$ (NDS regime),
it is in the order of $R^{1/2}$.}
To address this challenge, we explore the property of the idle probability $\Prob(X_\infty \leq N)$. 
Using Lemma~\ref{lem:ip_dtmc} below, 
we have that
\begin{align}
(1+1/\abs\zeta)\Prob(X_\infty\leq N)& \leq 1+
(1/\abs\zeta)\Prob(X_\infty\leq N)\notag\\
&\leq \begin{cases}
1+1/\beta,\quad q\geq 1/2,\\
1 + \frac{2 + \delta}{1-\mu} (1 + \frac{\gamma}{\beta} R^{1 - s - q} ),\quad q< 1/2.
\end{cases}
\label{eq:zeta_ip_}
\end{align} 
Then, part (b) of Lemma \blue{1} follows from~\eqref{eq:lem-1-part-b-intermed}, \eqref{eq:zeta_ip_}, 
and some algebra involving the characterization of $\abs{\zeta}$ as given in Equations~\eqref{eq:abs_zeta}-\eqref{eq:mu_zeta} 
in the online supplement~\cite{FengShi2016a}. \hfill\qedsymbol

\begin{lemma-sec}[Idle probability]
\label{lem:ip_dtmc}
For all $\Lambda, N$, and $\mu$ satisfying $N \geq 1$, $0 < \Lambda < N\mu $, and $0<\mu<1$,
\begin{equation}
    \Prob(X_\infty \leq N) \leq \frac{1}{1-\mu}\left[(2+\delta) (\abs\zeta + \mu\sqrt R)\right].
    \label{eq:ip}
\end{equation} 
\end{lemma-sec}
We outline a coupling argument to prove Lemma~\ref{lem:ip_dtmc}, 
and leave the complete details of the proof in~\ref{sec:pf_ip}.
First, we construct a middle system that bridges our discrete queue and an $M/M/N$ queue.
This middle system is a continuous-time queue with a time-homogeneous Poisson arrival process
and a \emph{two-time-scale} service time component introduced in~\cite{daishi2015c}, 
and we use $Y^M=\{Y^M_k\}$ to denote the \emph{midnight} customer count of this middle system. 
We show that
\begin{equation}
\Prob(X_\infty \leq N) = \Prob(Y^M_\infty \leq N), 
\label{eq:idle-1}
\end{equation} 
using the fact that $Y^M$ has the same dynamics as that of the DTMC $X$ in our discrete queue. 
Then, we construct a $M/M/N$ queue 
such that its \emph{midnight} customer count $Y^C$ 
is always stochastically smaller than $Y^M$, which gives  
\begin{equation}
\Prob(Y^M_\infty \leq N) \leq \Prob(Y^C_\infty \leq N).
\label{eq:idle-2}
\end{equation}
Then, we prove that the stationary distribution of the midnight count $Y^C$ 
is the same as that of the regular customer count process, $X^C=\{X^C(t), t\geq 0\}$, 
in the $M/M/N$ queue, for which we can obtain an upper bound on the idle probability 
$\Prob(X^C(\infty) \leq N)$
\begin{align}
\Prob(X^C(\infty) \leq N) &\leq \left( 2 + \frac{1}{\sqrt{\Lambda/\mu^C}}\right)
\frac{1}{\sqrt{\Lambda/\mu^C}}(N - \Lambda/\mu^C)
\notag\\
&\leq \frac{1}{1 - \mu} ( 2+ \delta)(\abs \zeta + \mu \sqrt R),
\label{eq:ip_c_}
\end{align}
Here, $\mu^C$ denotes the service rate of the $M/M/N$ queue, 
the first inequality follows (3.15) of Lemma 2 in~\cite{BravDaiFeng2015}, 
and the second inequality follows from some algebras detailed in~\ref{sec:pf_ip}.
We eventually have~\eqref{eq:ip} from~\eqref{eq:idle-1} to~\eqref{eq:ip_c_}.

Finally, we state Lemma~\ref{lem:rvmoments_}, the proof of which is given in Section~1 of the online supplement.

\begin{lemma-sec}[Random variable moments]
\label{lem:rvmoments_}

Let $A\sim $Poisson$(\lambda)$, and $D\sim$Binomial$(M,r)$. Then
\begin{align}
\E A^3&=\lambda + 3\lambda^2 +\lambda^3,
\notag\\
\E D^3&=Mr(1-3r+3Mr+2r^2-3Mr^2+M^2r^2).
\notag
\end{align}
\end{lemma-sec}

\subsection{Proof of Lemma~\ref{lem:gradboundsCW}}
\label{sec:proof-gradient}
To establish the bounds on $f'''_h$ as in~\eqref{eq:WCder3},  
we need to first bound $f'_h$ and $f''_h$ since 
\begin{align}
f_h'''(x)&=
-\left(\frac{2b(x)}{a(x)}\right)' f'_h(x) - \left(\frac{2b(x)}{a(x)}\right) f_h''(x) -\frac{2}{a(x)}h'(x)
-\frac{2a'(x)}{a^2(x)} \left[\E h(Y_\infty)
- h(x)\right],
\label{eq:fprimeprimeprime_}
\end{align}
derived from the Poisson equation~\eqref{eq:poisson}, 
where $a'(x)$ is interpreted as the left derivative at the points $x = -1/\delta$ and $x = -\zeta$. 
To bound $\abs{f_h'}$ and $\abs{f_h''}$, we use 
\begin{align}
f'_h(x)&= 
\frac{1}{q(x)} 
\int^x_{-\infty} \frac{2}{a(y)}
\left(\E h(Y_\infty) - h(y) \right)
q(y) \mathrm dy , 
\label{eq:fprime_1_} \\
f'_h(x)&= -
\frac{1}{q(x)} 
\int_x^{\infty} \frac{2}{a(y)}
\left(\E h(Y_\infty) - h(y) \right)
q(y) \mathrm dy , 
\label{eq:fprime_2_}\\
f''_h(x) &= \frac{1}{q(x)}
\int^x_{-\infty} \left\{
-\frac{2}{a(y)}h'(y)
-\frac{2a'(y)}{a^2(y)} \left[\E h(Y_\infty)
- h(y)\right] 
-r'(y) f'_h(y)
\right\} q(y)\mathrm dy ,
\label{eq:fprimeprime_1_}\\
f''_h(x) &= \frac{1}{q(x)}
\int_x^{\infty} \left\{
-\frac{2}{a(y)}h'(y)
-\frac{2a'(y)}{a^2(y)} \left[\E h(Y_\infty)
- h(y)\right] 
-r'(y) f'_h(y)
\right\} q(y)\mathrm dy ,
\label{eq:fprimeprime_2_}
\end{align}
which are also derived from the Poisson equation~\eqref{eq:poisson}, 
with $q(x)$ being defined as 
\begin{align*}
q(x) &\equiv \exp \Big({\int_{0}^{x} \frac{2 b(y)}{a(y)} \mathrm dy} \Big), 
\quad x\in\mathbb{R},
\end{align*}
and $r(x)$ as 
\begin{align*}
 r(x)&\equiv \frac{2b(x)}{a(x)} = 
 \begin{cases}
 \frac{-2x}{1 + \Lambda}, \quad x \leq -1/\delta, \\
 \frac{-2x}{2-\mu + \delta(1-\mu) x + \mu x^2}, \quad x \in [-1/\delta, \abs{\zeta}], \\
 \frac{-2\abs{\zeta}}{2-\mu + \delta(1-\mu) \abs{\zeta} + \mu \zeta^2}, \quad x \geq \abs{\zeta}.
 \end{cases}
 \end{align*} 
Lemma~\ref{lem:ggradbounds-1st-2nd} below states the gradient bounds for $f_h'$ and $f_h''$;
we give an outline of its proof at the end of this section.  
\begin{lemma-sec}
\label{lem:ggradbounds-1st-2nd}
For all $\Lambda>0$, $N\geq 1$, and $\mu\in (0,1)$ satisfying $1 \leq R < N$, 
\begin{align}
\abs{f_h'(x)} &\leq \begin{cases}
\frac{\tilde C_1}{\mu} (1 + 1/\abs{\zeta}),
\quad x\leq -\zeta ,\\
\frac{1}{2} + \frac{1}{\mu\abs{\zeta}}
\left[x+\left(\tilde C + \frac{\delta}{2}\right)
+\left(\tilde C + 1\right) \frac{1}{\abs \zeta}
\right]
,\quad x\geq -\zeta,
\end{cases}
\label{eq:grad_bdd_1_main}\\
\abs{f_h''(x)} & \leq 
\begin{cases}
\frac{\tilde C_2}{\mu}(1 + 1/\abs{\zeta}),
\quad x\leq -\zeta , \\
\frac{1}{\mu \abs{\zeta}},\quad x\geq -\zeta,
\end{cases}
\label{eq:grad_bdd_2_main} 
\end{align}
where $\tilde C_1$, $\tilde C_2$, and $\tilde C$ are constants that only depend on $\mu\delta^{-1}$ when $s < 1/2$
and only depend on $\gamma$ when $s\geq 1/2$. 
\end{lemma-sec}

Now, we are ready to prove Lemma~\ref{lem:gradboundsCW}. From~\eqref{eq:fprimeprimeprime_}
and the following property of function $h$
\begin{align} 
    \abs{\E h( Y_\infty ) - h( y )}
    &\leq \E \abs{ h(Y_\infty) - h(y) }
    \leq \E \abs{Y_\infty - y}
    \leq E\abs{Y_\infty} + \abs{y} , 
    \label{eq:lip}
\end{align} 
we have that
\begin{align*}
\abs{f_h'''(x)}&\leq 
\abs{r'(x) f'_h(x)} + \abs{r(x) f_h''(x)} +\frac{2}{a(x)}
+\frac{2\abs{a'(x)}}{a^2(x)} \E \abs{Y_\infty}
+\frac{2\abs{xa'(x)}}{a^2(x)} . 
\end{align*} 
Applying the bounds on $\abs{f'_h(x)}$ in~\eqref{eq:grad_bdd_1_main}, 
the bound on $\abs{r'(x)}$ in~\eqref{eq:abdd_4_} and the bounds~\eqref{eq:abdd_1_} --~\eqref{eq:abdd_3_} established in Lemma~\ref{lem:a_} below, 
we have that
\begin{align}
\abs{f_h'''(x)}&\leq
2\frac{1}{\mu}
\left[
1 + (\hat C - 1)(1+1/\abs \zeta) \mathds{1}_{\{x\leq -\zeta\}}
\right]
+\abs{r(x) f_h''(x)}
,
\label{eq:abs_fprimeprimeprime_main}
\end{align}
where 
\begin{align*}
\hat C &= 1 + (1 +\tilde{C})\left(\mu\delta^{-1}
\vee 2\right)
+ 3\tilde{C_1} 
.
\end{align*}

For the last term on the right hand side, $\abs{r(x) f_h''(x)}$, via some algebra 
(see intermediate steps in the proof of Lemma~\ref{lem:ggradbounds-1st-2nd} in the online supplement), we can show 
\begin{align}
\abs{r(x) f_h''(x)} &\leq \begin{cases}
2\frac{\hat C}{\mu} (1 + 1/\abs \zeta),
\quad x\leq -\zeta,\\
2\frac{1}{\mu},
\quad x\geq -\zeta .
\end{cases}
\label{eq:abs_fprimeprimeprime_1_main}
\end{align}
Combining~\eqref{eq:abs_fprimeprimeprime_main} and~\eqref{eq:abs_fprimeprimeprime_1_main} 
completes the proof for Lemma~\ref{lem:ggradbounds-1st-2nd}. \hfill\qedsymbol

We give an outline of proving the bounds on $\abs{f_h'}$, $\abs{f_h''}$, 
and leave the complete proof to Section 3 of the online supplement~\cite{FengShi2016a}. 
First, we establish three sets of bounds in Lemmas~\ref{lem:ints_for_grad_bdds_} to~\ref{lem:a_} stated below.  
Then, to bound $\abs{f'_h(x)}$, we apply \eqref{eq:int-1_} to~\eqref{eq:int-4_} in Lemma~\ref{lem:ints_for_grad_bdds_}, 
along with the bound on $\E\abs{Y_\infty}$ in Lemma~\ref{lem:ey_}, to~\eqref{eq:fprime_1_} and~\eqref{eq:fprime_2_}. 
Finally, to bound $\abs{f''_h(x)}$, we apply the bound derived on $\abs{f'_h(x)}$ and Lemma~\ref{lem:a_} 
to deal with the terms $-\frac{2a'(y)y}{a^2(y)}\left[\E h(Y_\infty)- h(y)\right]$ and $-r'(y)f'_h(y)$ in~\eqref{eq:fprimeprime_1_} and~\eqref{eq:fprimeprime_2_}.

\begin{lemma-sec}
\label{lem:ints_for_grad_bdds_}
For all $\Lambda>0$, $N\geq 1$ and $\mu\in(0,1)$ satisfying $1\leq R<N$,
	\begin{align}
	\frac{1}{q(x)} \int^x_{-\infty} \frac{2}{a(y)} q(y) \mathrm dy &\leq 
	\begin{cases}
	\frac{1}{\mu}, \quad x \leq - 1,\\
	\left(1 + 2 e \right)\frac{1}{\mu}, \quad x\in [-1,0],\\
	e^{\zeta^2}  \left(1 + 2e	+  2	\abs{\zeta}	\right)\frac{1}{\mu},
	\quad x\in [0,-\zeta].
	\end{cases} 
	\label{eq:int-1_}\\
	\frac{1}{q(x)} \int_x^\infty \frac{2}{a(y)} q(y) \mathrm dy &\leq
	\begin{cases}
	\left(\frac{1}{\eta}+2\eta e^{\eta^2}\right) \frac{1}{\mu},\quad x \in [0,\eta], \, \eta \leq -\zeta ,\\
	\frac{1}{\mu \abs{\zeta}},\quad x \geq -\zeta.
	\end{cases}
	\label{eq:int-2_}\\
	\frac{1}{q(x)} \int^x_{-\infty} \frac{2\abs{y}}{a(y)}q(y)\mathrm dy &\leq 
	\begin{cases}
	\frac{1}{\mu},\quad x\leq 0,\\
	2 e^{\zeta^2}\frac{1}{\mu},\quad x \in [0,-\zeta].
	\end{cases}
	\label{eq:int-3_}\\
	\frac{1}{q(x)} \int^\infty_x \frac{2 \abs{y}}{a(y)} q(y) \mathrm dy &\leq 
	\begin{cases}
	\frac{3}{2}\frac{1}{\mu} + \frac{\delta}{2}\frac{1}{\mu \abs{\zeta}}
	+ \frac{1}{\mu\zeta^2},
	\quad x \in [0,-\zeta], \\
	\frac{x}{\mu \abs{\zeta}} + \frac{\delta}{2}\frac{1}{\mu \abs{\zeta}}
	+ \frac{1}{\mu\zeta^2}+\frac{1}{2},
	\quad x \geq -\zeta .
	\end{cases}
	\label{eq:int-4_}\\
	\frac{\abs{r(x)}}{q(x)} \int^x_{-\infty}\frac{2}{a(y)}q(y)\mathrm dy &\leq 
	2\frac{1}{\mu},\quad x \leq 0.
	\label{eq:int-5_}\\
	\frac{\abs{r(x)}}{q(x)} \int^\infty_x \frac{2}{a(y)} q(y) \mathrm dy 
	&\leq 2\frac{1}{\mu}, \quad x\geq 0,
	\label{eq:int-6_}
	\end{align}
and when 
$\abs{\zeta} \geq 1$,
\begin{align}
\frac{1}{q(x)} \int_x^\infty \frac{2}{a(y)} q(y) \mathrm dy &\leq 
\begin{cases}
\left(1 + 2 e\right)\frac{1}{\mu},\quad x\in [0,1],\\
\frac{1}{\mu},\quad x \geq 1.
\end{cases}
\label{eq:int-22_}
\end{align}
\end{lemma-sec}

\begin{lemma-sec}
\label{lem:ey_}
For any $s\geq 1/2$ such that~\eqref{eq:char} is satisfied,
\begin{align}
\E \abs{Y_\infty}&\leq 
 \tilde C \left(1 + 1/\abs{\zeta}\right)
 ,
\label{eq:ey__}
\end{align}
where $\tilde{ C } = \tilde C(\gamma)$ is some constant depending only on $\gamma$. 
\end{lemma-sec}

\begin{lemma-sec}
\label{lem:a_}
For all $\Lambda>0$, $N\geq 1$ and $\mu\in(0,1)$ satisfying $1\leq R<N$, 
\begin{align}
a(x) &\geq \mu ,\quad x\in \mathbb{R}, 
\label{eq:abdd_1_}
\\
\frac{\abs{x a'(x)}}{a(x)} &\leq \left((1 - \mu )\vee 2\right) \mathds{1}_{\{x\in (-1/\delta,-\zeta]\}},
\quad x\in\mathbb{R},
\label{eq:abdd_2_}
\\
\E\abs{Y_\infty}\frac{\abs{a'(x)}}{a(x)}&\leq \tilde C \left(\mu\delta^{-1} 
\vee 1\vee 2\right)(1+1/\abs \zeta)
\mathds{1}_{\{x\in (-1/\delta,-\zeta]\}},
\quad x\in \mathbb{R},
\label{eq:abdd_3_}
\\
\abs{r'(x) a(x) } &\leq 
2 \left((2 -\mu )\vee 3\right)\mu 
\mathds{1}_{\{x \leq -\zeta\}} 
,\quad x\in \mathbb {R},
\label{eq:abdd_4_}
\end{align}
where $a'(x)$ and $r'(x)$ are interpreted as the left derivative at $x = 1/\delta$
and $x = -\zeta$.
\end{lemma-sec}

The proofs for Lemmas~\ref{lem:ints_for_grad_bdds_} to~\ref{lem:a_} are lengthy 
and are also deferred to the online supplement~\cite{FengShi2016a}. 
It is worth mentioning that if the diffusion coefficient were a constant such as the one used in Dai and Shi~\cite{daishi2015c}, 
the density $q(x)$ would have a simple Gaussian or exponential form, 
in which case one could directly bound~\eqref{eq:int-1_} -- \eqref{eq:int-22_} by evaluating the integrals.
However, in this paper, $a(x)$ is not only state-dependent, 
but also has a complicated form as given in~\eqref{eq:aform}.  
As a result, $q(x)$ becomes very complicated and the direct method no longer works.
To overcome the difficulty, we build upon a new framework established in~\cite{Brav2017} 
on the gradient bounds for Poisson equation with general $b(x)$ and $a(x)$. 
We illustrate the main idea with~\eqref{eq:int-1_}:
for $x\leq -1$,  since $b(x)/b(y) \geq 1$
\begin{align*}
\frac{1}{q(x)} \int_{-\infty}^{x} \frac{2}{a(y)}q(y) \mathrm dy \leq&\  \frac{1}{q(x)} \int_{-\infty}^{x} \frac{2 b(y)}{a(y)} \frac{1}{b(x)} q(y) \mathrm dy \notag \\
=&\  \frac{1}{q(x)} \frac{1}{b(x)} \int_{-\infty}^{x} \frac{2b(y)}{a(y)} e^{\int_{0}^{y}  \frac{2b(u)}{a(u)} du } \mathrm dy \notag  \\
=&\  \frac{1}{q(x)} \frac{1}{b(x)} \big(q(x) - q(-\infty)\big)  \notag \\
\leq&\ \frac{1}{b(x)} \leq \frac{1}{\mu}. 
\end{align*}
There are many more integrals to deal with, 
which come from expressions of $f_h'$, $f''_h$ and $f'''_h$. 
For each of these integrals, one may need to take different approaches for $x$ in different ranges.

\section{Numerical results}
\label{sec:num}

In this section, we perform extensive numerical experiments to 
(i) demonstrate the quality of using $Y_\infty$ to approximate the stationary distribution of the DTMC $X$;
and (ii) provide support for findings from Theorems~\ref{thm:dtmc_1} and~\ref{thm:dtmc_2}.

\textbf{Parameter settings. }  
\blue{First, we choose a ``baseline'' experimental setting with parameters estimated from the hospital data used in~\cite{daishi2015c}, 
since the discrete queue we study in this paper is motivated from hospital inpatient settings.  
We set $\Lambda=90.95$, corresponding to the daily average number of bed requests for the inpatient wards, 
and $\mu=1/5.3$, which is estimated from the average length-of-stay (5.3 days) from the patients. 
The total number of inpatient beds in the hospital is around 500, 
and we set $N=530$, $504$, and $495$ for the baseline of QE, QED, and NDS regimes, respectively,
which leads to a utilization of $91.0\%$, $95.6\%$, and $97.4\%$, respectively. 
We choose these different values of $N$ such that the system load is roughly light, moderate, and heavy, 
corresponding to the meaning of QD, QED, and NDS regimes, respectively.}

\blue{Then, to demonstrate the convergence of the error bounds under different system load conditions, 
we vary $N$, $R$, and $\mu$ from the baseline setting according to~\eqref{eq:char}. 
That is, we first calculate the value of $\beta$ from 
\begin{align*}
    \text{QD:}&\quad N = R + \beta R, \\
    \text{QED:}&\quad N = R + \beta \sqrt{R}, \\
    \text{NDS:}& \quad N = R + \beta,
\end{align*}
using the baseline parameters, and get $\beta = 0.0995$ for the QD regime, 
$\beta = 0.9994$ for the QED regime, and $\beta = 13$ for the NDS regime. 
We fix $\beta$ for each regime, and increase the values of $R$ by multiplying the baseline $R$ by a factor between 2 and 8. 
Then, we determine the values of $N$ from the above three equations.  
Finally, for each given $s$ (ranging from 0 to 3), we calculate the corresponding $\mu$ from the first equation in~\eqref{eq:char}, 
and $\Lambda$ is determined automatically. 
For scenarios with constant service rate (i.e., $s=0$), $\mu$ is fixed at the baseline value $1/5.3$. 
}

\textbf{Performance measures. }
In the experiments, we evaluate the following performance measures:
the expected queue length $\E(X_\infty-N)^+$,
its adjusted version $\E(X_\infty-R)^+$,
and the expected number of busy servers $\E(X_k \wedge N)$.
One can easily check that the three measures
can be represented by $\sqrt R \,\E[h(\tilde X_\infty)]$
with $h(x) = (x + \zeta)^+$, 
$h(x) = x^+$, and $h(x) = \delta N - (x+\zeta)^-$, respectively,
all satisfying $h(x)\in \textrm{Lip}(1)$.
We compare each of the performance measures
calculated from (i) $\pi$ solved from the exact Markov chain analysis,
and (ii) $\pi$ approximated from using the density function $p$ of $Y_\infty$ given by~\eqref{eq:nu}.
In the interest of space, we mainly report the expected queue length for QED and NDS regimes below. 
We leave supporting tables on the QD regime in~\ref{sec:add_num_res}. 
In the last column of each table, we report the \emph{scaled error},
which equals the absolute error (between the performance calculated from (i) and (ii))
scaled by $1/\sqrt{R}$. 
In other words, the scaled error corresponds to the left side of \eqref{eq:dtmc_1} and~\eqref{eq:dtmc_2}.
Also note that in a given regime, 
we choose the same set of $N$ and $R$ and change $\mu$ using different values of $s$.

\textbf{Constant service rate. }
Tables~\ref{tab:qed_mu_fixed} and~\ref{tab:nds_mu_fixed} demonstrate the results
for the QED and NDS regimes, respectively, with $\mu$ being fixed at $1/5.3$.
Clearly, we can see that the scaled error does not converge to $0$ when we let $R$ and $N$ grow towards infinity.
It is particularly evident for the NDS regime.
Indeed, this observation has motivated us to realize the important role that $\mu$ plays in the convergence of the error bounds,
and led us to identify the important asymptotic regime, i.e., to ensure the convergence of the error bounds,
$\mu$ also needs to converge to 0 as $N, R \rightarrow \infty$.

\begin{table}[!htbp]
\parbox{.45\linewidth}{
\centering
\begin{tabular}{cc|ccc}
 &  	&\multicolumn{2}{c}{$\E(X_\infty-N)^+$} &	\\ 
 N	&	R	&	Stein	&	Exact	&	scaled error	\\ \hline
504	&	482.06	&	4.78	&	4.57	&	0.93\%	\\ \hline
995	&	963.97	&	6.71	&	6.44	&	0.89\%	\\ \hline
1484	&	1446.00	&	8.20	&	7.85	&	0.90\%	\\ \hline
1972	&	1928.12	&	9.45	&	9.05	&	0.91\%	\\ \hline
2946	&	2892.25	&	11.55	&	11.05	&	0.92\%	\\ \hline
3919	&	3856.93	&	13.32	&	12.74	&	0.92\%	\\ 
\end{tabular}
\caption{QED regime with $\beta = 0.9994$ and $\mu$ fixed at $1/5.3$.}
\label{tab:qed_mu_fixed}
}
\hfill
\parbox{.45\linewidth}{
\centering
\begin{tabular}{cc|ccc}
 &  	&\multicolumn{2}{c}{$\E(X_\infty-N)^+$} &	\\ 
 N	&	R	&	Stein	&	Exact	&	scaled error	\\ \hline
495	&	482	&	14.94	&	15.02	&	0.37\%	\\ \hline
977	&	964	&	38.23	&	39.00	&	2.46\%	\\ \hline
1459	&	1446	&	63.82	&	65.49	&	4.41\%	\\ \hline
1941	&	1928	&	90.63	&	93.50	&	6.55\%	\\ \hline
2905	&	2892	&	146.38	&	152.56	&	11.48\%	\\ \hline
3869	&	3856	&	203.91	&	214.63	&	17.27\%	\\ 
\end{tabular}
\caption{NDS regime with $\beta = 13$ and $\mu$ fixed at $1/5.3$.}
\label{tab:nds_mu_fixed}
}
\end{table}

\textbf{Service rate converging to 0. }
Next, instead of fixing $\mu$, we let $\mu$ decrease to $0$ at different rates.
Table~\ref{tab:qed_mu_0point5_1} and~\ref{tab:nds_mu_0point5_1} report the results
for the QED and the NDS regime, respectively.
Table~\ref{tab:qd_mu_0point5_1} in the appendix reports the results in the QD regime.
In each table, the left part corresponds to the case where $\mu\sim 1/R^{1/2}$ (i.e., $s=1/2$),
and the right part corresponds to the case where $\mu \sim 1/R$ (i.e., $s=1$).
We observe that when $\mu \sim 1/R^{1/2}$ and $\mu \sim 1/R$,
the scaled error converges to $0$ as $R, N$ increase under both the QD and QED regime.
In contrast, Table~\ref{tab:nds_mu_0point5_1} clearly shows that in the NDS regime, when $\mu \sim 1/R^{1/2}$,
the scaled error does \emph{not} converge to $0$; in fact, it keeps increasing as $R$ and $N$ become large.
When $\mu \sim 1/R$, the scaled error starts to converge to 0.
These observations are also consistent with our theorems,
that is, not only does $\mu$ need to converge, but it needs to converge at a fast enough rate
to ensure the convergence of the error bounds. 
\blue{We also examine the case with $\mu\sim 1/R^{3/2}$ in line of Theorem~\ref{thm:dtmc_2};
see results in the QED and the NDS regime in Table~\ref{tab:qed_mu_invR1.5} and~\ref{tab:nds_mu_invR1.5} of~\ref{sec:add_num_res}. 
The convergence rates calculated from these numerical results are roughly in the ballpark of 
the corresponding theoretical convergence rates, 
suggesting our error bounds proved in Theorem 1 and 2 are reasonably tight among the tested experiments. 
}

\begin{table}[!htbp]
\begin{center}
\begin{tabular}{cc|cccc|cccc}
 &  	& &\multicolumn{2}{c}{$\E(X_\infty-N)^+$} &	& &\multicolumn{2}{c}{$\E(X_\infty-N)^+$} &\\ 
 N	&	R	&$\mu$&	Stein	&	Exact	&	scaled error	& $\mu$ &	Stein	&	Exact	&	scaled error\\ \hline
504	&	482.06	&	0.189	&	4.78	&	4.57	&	0.93\%	&	0.189	&	4.78	&	4.57	&	0.93\%	\\ \hline
995	&	963.97	&	0.133	&	6.83	&	6.62	&	0.66\%	&	0.094	&	6.90	&	6.76	&	0.48\%	\\ \hline
1484	&	1446.00	&	0.109	&	8.39	&	8.18	&	0.55\%	&	0.063	&	8.50	&	8.38	&	0.33\%	\\ \hline
1972	&	1928.12	&	0.094	&	9.71	&	9.50	&	0.48\%	&	0.047	&	9.84	&	9.73	&	0.25\%	\\ \hline
2946	&	2892.25	&	0.077	&	11.92	&	11.71	&	0.40\%	&	0.031	&	12.07	&	11.98	&	0.17\%	\\ \hline
3919	&	3856.93	&	0.067	&	13.79	&	13.57	&	0.35\%	&	0.024	&	13.95	&	13.87	&	0.13\%	\\ 
\end{tabular}
\caption{QED regime with the parameter $\beta = 0.9994$. On the left, $\mu = 4.1426/R^{1/2}$ and
on the right, $\mu = 90.9542/R$; in both cases the scaled error decreases towards $0$,
but the right side decreases much faster.
}
\label{tab:qed_mu_0point5_1}
\end{center}
\end{table}

\begin{table}[!htbp]
\begin{center}
\begin{tabular}{cc|cccc|cccc}
 &  	& &\multicolumn{2}{c}{$\E(X_\infty-N)^+$} &	& &\multicolumn{2}{c}{$\E(X_\infty-N)^+$} &\\ 
 N	&	R	&$\mu$&	Stein	&	Exact	&	scaled error	& $\mu$ &	Stein	&	Exact	&	scaled error\\ \hline
495	&	482	&	0.189	&	14.94	&	15.02	&	0.37\%	&	0.189	&	15.02	&	14.95	&	0.37\%	\\ \hline
977	&	964	&	0.133	&	39.46	&	39.60	&	0.45\%	&	0.094	&	40.32	&	40.42	&	0.31\%	\\ \hline
1459	&	1446	&	0.109	&	66.79	&	67.01	&	0.56\%	&	0.063	&	68.51	&	68.63	&	0.31\%	\\ \hline
1941	&	1928	&	0.094	&	95.62	&	95.89	&	0.60\%	&	0.047	&	98.12	&	98.25	&	0.29\%	\\ \hline
2905	&	2892	&	0.077	&	155.91	&	156.23	&	0.60\%	&	0.031	&	159.80	&	159.93	&	0.24\%	\\ \hline
3869	&	3856	&	0.067	&	218.35	&	218.72	&	0.58\%	&	0.024	&	223.46	&	223.58	&	0.20\%	\\ 
\end{tabular}
\caption{NDS regime with the parameter $\beta = 13$. On the left, $\mu = 4.1424/R^{1/2}$ and the scaled error increases;
on the right, $\mu = 90.9434/R$ and the scaled error decreases towards $0$.}
\label{tab:nds_mu_0point5_1}
\end{center}
\end{table}

\textbf{Approximation quality. } 
From all the results reported earlier in this section, we can see that the approximation based on $Y_\infty$ works remarkably well 
under a variety of system load conditions and assumptions on $\mu$ when $N$ is large ($>500$). 
In addition, even for small to moderately sized systems, the approximation still works well. 
Tables~\ref{tab:N18} and~\ref{tab:N66} summarize the results 
for $N = 18$ and $N = 66$ respectively, with $\mu=1/5.3$ 
and the system utilization $\rho\equiv R/N$ varying between $88\%$ to $96\%$. 
\blue{We choose to test the two $N$'s
because the number of beds in a single ward is usually between these two values. 
For comparison purpose, we also display the approximation results using $Y^0_\infty$ --   
the r.v. with the diffusion coefficient being a constant $2\mu$ used in~\cite{daishi2015c}. 
Note that in Tables~\ref{tab:N18} and~\ref{tab:N66} we report the relative approximation error,
i.e., absolute error divided by the exact value, 
which differs from the scaled error only in the denominator. 
This is because (i) we are comparing the approximation qualities between using $Y_\infty$ and using $Y^0_\infty$ for given $N$ and $\mu$, 
instead of examining the convergence in the error bound for a sequence of approximations, 
so that the relative error is a more natural measure; 
(ii) the relative error is more prominent to show the difference between the two approximations, 
since the scaled error is usually too small.}

As we can observe from Tables~\ref{tab:N18} and~\ref{tab:N66}, 
$Y_\infty$ significantly improves the approximation quality comparing to $Y^0_\infty$.  
For example, when $N=18$ and $\rho = 90\%$, the relative error is only $0.51\%$ using $Y_\infty$, but is 5.62\% using $Y^0_\infty$. 
\blue{We also perform experiments where $N, R$ are fixed (with $\rho = 90\%$) but $\mu$ varies between $1/2$ and $1/10$, 
corresponding to an average service time of 2 days and 10 days; see Tables~\ref{tab:N18-1} and~\ref{tab:N66-1} in~\ref{sec:add_num_res}. Even though our limit regime requires $\mu\to 0$, we can see that the approximation based on $Y_\infty$ still performs remarkably well when $\mu$ is $1/2$ or $1/3$. 
In addition, when $\mu$ decreases ($1/\mu$ increases), we observe that the approximation values remain the same when using $Y_\infty^0$.   
This is because the density of the constant-diffusion approximation $Y^0_\infty$,}  
\begin{align*} 
p^0(x) &\sim \exp\left(\int^x_0
 -\frac{1}{2} \big( y\wedge ((N-R)/\sqrt{R}) \big) \mathrm{dy}\right),
\quad x\in \mathbb{R} 
\end{align*}
is independent of $\mu$, as both $b(x)$ and $a(x)$ contain a factor $\mu$ and cancel out each other. 
\blue{We report additional results corresponding to moderate ($N=132$) and large ($N=504$) systems in Tables~\ref{tab:N132} and~\ref{tab:N504} of~\ref{sec:add_num_res}. 
We can see that the two approximations have similar performance, where the state-dependent approximation performs slightly better when $N = 132$ 
and has not much difference when $N = 504$. } 

\begin{table}[!htbp]
    \centering
    \begin{tabular}{cc|cc|cc}
$\rho$ &$\E (X_\infty-N)^+$& $\sqrt{R}\E (Y_\infty+\zeta)^+$& relative error&
$\sqrt{R}\E (Y_\infty^0 + \zeta)^+$ & relative error \\
\hline
88\%	&	3.33	&	3.32	&	0.40\%	&	3.47	&	4.10\%	\\ \hline
90\%	&	4.65	&	4.62	&	0.51\%	&	4.91	&	5.62\%	\\ \hline
92\%	&	6.67	&	6.68	&	0.02\%	&	7.18	&	7.51\%	\\ \hline
94\%	&	9.93	&	10.23	&	3.00\%	&	11.10	&	11.81\%	\\ \hline
96\%	&	15.11	&	17.55	&	16.11\%	&	19.19	&	26.99\%	
    \end{tabular}
    \caption{Approximations of the expected queue length using $Y_\infty$ and $Y_\infty^0$ for $N = 18$
    and $\mu = 1/5.3$.}
    \label{tab:N18}
\end{table}
\begin{table}[!htbp]
    \centering
    \begin{tabular}{cc|cc|cc}
$\rho$ &$\E (X_\infty-N)^+$& $\sqrt{R}\E (Y_\infty+\zeta)^+$& relative error&
$\sqrt{R}\E (Y_\infty^0 + \zeta)^+$ & relative error \\
\hline
88\%	&	1.50	&	1.57	&	4.19\%	&	1.53	&	1.46\%	\\ \hline
90\%	&	2.48	&	2.53	&	1.72\%	&	2.58	&	3.77\%	\\ \hline
92\%	&	4.18	&	4.19	&	0.24\%	&	4.42	&	5.69\%	\\ \hline
94\%	&	7.34	&	7.30	&	0.51\%	&	7.87	&	7.26\%	\\ \hline
96\%	&	14.24	&	14.14	&	0.70\%	&	15.45	&	8.54\%	
    \end{tabular}
    \caption{Approximations of the expected queue length using $Y_\infty$ and $Y_\infty^0$ for $N = 66$
    and $\mu = 1/5.3$.}
    \label{tab:N66}
\end{table}

\section{Conclusion}
\label{sec:conclusion}

In this paper, we apply the Stein's method framework
to identify a continuous random variable $Y_\infty$ to approximate
the stationary distribution of the scaled customer count, $\tilde X_\infty$,
in a discrete-time queueing system.
Using this framework, we characterize the convergence rate of the error bounds
between $\tilde X_\infty$ and $Y_\infty$ under different system load conditions.
Different from the continuous-time systems, we identify the important role of $\mu$
in the converge rate of the error bounds.
The numerical results support our theoretical findings.

This work could be extended in several directions. First, under the current queueing setting,
it remains to identify the accurate ``cutoff'' point for $s$ 
that is required to ensure the convergence of the error bounds in each operation regime. 
Second, a major limitation of this paper is the geometric service time assumption. 
Following~\cite{BravDai2015}, one could adapt the approximation developed in this paper 
to queueing systems with discrete phase-type service distributions 
(i.e., replacing the exponential with geometric distributions in the regular phase-type distributions). 
This potentially could lead to efficient algorithms to analyze systems with non-geometric service time distributions.
Third, since our discrete queue is motivated from the hospital inpatient flow management,
a variety of model features that are important in the healthcare context
can be added to the current system, for example, including the day-of-week phenomenon
(which requires a discrete version of time-varying arrival process),
and multiple customer classes to represent patients with different characteristics.

\section*{Acknowledgement}
{
We thank the AE and two anonymous referees for their invaluable feedback to help us improve this paper. 
We also thank Jim Dai at Cornell University and Anton Braverman at Northwestern University 
for their useful comments to this paper.
}

\renewcommand{\baselinestretch}{0.9}
\footnotesize
\def\cprime{$'$} \def\cprime{$'$} \def\cprime{$'$} \def\cprime{$'$}
  \def\cprime{$'$} \def\cprime{$'$} \def\cprime{$'$}

\appendix

\section{Proof of Lemma~\ref{lem:ip_dtmc}}
\label{sec:pf_ip}

\begin{proof}[Proof. ]
As mentioned in the main paper, 
we first construct a ``middle system'', denoted as \emph{System $M$}, to bridge our discrete queue and an $M/M/N$ queue. 
This middle system has $N$ identical servers and a buffer of infinite size.
Customers arrive according to a homogeneous Poisson process with rate $\Lambda$, 
and for each customer, the service time $S^M$ follows a ``two-time-scale'' form~\cite{daishi2015c} 
\begin{align}
S^M = \begin{cases}
\text{LOS}^M + \left(1 - h_{\text{adm}}^M\right),\quad 0<h_{\text{adm}}^M<1,\\
\text{LOS}^M ,\quad h_{\text{adm}}^M = 0.
\end{cases}
\label{eq:sm}
\end{align}
Here, $\text{LOS}^M$ denotes the \emph{number of discrete time epochs} that the customer occupies a server
which takes values on $1,2,\dotsc$, and we assume it follows a geometric distribution with mean $1/\mu$; 
$h_{\text{adm}}^M$ is a number between $0$ (inclusive) and $1$ (exclusive) 
that denotes the instant within a discrete time epoch when the customer is admitted. Mathematically,
\begin{equation}
h_{\text{adm}}^M= {\text{adm}}^M - \lfloor {\text{adm}}^M \rfloor,
\label{eq:hadm}
\end{equation}
where ${\text{adm}}^M$ is the admission time of the customer.

Now, let $X^M = \{X^M(t):t\geq 0\}$ denote the customer count process of system $M$ 
and define its \emph{discrete-time-epoch count} $Y^M = \{Y^M_k:k=0,1,\dotsc\}$ as 
\begin{equation}
Y^M_k = X^M(k) = X^M(k-1) + A^M(k-1,k]-D^M(k-1,k],
\label{eq:ym}
\end{equation}
where $A^M(k-1,k]$ and $D^M(k-1,k]$ denote the total numbers of arrivals and departures 
occurred between time $k-1$ (exclusive) and $k$ (inclusive), respectively.

This process $Y^M$ is referred to as the ``midnight count process'' in~\cite{daishi2015c}, 
and indeed, $Y^M$ has exactly the same dynamics as our DTMC $X$ characterized in \eqref{eq:evol}. 
Because (i) $A^M(k-1,k]$ follows a Poisson distribution with mean $\Lambda$;
and (ii) for $D^M(k-1,k]$, because of the geometric assumption on $\text{LOS}^M$, 
using the coin-toss argument in~\cite{daishi2015c}, 
we can see $D^M(k-1,k]$ follows a binomial distribution with parameters $(Z^M_{k-1},\mu)$, 
with $Z^M_{k-1} = Y^M_{k-1}\wedge N$ denoting the number of busy servers at time $k-1$.
As a result, when the same stability condition $\Lambda < N\mu$ holds, 
the stationary distribution of $Y^M$ uniquely exists and equals $\pi$. 
The corresponding steady-state random variable, $Y^M_\infty$, satisfies
\begin{equation}
\Prob(Y^M_\infty \leq N) = \Prob(X_\infty \leq N).
\label{eq:ymx}
\end{equation}

Next, we consider the M/M/N queue (Erlang-C model) 
where customers arrive according to a homogeneous Poisson process at rate $\Lambda$,
and the service time for each customer, $S^C$, has an exponential distribution with rate $\mu^C = -\log (1-\mu) > 0$.
%

Let $X^C = \{X^C(t):t\geq 0\}$ denote the customer count process of this Erlang-C system, 
and define its discrete-time-epoch count $Y^C = \{Y^C_k:k=0,1,\dotsc\}$ as
\begin{equation}
Y^C_k = X^C(k) = X^C(k-1) + A^C(k-1,k]-D^C(k-1,k],
\label{eq:yc}
\end{equation}
where $A^C(k-1,k]$ and $D^C(k-1,k]$ are the total numbers of arrivals and departures occurred between time $k-1$ (exclusive) and $k$ (inclusive), respectively,
in this Erlang-C system.

Next, we use a \emph{coupling argument} to show that, 
on any given sample path, the discrete-time-epoch count of this Erlang-C system is
always less than that of System $M$, which gives  
\begin{equation}
\Prob(Y^M_\infty \leq N) \leq \Prob(Y^C_\infty \leq N).
\label{eq:idle-3}
\end{equation}
To do so, for a given sample path, we construct a stream of customers, with index $i=1, 2, \dots$ 
to arrive to the Erlang-C system at time $t_1,t_2,\dotsc$. 
These customers are pre-designated with service times $s_1,s_2,\dotsc$, 
sampled from the exponential distribution with rate $\mu^C$. 
Denote these customers' admission and departure times as
$\text{adm}^C_1, \text{adm}^C_2, \dots$, and $\text{dis}^C_1, \text{dis}^C_2, \dots$, respectively.

To couple system $M$ with the Erlang-C system, 
we construct another stream of customers to arrive to System $M$, 
also with index $i=1, 2, \dots$ and arrive at the exact same time $t_1, t_2, \dots$. 
Let their $\text{LOS}^M$ be $\lceil s_1 \rceil  , \lceil s_2 \rceil  \dots$,
and their service times be calculated from~\eqref{eq:sm}. 
According to Section 5.2.3 of~\cite{Baro2013}, 
$\lceil S^C \rceil$ follows a geometric distribution with success probability $ 1 - \exp(-\mu^C)$, which exactly equals $\mu$.
Hence, $\lceil s_1 \rceil  , \lceil s_2 \rceil  \dots$
are indeed generated from a geometric distribution whose mean is $1/\mu$ (and takes values on $1,2,\dotsc$). 
Denote these customers' admission and departure times in System $M$ as
$\text{adm}^M_1, \text{adm}^M_2, \dots$, and $\text{dis}^M_1, \text{dis}^M_2, \dots$, respectively.

The following lemma shows that for each customer, 
the admission and departure times in System $M$ are always earlier than
those in the Erlang-C system. 
\begin{lemma-sec}
\label{lem:ad_ds}
On any given sample path, for each customer $i = 1,2,\dotsc$,
\begin{equation}
\text{adm}_i^C \leq \text{adm}_i^M,\quad \text{dis}_i^C\leq \text{dis}_i^M.
\label{eq:ad_ds}
\end{equation}
\end{lemma-sec}
The proof is given in~\ref{sec:pfLem_ad_ds}. 
We know that for each $i$, customer $i$ arrives to both systems at time $t_{i}$,
and departs at time $\text{dis}_i^C$ from the Erlang-C system. Then, this customer is included in $Y^C$ at
and only at the discrete time epochs $\lceil t_i\rceil, \lceil t_i\rceil+1,\dotsc, \lceil \text{dis}_i^C\rceil -1$.
By Lemma~\ref{lem:ad_ds}, $\text{dis}_i^C\leq \text{dis}_i^M$. This implies that the customer is included in $Y^M$
at least at the discrete time epochs $\lceil t_i\rceil, \lceil t_i\rceil+1,\dotsc, \lceil \text{dis}_i^C\rceil -1$.
Since this observation is true for every customer, we conclude that $Y^C_k\leq Y^M_k$ for all $k = 0,1,\dotsc$,
which proves \eqref{eq:idle-3}.

For the Erlang-C system, one can show that the discrete count process $Y^C$ 
is an irreducible Markov chain which is positive recurrent with a unique stationary distribution, $\pi^Y$,
under the condition that 
\begin{equation}
\Lambda/\mu^C  < N,
\label{eq:cond_stat}
\end{equation}
which always holds when \eqref{eq:stable-condition} is satisfied due to the inequality $\mu <  -\log (1 - \mu)$.
Then, from~\eqref{eq:idle-3} and~\eqref{eq:ymx}, 
\begin{equation}
\Prob(X_\infty \leq N) \leq \Prob(Y^C_\infty \leq N),
\label{eq:ip_cmp}
\end{equation}
where $Y^C_\infty$ denotes the steady-state random variable of $Y^C$.

Finally, we show that in the Erlang-C queue, 
\begin{equation}
\Prob(Y_\infty^C \leq N) = \Prob(X^C(\infty) \leq N).
\label{eq:pasta}
\end{equation} 
To do so, note that $Y^C$ is aperiodic. 
Thus, we have the following relationship according to Proposition 2.9 of~\cite{Sigm2009}
\begin{align}
\pi^Y_n  =\lim_{k \to \infty} \Prob(Y^C_k = n |Y^C_0 = m), \quad m, n \in \mathbb N,
\label{eq:conv_stat_dtmc}
\end{align}
where the type of convergence above is weak convergence. 
Also note that~\eqref{eq:cond_stat} guarantees that $X^C$ is positive recurrent 
with a unique stationary distribution. 
Denote this distribution with $\pi^X$, and the corresponding random variable with $X^C(\infty)$. 
According to Proposition 1.1 of~\cite{Sigm2009a}, we have
\begin{align}
\pi^X_n  =\lim_{t \to \infty} \Prob(X^C(t) = n | X^C(0) = m), \quad m,n \in \mathbb N,
\label{eq:conv_stat_ctmc}
\end{align}
where the type of convergence above is weak convergence. 

Combining~\eqref{eq:conv_stat_dtmc} and~\eqref{eq:conv_stat_ctmc}, 
it is immediately seen that the stationary distributions of $X^C$ and $Y^C$ are the same, 
which implies~\eqref{eq:pasta}.  
Applying Lemma 2 of~\cite{BravDaiFeng2015} to the right side of \eqref{eq:pasta}, we obtain
\begin{align}
\Prob(X^C(\infty) \leq N) &\leq \left( 2 + \frac{1}{\sqrt{\Lambda/\mu^C}}\right)
\frac{1}{\sqrt{\Lambda/\mu^C}}(N - \Lambda/\mu^C)
\notag\\
&\leq \left( 2 + \frac{1}{\sqrt{\Lambda/\nu}}\right)
\frac{1}{\sqrt{\Lambda/\nu}}(N - \Lambda/\nu)
\notag\\
&= 2 \frac{1}{\sqrt{R(1 - \mu)}} \left[N - R(1 - \mu)\right] + \frac{1}{R (1 - \mu)} \left[N - R(1 - \mu)\right]
\notag\\
&= \frac{2}{\sqrt{1 - \mu}} (\abs \zeta + \mu \sqrt R) + \frac{1}{1 - \mu} (\delta \abs \zeta + \mu )
\notag\\
&\leq \frac{1}{1 - \mu} ( 2+ \delta)(\abs \zeta + \mu \sqrt R),
\label{eq:ip_c}
\end{align}
where $\nu = \mu/(1-\mu)$, and the second inequality comes from the fact that for $\mu\in (0,1)$, $-\log (1 - \mu) < \mu/(1-\mu)$.

Combining \eqref{eq:ip_cmp}, \eqref{eq:pasta}, and \eqref{eq:ip_c} establishes Lemma~\ref{lem:ip_dtmc}.
\end{proof}

\subsection{Proof of Lemma~\ref{lem:ad_ds}}
\label{sec:pfLem_ad_ds}

We prove~\eqref{eq:ad_ds} by induction, starting with $i = 1$.

For customer $1$ in both systems, she arrives at time $t_1$ and is admitted immediately. That is,
\begin{equation}
\text{adm}_1^C = \text{adm}_1^M = t_1.
\label{eq:ad1}
\end{equation}

For the departure time of this customer, in the Erlang-C system,
\begin{equation}
\text{dis}_1^C = \text{adm}_1^C + s_1 = t_1 + s_1.
\label{eq:ds1}
\end{equation}
In System $M$, $\text{LOS}^M$ for this customer is $\lceil s_1\rceil$. Since the service time $S^M$
is always greater than or equal to $\text{LOS}^M$, we have
\begin{equation}
\text{dis}_1^M \geq  \text{adm}_1^M + \lceil s_1\rceil \geq  t_1 + s_1.
\label{eq:ds1M}
\end{equation}

Together, \eqref{eq:ad1}, \eqref{eq:ds1}, and \eqref{eq:ds1M} proves~\eqref{eq:ad_ds} for $i = 1$.

Now, assume~\eqref{eq:ad_ds} holds for the first $i$ customers, where $1 \leq i \in \mathbb N$.
Consider the next customer, $i+1$.

We claim that the admission time for this customer satisfies
\begin{equation}
t_{i+1} \leq \text{adm}_{i+1}^C \leq \text{adm}_{i+1}^M.
\label{eq:ad_iplus1}
\end{equation}
Suppose the otherwise. Then at time $\text{adm}_{i+1}^M$, $N$ different customers, $j_1,j_2,\dotsc,j_N\leq i$, are being served by the $N$ servers in the Erlang-C system, whereas at least one of them, $j^*$, has departed from System $M$. This leads to a contradiction, since
we are assuming that~\eqref{eq:ad_ds} holds for all $j \leq i$, and thus in particular for $j^*$.
%
%
%
%

With \eqref{eq:ad_iplus1} in hand, the departure time of customer $i+1$ satisfies
\begin{equation}
\text{dis}_{i+1}^C = \text{adm}_{i+1}^C + s_{i+1} \leq \text{adm}_{i+1}^M + \lceil s_{i+1}\rceil \leq \text{dis}_{i+1}^M.
\label{eq:ds_iplus1}
\end{equation}

Combining \eqref{eq:ad_iplus1} and \eqref{eq:ds_iplus1}, we have shown that \eqref{eq:ad_ds} holds for the first $i+1$ customers. This finishes the induction step and
proves Lemma~\ref{lem:ad_ds}.


\section{Computational time to evaluate $\pi$ using the Markov chain analysis}
\label{sec:add_num}

All the experiments are implemented in Matlab and run on a Linux $64$-bit cluster hosted either by a Dell R620 server 
or by a R720 server (indicated by an ``*'') with $1$ processor.
The virtual size supplied by the cluster is $8192$ MiB for each experiment with $N = 4242$ or $977$.

The computational time of these experiments is shown in Table~\ref{tab:comp_time_1}.

\begin{table}[!htbp]
\small{
\parbox{.45\linewidth}{
\centering
\begin{tabular}{c|c}
$\mu$	&	elapsed time (day)\\ \hline
0.189	&	1.25\\ \hline
0.067	&	0.84*\\ \hline
0.040	&	0.84*\\ \hline 
0.024	&	0.87\\ \hline 
0.008	&	1.13
\end{tabular}
\caption*{\small{$N= 4242$, $\rho = 90\%$}}
}
\hfill
\parbox{.45\linewidth}{
\centering
\begin{tabular}{c|c}
$\mu$	&	elapsed time (hour)\\ \hline
0.189	&	1.34\\ \hline
0.133	&	0.96*\\ \hline
0.112	&	0.93*\\ \hline 
0.094	&	1.18*\\ \hline
0.067	&	0.97*
\end{tabular}
\caption*{\small{$N = 977$, $\rho = 98\%$}}
}
\caption{\small{Computational time of $\pi$ for large systems with moderately high utilization (left) and for moderately large systems with high utilization (right).}}
\label{tab:comp_time_1}
}
\end{table}

\section{
Additional numerical results
}
\label{sec:add_num_res}


This section includes numerical results for the QD regime, for the QED and the NDS regime with $\mu\sim 1/R^{3/2}$, 
and for small to moderately sized systems. 

In QD regime, the expected queue length is very close to $0$ because of the light system load. Hence we report an adjusted version of the queue length, $\E(X_\infty-R)^+$, in Table~\ref{tab:qd_mu_fixed} and~\ref{tab:qd_mu_0point5_1}.

\begin{table}[ht]
\begin{center}
\small{
\begin{tabular}{cc|ccc}
 &  	&\multicolumn{2}{c}{$\E(X_\infty-R)^+$} &	\\ 
 N	&	R	&	Stein	&	Exact	&	scaled error	\\ \hline
530	&	482.06	&	8.80	&	8.91	&	0.48\%	\\ \hline
1061	&	965.02	&	12.31	&	12.40	&	0.87\%	\\ \hline
1591	&	1447.08	&	14.80	&	15.18	&	0.98\%	\\ \hline
2121	&	1929.14	&	17.08	&	17.52	&	1.01\%	\\ \hline
3182	&	2894.16	&	20.91	&	21.50	&	1.10\%	\\ \hline
4242	&	3858.28	&	24.15	&	24.78	&	1.02\%	\\ 
\end{tabular}
\caption{\small{QD regime with the parameter $\beta = 0.0995$ and $\mu$ fixed at $1/5.3$.}}
\label{tab:qd_mu_fixed}
}
\end{center}
\end{table}

\begin{table}[t]
\begin{center}
\small{
\begin{tabular}{cc|cccc|cccc}
 &  	& &\multicolumn{2}{c}{$\E(X_\infty-R)^+$} &	& &\multicolumn{2}{c}{$\E(X_\infty-R)^+$} &\\ 
 N	&	R	&$\mu$&	Stein	&	Exact	&	scaled error	& $\mu$ &	Stein	&	Exact	&	scaled error\\ \hline
530	&	482.06	&	0.189	&	8.80	&	8.91	&	0.48\%	&	0.189	&	8.80	&	8.91	&	0.48\%	\\ \hline
1061	&	965.02	&	0.133	&	12.21	&	12.40	&	0.61\%	&	0.094	&	12.27	&	12.40	&	0.43\%	\\ \hline
1591	&	1447.08	&	0.109	&	14.97	&	15.18	&	0.55\%	&	0.063	&	15.06	&	15.18	&	0.32\%	\\ \hline
2121	&	1929.14	&	0.094	&	17.31	&	17.52	&	0.49\%	&	0.047	&	17.42	&	17.52	&	0.24\%	\\ \hline
3182	&	2894.16	&	0.077	&	21.25	&	21.46	&	0.40\%	&	0.031	&	21.38	&	21.46	&	0.16\%	\\ \hline
4242	&	3858.28	&	0.067	&	24.57	&	24.78	&	0.34\%	&	0.024	&	24.71	&	24.78	&	0.12\%	\\ 
\end{tabular}
}
\caption{\small{QD regime with the parameter $\beta = 0.0995$. On the left, $\mu = 4.1426/R^{1/2}$ and
on the right, $\mu = 90.9542/R$; in both cases the scaled error decreases towards $0$
at a very similar rate.} 
}
\label{tab:qd_mu_0point5_1}
\end{center}
\end{table}

\FloatBarrier

Also we include the results for the QED and the NDS regime with $\mu \sim 1/R^{3/2}$. In both tables the scaled errors converge to $0$, as we expect from Theorem~\ref{thm:dtmc_2}. Moreover, one can see that the convergence rates are slower than $1/R$, which is the convergence rate of the approximation errors proven in Theorem 3.1 of~\cite{Brav2017}. This illustrates the different nature of our discrete queue from that of the continuous-time queueing system studied there.

\begin{table}[!htbp]
\parbox{.45\linewidth}{
\centering
\begin{tabular}{cc|ccc}
 &  	&\multicolumn{2}{c}{$\E(X_\infty-N)^+$} &	\\ 
 N	&	R	&	Stein	&	Exact	&	scaled error	\\ \hline
504	&	482.06	&	4.78	&	4.57	&	0.93\%	\\ \hline
995	&	963.97	&	6.96	&	6.85	&	0.34\%	\\ \hline
1484	&	1446.00	&	8.56	&	8.49	&	0.19\%	\\ \hline
1972	&	1928.12	&	9.90	&	9.85	&	0.13\%	\\ \hline
2946	&	2892.25	&	12.13	&	12.09	&	0.07\%	\\ \hline
3919	&	3856.93	&	14.00	&	13.98	&	0.05\%	
\end{tabular}
\caption{
QED regime with $\beta = 0.9994$ and $\mu
= 1996.9729/R^{3/2}
$ 
.
}
\label{tab:qed_mu_invR1.5}
}
\hfill
\parbox{.45\linewidth}{
\centering
\begin{tabular}{cc|ccc}
 &  	&\multicolumn{2}{c}{$\E(X_\infty-N)^+$} &	\\ 
 N	&	R	&	Stein	&	Exact	&	scaled error	\\ \hline
495	&	482	&	14.94	&	15.02	&	0.37\%	\\ \hline
977	&	964	&	40.93	&	41.00	&	0.21\%	\\ \hline
1459	&	1446	&	69.50	&	69.56	&	0.18\%	\\ \hline
1941	&	1928	&	99.37	&	99.43	&	0.14\%	\\ \hline
2905	&	2892	&	161.39	&	161.44	&	0.10\%	\\ \hline
3869	&	3856	&	225.27	&	225.31	&	0.07\%	
\end{tabular}
\caption{NDS regime with $\beta = 13$ and $\mu
= 1996.6166/R^{3/2}
$ 
.
}
\label{tab:nds_mu_invR1.5}
}
\end{table}

Next we report the numerical results for moderate to large systems with $N = 132$ and $N = 504$, with the system utilization varying between $88\%$ to $96\%$.


\begin{table}[!h]
    \centering
    \begin{tabular}{cc|cc|cc}
$\rho$ & $\E (X_\infty-R)^+$ & $\sqrt{R}\E (Y_\infty)^+$ & relative error&
$\sqrt{R}\E (Y_\infty^0)^+$ & relative error \\
\hline
88\%	&	4.87	&	4.86	& 0.13\% &	4.87	&	0.02\% \\ \hline
90\%	&	5.48	&	5.48	&	0.03\%	&	5.52	&	0.72\% \\ \hline
92\%	&	6.65	&	6.63	&	0.20\%	&	6.78	&	2.07\% \\ \hline
94\%	&	9.04	&	8.99	&	0.47\%	&	9.41	&	4.18\% \\ \hline
96\%	&	14.80	&	14.70	&	0.72\%	&	15.79	&	6.67\%
    \end{tabular}
    \caption{Approximations of the expected queue length (adjusted) using $Y_\infty$ and $Y_\infty^0$ for $N = 132$
    and $\mu = 1/5.3$.}
    \label{tab:N132}
\end{table}
\begin{table}[!h]
    \centering
    \begin{tabular}{cc|cc|cc}
$\rho$ &$\E (X_\infty-R)^+$& $\sqrt{R}\E (Y_\infty)^+$& relative error&
$\sqrt{R}\E (Y_\infty^0)^+$ & relative error \\
\hline
88\%	&	8.42	&	8.25	&	2.02\%	&	8.42	&	0.04\% \\ \hline 
90\%	&	8.58	&	8.46	&	1.45\%	&	8.58	&	0.05\% \\ \hline 
92\%	&	8.95	&	8.89	&	0.66\%	&	8.96	&	0.06\% \\ \hline 
94\%	&	10.04	&	10.04	&	0.03\%	&	10.12	&	0.84\% \\ \hline 
96\%	&	13.69	&	13.67	&	0.13\%	&	14.14	&	3.30\%
    \end{tabular}
    \caption{Approximations of the expected queue length (adjusted) using $Y_\infty$ and $Y_\infty^0$ for $N = 504$
    and $\mu = 1/5.3$.}
    \label{tab:N504}
\end{table}

Lastly we include some additional numerical results for small systems with $N = 18$ and $N = 66$, where for each $N$, we fix $R$ such that utilization $\rho$ stays fixed at $90\%$, and let $\mu$ vary. 


\begin{table}[h!]
    \centering
    \begin{tabular}{cc|cc|cc}
$1/\mu$ &$\E (X_\infty-N)^+$& $\sqrt{R}\E (Y_\infty+\zeta)^+$& relative error&
$\sqrt{R}\E (Y_\infty^0 + \zeta)^+$ & relative error \\
\hline
 2	&	3.88	&	3.76	&	3.22\%	&	4.91	&	26.41\% \\ \hline 
 3	&	4.29	&	4.22	&	1.68\%	&	4.91	&	14.34\% \\ \hline 
 4	&	4.50	&	4.45	&	1.06\%	&	4.91	&	9.08\% \\ \hline 
 5	&	4.63	&	4.59	&	0.74\%	&	4.91	&	6.13\% \\ \hline 
 6	&	4.71	&	4.68	&	0.54\%	&	4.91	&	4.25\% \\ \hline 
 7	&	4.77	&	4.75	&	0.40\%	&	4.91	&	2.94\% \\ \hline 
 8	&	4.81	&	4.80	&	0.31\%	&	4.91	&	1.97\% \\ \hline 
 9	&	4.85	&	4.84	&	0.23\%	&	4.91	&	1.24\% \\ \hline 
  10	&	4.88	&	4.87	&	0.18\%	&	4.91	&	0.65\%
    \end{tabular}
    \caption{Approximations of the expected queue length using $Y_\infty$ and $Y_\infty^0$ for $N = 18$
    and $R = 16.20$.}
    \label{tab:N18-1}
\end{table}

\begin{table}[h!]
    \centering
    \begin{tabular}{cc|cc|cc}
$1/\mu$ &$\E (X_\infty-N)^+$& $\sqrt{R}\E (Y_\infty+\zeta)^+$& relative error&
$\sqrt{R}\E (Y_\infty^0 + \zeta)^+$ & relative error \\
\hline
 2	&	2.09	&	2.16	&	3.26\%	&	2.58	&	23.32\% \\ \hline 
 3	&	2.30	&	2.36	&	2.61\%	&	2.58	&	12.18\% \\ \hline 
 4	&	2.41	&	2.46	&	2.12\%	&	2.58	&	7.21\% \\ \hline 
 5	&	2.47	&	2.51	&	1.79\%	&	2.58	&	4.39\% \\ \hline 
 6	&	2.51	&	2.55	&	1.56\%	&	2.58	&	2.58\% \\ \hline 
 7	&	2.54	&	2.58	&	1.38\%	&	2.58	&	1.32\% \\ \hline 
 8	&	2.57	&	2.60	&	1.25\%	&	2.58	&	0.39\% \\ \hline 
 9	&	2.59	&	2.62	&	1.14\%	&	2.58	&	0.33\% \\ \hline 
  10	&	2.60	&	2.63	&	1.05\%	&	2.58	&	0.89\% 
    \end{tabular}
    \caption{Approximations of the expected queue length using $Y_\infty$ and $Y_\infty^0$ for $N = 66$
    and $R = 59.40$.}
    \label{tab:N66-1}
\end{table}

\renewcommand{\thesection}{\arabic{section}}
\setcounter{section}{0}

\begin{frontmatter}

\title{Online Supplement of \\
``Steady-state Diffusion Approximations for Discrete-time Queue in Hospital Inpatient Flow Management''} 



%



\end{frontmatter}

\linenumbers

\normalsize

This document serves as the online supplement for Feng and Shi~\cite{FengShi2016}, 
which we refer to as the ``main paper.'' 
In the main paper,  
we analyze a $GI/Geo/N$ discrete-time queue (or \emph{discrete queue} in short), 
and use the Stein's method framework to develop steady-state diffusion approximations
for the customer count process, with a focus on the Poisson arrival case. 
We establish the error bounds of the approximations in Theorems 1 and 2 there.   
The proof of these two theorems rely on Lemmas 1 and 2. 
Section~3.1 of the main paper proves Lemma 1, 
which further depends on several additional lemmas proved in Section~\ref{sec:addLem} of this document. 
We also give the complete proof of Lemma 2 in Section~\ref{sec:general-gradient-bounds} of this document.  
Section~\ref{sec:general_arr} of this document extends the results in the main paper 
by considering a general arrival distribution, 
where we develop analogous steady-state approximations 
and establish the corresponding error bounds.  

\normalsize 



\section{Proof of additional lemmas}
\label{sec:addLem}

\subsection{Moments of random variables}

\begin{customthm}{3.4 (Feng \& Shi '17)}[Random variable moments]
\label{lem:rvmoments}

Let $A\sim $Poisson$(\lambda)$, and $D\sim$Binomial$(M,r)$. Then
\begin{align}
\E A &= \lambda,\,\E A^2 = \lambda +\lambda^2,\,\E A^3=\lambda + 3\lambda^2 +\lambda^3,\,
\E A^4=\lambda + 7\lambda^2 +6\lambda^3+\lambda^4.
\label{eq:poissonmoments}\\
\E D&=Mr,\,\E D^2=Mr(1-r+Mr),\,
\E D^3=Mr(1-3r+3Mr+2r^2-3Mr^2+M^2r^2),\notag\\
\E D^4&=Mr(1-7r+7Mr+12r^2-18Mr^2+6M^2r^2-6r^3+11Mr^3-6M^2r^3+M^3r^3).
\label{eq:binomialmoments}
\end{align}
\end{customthm}
\begin{proof}
%
See \cite{WeisPois} and \cite{WeisBino}.
\end{proof}
\begin{remark}

Note that for $r \in [0,1]$, 
\begin{align*}
1 - 3r + 2r^2 &\leq 1,\\
1-7r+12r^2-6r^3 &\leq 1-7r+12r^2\leq 6,\\
7Mr - 18 Mr^2 + 11 Mr^3 &= 7Mr - Mr^2(18-11r)\leq 7Mr.
\end{align*}
Hence, \eqref{eq:binomialmoments} implies
\begin{equation}
\E D^3 \leq 5 \max\left\{Mr,\,(Mr)^3\right\},\quad \E D^4 \leq 20 \max \left\{Mr,\,(Mr)^4\right\}.
\label{eq:binomial3rdmom}
\end{equation}
\end{remark}
%
%

\begin{remark}
Applying the $c_r$-inequality and Lemma~\ref{lem:rvmoments}, we have the following inequalities
\begin{align}
\E |A - D|^3 &\leq 4 \E A^3 + 4\E D^3
\leq 20 \max(\lambda, \lambda^3) + 20 \max(Mr, (Mr)^3).
\label{eq:e3}
\end{align}
and
\begin{align}
\E_n |A_0-D_0| &\leq \E A_0 + \E_n D_0 \leq 2 N\mu,
\label{eq:en1}
\\
\E_n (A_0-D_0)^2 &\leq 2\E A_0^2+2 \E_n D_0^2\leq 8\left[1\vee (N\mu)\right]N\mu ,
\label{eq:en2}
\\
\E_n |A_0-D_0|^3&\leq 4 \E A_0^3 + 4 \E_n D_0^3\leq 40\left[1\vee (N\mu)^2\right]N\mu,
\label{eq:en3}
\\
\E_n(A_0-D_0)^4&\leq 8\E A_0^4 + 8 \E_n D_0^4\leq 280 \left[1\vee (N\mu)^3\right]N\mu.
\label{eq:en4}
\end{align}
where $\E_n$ is the expectation under $\Prob_n$, the conditional probability distribution given that the starting customer count equals $n$.
\end{remark}

\begin{customthm}{3.1 (Feng \& Shi '17)}[Random variable absolute central moments]
\label{lem:rvcentralmoments}
Let $A\sim $Poisson$(\lambda)$, and $D\sim$Binomial$(M,r)$. Then
\begin{align}
\E[\abs{A-\lambda}^3] &\leq 3 \left[ \lambda \mathds{1}_{\{\lambda <1\}} +\lambda^{3/2}\mathds{1}_{\{\lambda \geq 1\}}\right]
=3\max\left( \lambda ,\, \lambda^{3/2} \right),
\label{eq:poissoncm}\\
\E[\abs{D-Mr}^3]& \leq 3 \left[ Mr \mathds{1}_{\{Mr <1\}} +(Mr)^{3/2}\mathds{1}_{\{Mr \geq 1\}}\right]
=3\max\left( Mr ,\, (Mr)^{3/2} \right).
\label{eq:binomialcm}
\end{align}
\end{customthm}

\begin{proof}
From \cite{WeisPois}, the central moments of $A$ are
\begin{equation*}
\E \left[(A-\lambda)^4\right] = \lambda + 3\lambda^2,\quad\E \left[(A-\lambda)^3\right] = \lambda.
\end{equation*}
When $\lambda \geq 1$, $\E \left[(A-\lambda)^4\right] \leq 4\lambda^2$, and Jensen's inequality implies that
\begin{equation*}
\E \left[\abs{A-\lambda}^3\right]\leq \left\{\E \left[(A-\lambda)^4\right]\right\}^{3/4}\leq 3\lambda^{3/2}.
\end{equation*}
When $0<\lambda<1$,
\begin{equation*}
\E \left[\abs{A-\lambda}^3\right] = \E\left [(A-\lambda)^3\right] + 2 \lambda^3 \Prob(A = 0) = \lambda + 2 \lambda^3 e^{-\lambda} \leq 3\lambda.
\end{equation*}

From \cite{WeisBino},
the central moments of $D$ are
\begin{equation*}
\E \left[(D - Mr)^4\right] = 3M^2 r^2 (1-r)^2 + M r(1-r)[1-6r(1-r)],\quad
\E \left[(D-Mr)^3\right] = Mr (1-r)(1-2r).
\end{equation*}
When $Mr \geq 1$,
\begin{equation*}
\E \left[(D - Mr)^4\right]\leq 3M^2 r^2 + Mr \leq 4(Mr)^2.
\end{equation*}
Jensen's inequality implies
\begin{equation*}
\E \left[\abs{D-Mr}^3\right]\leq \left\{\E \left[(D-Mr)^4\right]\right\}^{3/4} \leq 3(Mr)^{3/2}.
\end{equation*}
When $0<Mr<1$,
\begin{equation*}
\E \left[\abs{D-Mr}^3\right] = \E \left[(D-Mr)^3\right] + 2 (Mr)^3 \Prob(D = 0) \leq Mr + 2 (Mr)^3\leq 3Mr.
\end{equation*}
\end{proof}


\subsection{
Moment bounds of $\tilde{X}_\infty$
}
\label{subsec:mb}
\begin{customthm}{3.3 (Feng \& Shi '17)}[Moment bounds]
\label{lem:mb_dtmc}
For all $\Lambda, N$, and $\mu$ satisfying $N \geq 1$, $0 < \Lambda < N\mu $, and $0<\mu<1$,
\begin{align}
&\E \Big[(\tilde X_\infty)^2 \mathds{1}_{\{\tilde X_\infty \leq -\zeta\}}\Big] \leq \frac{4}{3} + \frac{8}{3}\delta^2,
\label{eq:xsquaredelta}\\
&\E \Big[ \big|\tilde X_\infty \mathds{1}_{\{\tilde X_\infty \leq -\zeta\}}\big| \Big] \leq \sqrt{ \frac{4}{3} + \frac{8}{3}\delta^2},
\label{eq:xminusdelta}\\
&\E \Big[ \big|\tilde X_\infty \mathds{1}_{\{\tilde X_\infty \leq -\zeta\}}\big| \Big] \leq 2\abs{\zeta}
\label{eq:xminuszeta}\\
&\E \Big[\big|\tilde X_\infty \mathds{1}_{\{\tilde X_\infty \geq -\zeta\}}\big| \Big] \leq (\delta^2+1)\frac{1}{\abs{\zeta}}  + \delta,
\label{eq:xplus}\\
&\E \Big[\big|(\tilde X_\infty+\zeta)\mathds{1}_{\{\tilde X_\infty \leq -\zeta\}}\big| \Big] =\abs\zeta ,
\label{eq:xpluszeta}
\end{align}
where, for a set $F$, $\mathds{1}_{F}$ denotes the indicator function of $F$.
\end{customthm}
\begin{proof}[Proof of Lemma~\ref{lem:mb_dtmc}]
The proof is similar to that in Appendix \blue{E.3.1.} of Dai and Shi~\cite{daishi2015cc}, 
and thus is omitted here. 
\end{proof}

\subsection{Bounding partial third moment (3.11) (Feng \& Shi '17)}  
\label{sec:pfLem3rd_order_mom_QDQED}

Consider a function $V(x) = x^4 + a_1 x^3 + a_2 x^2$, where $a_1$ and $a_2$ are two constants that will be determined later.
Recalling the definition of $G_{\tilde X}$ in the main paper, there is
\begin{align*}
G_{\tilde X}V(x)& = 4x^3 \delta \E_n (A_0-D_0) +\left\{ 6\delta^2 \E_n [(A_0-D_0)^2]+3a_1 \delta \E_n(A_0-D_0) \right\} x^2\\
&+\left\{ 4\delta^3 \E_n[(A_0-D_0)^3] +3a_1\delta^2 \E_n [(A_0-D_0)^2]+2a_2 \delta \E_n(A_0-D_0) \right\} x\\
&+\delta^4 \E_n [(A_0-D_0)^4]+a_1\delta^3 \E_n[(A_0-D_0)^3] +a_2\delta^2 \E_n [(A_0-D_0)^2],
\end{align*}
where $n\in\mathbb N$ is such that $x=\delta (n-x_\infty)$.

Now we determine $a_1$ and $a_2$ by taking them to satisfy
\begin{align}
6\delta^2 \E_N [(A_0-D_0)^2] &+3a_1 \delta \E_N(A_0-D_0) = 0,
\notag\\
4\delta^3 \E_N[(A_0-D_0)^3] &+3a_1\delta^2 \E_N [(A_0-D_0)^2]+2a_2 \delta \E_N(A_0-D_0)= 0.
\label{eq:a12}
\end{align}
Then,
\begin{align}
a_1 &= 2\delta \E_N [(A_0-D_0)^2]/\E_N(D_0-A_0)\notag\\
&= 2\delta \left[1-\mu + (2 - \mu) \frac{\Lambda}{N\mu - \Lambda} + (N\mu - \Lambda)\right]\notag\\
&= 2 \left[(1 - \mu)\delta + \mu\abs\zeta + (2-\mu)\frac{1}{\abs \zeta}  \right],
\label{eq:a1formula}
\end{align}
and
\begin{align}
a_2 &= \left\{4\delta^2 \E_N[(A_0-D_0)^3] +3a_1\delta \E_N [(A_0-D_0)^2]\right\} /2\E_N(D_0-A_0)\notag\\
&\leq 2 \delta^2 40 \left[1 \vee (N\mu)^2\right] \frac{\delta N}{\abs \zeta} + \frac{3}{4}a_1^2\notag\\
&=80 \left[1 \vee (N\mu)^2\right] \frac{N \delta^3}{\abs\zeta}+ \frac{3}{4}a_1^2.
\label{eq:a2formula}
\end{align}
where the second line uses the inequality~\eqref{eq:en3}.

With the $a_1$ and $a_2$ chosen as above,
\begin{align}
G_{\tilde X}V(x)& = 4x^3 b(x) +\left\{ 6\delta^2 \E_n [(A_0-D_0)^2]+3a_1 \delta \E_n(A_0-D_0) \right\} x^2\mathds{1}_{\{x<-\zeta\}}\notag\\
&+\left\{ 4\delta^3 \E_n[(A_0-D_0)^3] +3a_1\delta^2 \E_n [(A_0-D_0)^2]+2a_2 \delta \E_n(A_0-D_0) \right\} x\mathds{1}_{\{x<-\zeta\}}\notag\\
&+\delta^4 \E_n [(A_0-D_0)^4]+a_1\delta^3 \E_n[(A_0-D_0)^3] +a_2\delta^2 \E_n [(A_0-D_0)^2],
\label{eq:gxVx}
\end{align}
where
\begin{align}
 4x^3 b(x) &= - 4 \mu x^4 \mathds{1}_{\{x \leq -\zeta\}} -4 \mu \abs\zeta x^3 \mathds{1}_{\{x > -\zeta\}}.
 \label{eq:4x3bx}
\end{align}

Now with a proof similar to that in Appendix E.4.2 of \cite{daishi2015cc},  
we have the following basic adjoint relation (BAR)
\begin{align}
\E \left[G_{\tilde X}V(\tilde X_\infty)\right] = 0.
\label{eq:bar}
\end{align}

Taking expectation with respect to $\tilde X_\infty$ on both sides of \eqref{eq:gxVx}, we obtain
\begin{align}
4 \mu \E[\tilde X_\infty^4 \mathds{1}_{\{\tilde X_\infty \leq -\zeta\}}] &\leq 2 N\mu \left\{3  a_1\delta\left(\frac{4}{3} + \frac{8}{3}\delta^2\right)+ 2a_2 \delta \left[\sqrt{\frac{4}{3} + \frac{8}{3}\delta^2}\wedge 2\abs\zeta\right] \right\}
\notag\\
&+  8\left[1\vee (N\mu)\right]N\mu
\left\{6\delta^2 \left(\frac{4}{3} + \frac{8}{3}\delta^2\right) + 3 a_1 \delta^2 \left[\sqrt{\frac{4}{3} + \frac{8}{3}\delta^2}\wedge 2\abs\zeta\right]
+ a_2 \delta^2\right\}
\notag\\
&+ 40\left[1\vee (N\mu)^2\right]N\mu \left\{4\delta^3   \left[\sqrt{\frac{4}{3} + \frac{8}{3}\delta^2}\wedge 2\abs\zeta\right]
+  a_1\delta^3\right\}\notag\\
&+ 280 \left[1\vee (N\mu)^3\right]N\mu \delta^4,
\label{eq:lyapunovQDQED}
\end{align}
where we use the moment bounds~\eqref{eq:xsquaredelta},~\eqref{eq:xminusdelta},~\eqref{eq:xminuszeta} 
in this supplement and
the inequalities in Remark \blue{2} of Lemma~\ref{lem:rvmoments}. 

Applying the assumption $R\geq 1$, or $\delta \leq 1$, \eqref{eq:lyapunovQDQED} implies that
\begin{align}
\E[\tilde X_\infty^4 \mathds{1}_{\{\tilde X_\infty \leq -\zeta\}}] &\leq
\left\{48 \left[1\vee (N\mu)\right] + 80 \left[1\vee (N\mu)^2\right]\delta +70\left[1\vee (N\mu)^3\right]\delta^2\right\}N\delta^2 \notag\\
&+
a_1\left\{6\delta^{-1} + 12\left[1\vee (N\mu)\right](1\wedge \abs \zeta) + 10 \left[1\vee (N\mu)^2\right]\delta\right\}N\delta^2 \notag\\
&+ a_2 \left\{2(1\wedge \abs \zeta)\delta^{-1} + 2\left[1\vee (N\mu)\right]\right\} N\delta^2.
\label{eq:lyapunovQDQED1}
\end{align}

Now, recall the characterization of $\mu$ and $N$ in~\eqref{eq:char} of the main paper.
Under the settings of the theorems,
either $s \in [1/2,1]$, $q\in [1/2,1]$ or $s\geq 1$, $q\in [0,1]$. In both cases, we have
\begin{align}
\abs\zeta &= \beta R^{q-1/2} = \beta \delta^{1-2q},
\label{eq:abs_zeta}\\
N\delta^2&=(N-R+R)/R=\beta R^{q-1} + 1\leq \beta +1,
\label{eq:Ndeltasq}\\
N\mu &= \gamma (N\delta^2) \delta^{2s-2}\leq \gamma (\beta+1)\delta^{2s-2}.
\label{eq:Nmu}\\
\mu/\abs\zeta&=\frac{\gamma \delta^{2s}}{\beta \delta^{1-2q}}
=\frac{\gamma}{\beta}\delta^{2s+2q-1}
\leq \frac{\gamma}{\beta}\delta.
\label{eq:mu_zeta}
\end{align}

Substituting \eqref{eq:abs_zeta} - \eqref{eq:Nmu} into \eqref{eq:a1formula} and \eqref{eq:a2formula} implies that
\begin{align}
a_1 &\leq 2 \left(\delta + \gamma \delta^{2s}\beta \delta^{1-2q}+\frac{2}{\beta}\delta^{2q-1}\right)
\leq 2(1+\gamma \beta + 2/\beta):=C_4(\gamma,\beta),
\label{eq:a1QDQED}\\
a_2 &\leq 80 [1\vee \gamma^2 (\beta+1)^2] \delta^{4s-4} \frac{\delta}{\abs\zeta}(\beta+1)+\frac{3}{4}a_1^2\notag\\
&= 80\frac{1+\beta}{\beta} [1\vee \gamma^2 (\beta+1)^2] \delta^{4s+2q-4}+\frac{3}{4}C_4(\gamma,\beta)^2\notag\\
&\leq C_5(\gamma,\beta) \delta^{\min(4s+2q-4,0)},
\label{eq:a2QDQED}
\end{align}
where
\begin{align}
C_5(\gamma,\beta) = 80\frac{1+\beta}{\beta} [1\vee \gamma^2 (\beta+1)^2]+\frac{3}{4}C_4(\gamma,\beta)^2.
\label{eq:A2}
\end{align}

Next, substituting  \eqref{eq:abs_zeta} - \eqref{eq:Nmu}, and \eqref{eq:a1QDQED}, \eqref{eq:a2QDQED} into \eqref{eq:lyapunovQDQED1}, we obtain
 \begin{align}
\E[\tilde X_\infty^4 \mathds{1}_{\{\tilde X_\infty \leq -\zeta\}}] &\leq (1+\beta)\{48[1\vee \gamma(1+\beta)]\delta^{2s-2}
+ 80[1\vee \gamma^2(1+\beta)^2] \delta^{4s-3}+ 70[1\vee \gamma^3(1+\beta)^3]\delta^{6s-4}\}\notag\\
&+C_4(\gamma,\beta) (1+\beta) \{6+12 [1\vee \gamma(1+\beta)]+10 [1\vee \gamma^2(1+\beta)^2]  \}\delta^{-1} \notag\\
&+C_5(\gamma,\beta)(1+\beta) \{2+2[1\vee \gamma (1+\beta)]\} \delta^{\min(4s+2q-4,0)}\delta^{-1}\notag\\
&\leq C_6(\gamma,\beta)\delta^{\min(4s+2q-4,0)}\delta^{-1},
\label{eq:4th_order_mom}
\end{align}
where 
\begin{align}
C_6(\gamma,\beta)&=(1+\beta)\{48[1\vee \gamma(1+\beta)]+ 80[1\vee \gamma^2(1+\beta)^2]+ 70[1\vee \gamma^3(1+\beta)^3]\}\notag\\
&+C_4(\gamma,\beta) (1+\beta) \{6+12 [1\vee \gamma(1+\beta)]+10 [1\vee \gamma^2(1+\beta)^2]  \}\notag\\
&+C_5(\gamma,\beta)(1+\beta) \{2+2[1\vee \gamma (1+\beta)]\}.
\label{eq:A3}
\end{align}
Applying Jensen's inequality to~\eqref{eq:4th_order_mom} proves (3.11) and thus Lemma 3.1 in the main paper. 


\section{Discrete-time queue with general arrival distribution}
\label{sec:general_arr}

We consider the customer count process, $X = \{X_k: k=0,1,\dotsc\}$, with a general arrival distribution.
Specifically, comparing with the discrete queue studied in the main paper~\cite{FengShi2016},
we assume that the arrivals $\{A_k:k=0,1,\dotsc\}$ form an i.i.d. sequence
and follows a general distribution $G(\cdot)$ such that
\begin{itemize}
\item Variance of the distribution is $\sigma^2_A <\infty$;
\item Third non-central moment of the distribution is $\mu_3 < \infty$ .
\end{itemize}
Note that the Poisson arrival case studied in the main paper
is one special case satisfying the two conditions above.

Let $X_\infty$ and $\tilde X_\infty$ be the steady-state customer count and the scaled version of it,
defined in Section \blue{1} of the main paper. 
Let $Y_\infty$ be the continuous random variable having the following density
\begin{align}
    p(x) \propto \frac{2}{a(x)}\exp\left(\int^x_0 \frac{2b(y)}{a(y)}\mathrm dy\right),
    \quad x\in \R ,
\label{eq: py}
\end{align}
where $b(x)$ is the same as the one defined by (\blue{1.5}) of the main paper, and 
\begin{align}
 a(x)
 &=
 \begin{cases}
 \mu\left[(c_A - 1) + \Lambda\right],\quad x\leq -1/\delta,\\
 \mu\left(c_A - \mu + \delta(1-\mu)x + \mu x^2\right), \quad x\in [-1/\delta,\abs{\zeta}],
 \\
 \mu\left(c_A - \mu + \delta(1-\mu)\abs{\zeta} + \mu \zeta^2\right),
 \quad x\geq \abs \zeta .
 \end{cases}
 \label{eq:aform-1}
 \end{align}
Here,
\begin{align}
 c_A\equiv  \sigma^2_A /\Lambda + 1.
 \label{eq:ca}
 \end{align}
Note that when the arrival distribution is Poisson, $c_A = 2$. In that case,~\eqref{eq:aform-1} coincides with the definition of $a(x)$ in (\blue{1.6}) of the main paper.

Note that $p(x)$ is the stationary density of the diffusion process 
\begin{align}
G_Y f(x) = b(x) f'(x)  + \frac{1}{2} a(x) f''(x),\quad x\in \R,\,f\in C^2(\R).
\label{eq:ggy}
\end{align}

Next, we state the main theorem for this discrete-time queueing system with a general arrival distribution.

\begin{theorem-sec}
\label{thm:gdtmc}
Consider the DTMC $X$ with arrivals $\{A_k:k=0,1,\dotsc\}$ following distribution $G(\cdot)$ 
such that 
\begin{align}
\text{Var}(A_0)/\Lambda = c_A-1<\infty,\quad \E A_0^3 /\Lambda = v_A < \infty.
\label{eq:ca_va}
\end{align}
For all $N\geq 1$,
 $\,\mu\in (0,1)$
satisfying $1\leq R < N$ and
\begin{equation*} 
\mu = \gamma R^{-s}, \quad N - R = \beta R^q 
\end{equation*}
for some $s\geq 1,\,\frac{1}{2}\leq q\leq 1$. 
The Wasserstein distance
between $\tilde X_\infty$ and $Y_\infty$
\begin{equation}
d_{W}\left(\tilde X_\infty, \ Y_\infty \right)  \leq  C(\gamma,\beta, c_A, v_A)\delta ,
    \label{eq:gdtmc}
\end{equation}
where 
\begin{align}
C(\gamma,\beta, c_A, v_A) &=
\gbc\left(1 +\frac{1}{\beta}\right) 
\left[
\frac{2}{3} v_A +\frac{10}{3} (1+\beta)\left\{1\vee \left(\gamma(1+\beta)\right)^2\right\}
\right]
.
\label{eq:gcgammabetacv}
\end{align}
Here, $\gbc=\gbc(\gamma,c_A)$ is a constant depending only on $\gamma$ and $c_A$,
with the explicit form specified in Lemma~\ref{lem:ggradboundsCW}.
Note that in the Poisson arrival case, $c_A = 2$ and $\gbc$ coincides 
with its counterpart in Lemma \blue{2} of the main paper.
\end{theorem-sec}

\begin{proof}[Proof of Theorem~\ref{thm:gdtmc}]
The basic framework to prove Theorem~\ref{thm:gdtmc} is the same
as that in the main paper.


For any $h \in \lipone$,
let $f=f_h$ be a solution to the Poisson equation
 \begin{equation}
 G_Y f(x) = \E \left[h(Y_\infty)\right]-h(x),\quad x\in \mathbb R .
 \label{eq:gpoisson}
 \end{equation}

After the generator coupling via the Poisson equation
\begin{align}
\Big|\E h(\tilde X_\infty) - \E h(Y_\infty)|
&=
\Big|\E \left[G_{\tilde{X}}f(x) - G_Y f(x)\right]\Big|,
\label{eq:generator_couple}
\end{align}
we perform the following Taylor expansion for any given $x = \delta(n - x_\infty)$ and $n=0,1,\dotsc$,
\begin{align}
G_{\tilde X}f(x) &= \E_n [f(x+\delta (A_0 - D_0))] - f(x)\notag\\
&=f'(x) \delta \E_n(A_0 -D_0) +\frac{1}{2}f''(x)\delta^2 \E_n[(A_0-D_0)^2]+\frac{1}{6}\delta^3 \E_n [f'''(\xi)(A_0-D_0)^3] ,
\label{eq:gte} 
\end{align}
where \begin{equation*}
    \abs{\xi - x}\leq \delta \abs{A_0 - D_0}.
\end{equation*} 
It can be easily verified that  
\begin{align*}
\delta \E_n(A_0 -D_0)&= b(x),
\end{align*}
and
\begin{align}
\delta^2 \E_n[(A_0-D_0)^2]&=a(x),
\label{eq:eq_ax}
\end{align}
where~\eqref{eq:eq_ax} follows from the calculation below
\begin{align}
\delta^2 \E_n [(A_0 -D_0)^2] &=\delta^2 Var_n (A_0-D_0) +\delta^2 [\E_n(A_0-D_0)]^2 \nonumber \\
&=\delta^2 \sigma_A^2 +\delta^2 [N-(n-N)^-]\mu(1-\mu) +b^2(x)\nonumber \\
&=\delta^2 \sigma_A^2 +\left[\delta^2 N \mu -\delta (x+\zeta)^- \mu\right](1-\mu) +b^2(x)\nonumber \\
&=\delta^2 \sigma_A^2 +\left[-\delta b(x) + \delta^2 \Lambda\right](1-\mu) +b^2(x)\nonumber \\
&=\delta^2\sigma_A^2+ \delta^2 \Lambda(1-\mu) -(1-\mu)\delta b(x) +b^2(x).
 \label{eq:gsecondorder}
\end{align}

Combining~\eqref{eq:generator_couple} and~\eqref{eq:gte} implies that
\begin{align}
\abs{\E h(\tilde X_\infty) - \E h(Y_\infty)}
&=
\Big|\E \left[G_{\tilde{X}}f(x) - G_Y f(x)\right]\Big|\notag\\
&\leq 
\frac{1}{6}\delta^3 \E\left\{\E_{X_\infty} \left[\norm{f'''}\big|A_0-D_0\big|^3\right]\right\},
\label{eq:gerrors}
\end{align}
where $\norm{f'''} = \max_{x\in\R}\abs{f'''(x)}$.

For~\eqref{eq:gerrors}, one first uses the $c_r$-inequality 
and equation~\eqref{eq:binomial3rdmom} to obtain
\begin{align}
\E_n\abs{A_0-D_0}^3 \leq 4 \E A_0^3 + 20 \max \{1,\,(N\mu)^2\}\,N\mu
=4 \frac{\mu_3}{\Lambda}\,\Lambda +
20 \max \{1,\,(N\mu)^2\}\,N\mu ,
\label{eq:gaminusd}
\end{align}
for each $n=0,1,\dotsc$. 
Then, applying the gradient bound~\eqref{eq:grad_bdd_3}
stated in Lemma~\ref{lem:ggradboundsCW} below, we get 
\begin{align}
\frac{1}{6}\delta^3 \E\left\{\E_{X_\infty} \left[\norm{f'''}\big|A_0-D_0\big|^3\right]\right\}
&\leq \frac{1}{6}
\gbc(\mu\delta^{-1},c_A)\delta^3\frac{1}{\mu} \left(1 + \frac{1}{\abs\zeta}\right)
\left[4 \frac{\mu_3}{\Lambda}\,\Lambda +
20 \max \{1,\,(N\mu)^2 \}\,N\mu\right]\notag\\
&\leq \left[\frac{2}{3} \frac{\mu_3}{\Lambda} +
\frac{10}{3} \max \{1,\,(N\mu)^2\}\frac{N}{R}\right]
\gbc(\mu\delta^{-1},c_A)\left(1 +\frac{1}{\abs \zeta}\right) 
\delta . 
\label{eq:gerrors4th}
\end{align}


Recall the characterizations of $\abs\zeta$, $N\delta^2$, $N\mu$, 
and $\mu/\abs{\zeta}$ in~\eqref{eq:abs_zeta}-\eqref{eq:mu_zeta} of this supplement.  
Under the assumptions on $s$ and $q$ in Theorem~\ref{thm:gdtmc}, 
we have $\abs\zeta \geq \beta$ 
and $N\mu \leq \gamma (\beta + 1)$.
Applying these characterizations to~\eqref{eq:gerrors4th}, 
we obtain through~\eqref{eq:gerrors} that 
\begin{align}
\abs{\E h(\tilde X_\infty) - \E h(Y_\infty)} &= \abs{\E G_{\tilde X}f(\tilde X_\infty) - \E G_Yf(\tilde X_\infty)}\notag\\
&\leq 
\gbc(\gamma,c_A)\left(1 +\frac{1}{\beta}\right) 
\left[
\frac{2}{3} v_A +\frac{10}{3} (1+\beta)\left\{1\vee \left(\gamma(1+\beta)\right)^2\right\}
\right]
\delta
,
\label{eq:gbdd}
\end{align}
where the inequality comes from the observation that $C_0(\cdot,c_A)$ is an increasing function in its first variable, and that $\mu\delta^{-1} \leq \gamma$.

This proves Theorem~\ref{thm:gdtmc}.
\end{proof}

\begin{lemma-sec}[Gradient bounds]
\label{lem:ggradboundsCW}
Fix an $h\in \lipone$ with $h(0) = 0$.
There exists a solution $f_h$ to the Poisson equation,
\begin{align}
G_Y f(x)  = \E h(Y) - h(x), \quad x \in \R,
\label{eq:pois_eq}
\end{align} 
that is twice continuously differentiable, with an absolutely continuous second derivative, and for all $\Lambda>0$, $N\geq 1$, and $\mu\in (0,1)$ satisfying $1 \leq R < N$, 
\begin{align}
\abs{f_h'(x)} &\leq \begin{cases}
\frac{\tilde C_1}{\mu} (1 + 1/\abs{\zeta}),
\quad x\leq -\zeta ,\\
\frac{1}{2} + \frac{1}{\mu\abs{\zeta}}
\left[x+\left(\tilde C + \frac{\delta}{2}\right)
+\left(\tilde C + \frac{c_A}{2}\right) \frac{1}{\abs \zeta}
\right]
,\quad x\geq -\zeta,
\end{cases}
\label{eq:grad_bdd_1}\\
\abs{f_h''(x)} & \leq 
\begin{cases}
\frac{\tilde C_2}{\mu}(1 + 1/\abs{\zeta}),
\quad x\leq -\zeta , \\
\frac{1}{\mu \abs{\zeta}},\quad x\geq -\zeta,
\end{cases}
\label{eq:grad_bdd_2} \\
\abs{f_h'''(x)} &\leq \begin{cases}
\frac{\gbc}{\mu} (1 + 1/\abs{\zeta}),
\quad x\leq -\zeta , \\
\frac{4}{c_A - 1}\frac{1}{\mu},
\quad x \geq -\zeta, 
\end{cases},
\label{eq:grad_bdd_3}
\end{align}
where 
\begin{align}
\tilde C_1 &= \tilde{C}_1(\mu\delta^{-1},c_A)=
\tilde C　e^{\frac{1}{c_A - 1}}
\left(3 +\frac{2 }{c_A - 1}
+ \frac{2}{c_A - 1}e^{\frac{1}{c_A - 1}}\right),
\label{eq:gb_const_1} \\
\tilde C_2 &=\tilde{C}_2(\mu\delta^{-1},c_A)=
e^{\frac{1}{c_A - 1}}
\left(1 +\frac{2 }{c_A - 1}
+ \frac{2}{c_A - 1}e^{\frac{1}{c_A - 1}}\right)
\bigg[1+(1+\tilde C )\left(\frac{\mu\delta^{-1}}{c_A-1}
\vee\frac{1}{c_A-1}\vee 2\right) 
\notag\\
&\hspace{3cm}+  \tilde C_1 \left(\frac{c_A-\mu}{c_A-1}\vee 3\right)\bigg],
\label{eq:gb_const_2}\\
\gbc &= \gbc(\mu\delta^{-1},c_A)=\frac{4}{c_A-1}
\left[1+(1+\tilde C )\left(\frac{\mu\delta^{-1}}{c_A-1}
\vee\frac{1}{c_A-1}\vee 2\right)+  \tilde C_1 \left(\frac{c_A-\mu}{c_A-1}\vee 3\right)\right],
\label{eq:gb_const_3}
\end{align}  
and \begin{align}
\tilde C &= \tilde{C}(\mu \delta^{-1},c_A)
\notag\\
&= \frac{1 + \sqrt{2}}{2}  \mu \delta^{-1} \left(1\vee \mu \delta^{-1}\right)
+ \left[2 + \sqrt{\mu \delta^{-1}}
\left(1\vee \sqrt{\mu \delta^{-1}}\right)\right]
\left[c_A + \frac{3}{2} (1-\mu) \delta^2\right]
.
\label{eq:tilde_c}
\end{align}
\end{lemma-sec}

We leave the complete details of the proof for these gradient bounds 
to the last section due to its complexities.


\section{Gradient bounds for state-dependent diffusion process}
\label{sec:general-gradient-bounds}

To establish the gradient bounds in Lemma~\ref{lem:ggradboundsCW}, we first define the following useful quantity for notational convenience
\begin{align}
 r(x)&\equiv \frac{2b(x)}{a(x)} = 
 \begin{cases}
 \frac{-2x}{(c_A-1) + \Lambda}, \quad x \leq -1/\delta, \\
 \frac{-2x}{c_A-\mu + \delta(1-\mu) x + \mu x^2}, \quad x \in [-1/\delta, \abs{\zeta}], \\
 \frac{-2\abs{\zeta}}{c_A-\mu + \delta(1-\mu) \abs{\zeta} + \mu \zeta^2}, \quad x \geq \abs{\zeta}.
 \end{cases}
 \label{eq:rform-1}
 \end{align}
 Note that for $x\in [-1/\delta,0]$,
 \begin{align}
 r(x) &\leq \frac{-2x}{c_A - \mu + \delta(1-\mu) (-1/\delta)} 
 = \frac{-2x}{c_A - 1}	,
 \label{eq:r-ineq}
 \end{align}
 and for $x\in [0,-\zeta]$,
 \begin{align}
 -r(x) &\leq \frac{2x}{c_A - \mu} \leq \frac{2}{c_A - 1}x.
 \label{eq:r-ineq-1}
 \end{align}
These inequalities will turn out to be of use in the proof later.

Then, we needs the following three lemmas as preliminaries, which will be proven at the end of this section.
Note that these three lemmas are more general versions of Lemma \blue{3.7}, \blue{3.8}, and \blue{3.9} in the main paper, as they hold for any arrival distribution satisfying the conditions~\eqref{eq:ca_va} now.
\begin{lemma-sec}
\label{lem:ints_for_grad_bdds}
Recall that 
$q(x)$ is defined by \blue{(3.10)} of the main paper~\cite{FengShi2016}.
For all $\Lambda>0$, $N\geq 1$ and $\mu\in(0,1)$ satisfying $1\leq R<N$,
	\begin{align}
	\frac{1}{q(x)} \int^x_{-\infty} \frac{2}{a(y)} q(y) \mathrm dy &\leq 
	\begin{cases}
	\frac{1}{\mu}, \quad x \leq - 1,\\
	\left(1 + \frac{2}{c_A - 1} e^{\frac{1}{c_A - 1}}\right)\frac{1}{\mu}, \quad x\in [-1,0],\\
	e^{\frac{1}{c_A-1}\zeta^2}  \left(1 + \frac{2}{c_A-1}e^{\frac{1}{c_A-1}}	+  \frac{2}{c_A-1}	\abs{\zeta}	\right)\frac{1}{\mu},
	\quad x\in [0,-\zeta].
	\end{cases} 
	\label{eq:int-1}\\
	\frac{1}{q(x)} \int_x^\infty \frac{2}{a(y)} q(y) \mathrm dy &\leq
	\begin{cases}
	\left(\frac{1}{\eta}+\frac{2\eta}{c_A - 1}e^{\frac{\eta^2}{c_A-1}}\right) \frac{1}{\mu},\quad x \in [0,\eta], \, \eta \leq -\zeta ,\\
	\frac{1}{\mu \abs{\zeta}},\quad x \geq -\zeta.
	\end{cases}
	\label{eq:int-2}\\
	\frac{1}{q(x)} \int^x_{-\infty} \frac{2\abs{y}}{a(y)}q(y)\mathrm dy &\leq 
	\begin{cases}
	\frac{1}{\mu},\quad x\leq 0,\\
	2 e^{\frac{1}{c_A - 1 }\zeta^2}\frac{1}{\mu},\quad x \in [0,-\zeta].
	\end{cases}
	\label{eq:int-3}\\
	\frac{1}{q(x)} \int^\infty_x \frac{2 \abs{y}}{a(y)} q(y) \mathrm dy &\leq 
	\begin{cases}
	\frac{3}{2}\frac{1}{\mu} + \frac{\delta}{2}\frac{1}{\mu \abs{\zeta}}
	+ \frac{c_A}{2}\frac{1}{\mu\zeta^2},
	\quad x \in [0,-\zeta], \\
	\frac{x}{\mu \abs{\zeta}} + \frac{\delta}{2}\frac{1}{\mu \abs{\zeta}}
	+ \frac{c_A}{2}\frac{1}{\mu\zeta^2}+\frac{1}{2},
	\quad x \geq -\zeta .
	\end{cases}
	\label{eq:int-4}\\
	\frac{\abs{r(x)}}{q(x)} \int^x_{-\infty}\frac{2}{a(y)}q(y)\mathrm dy &\leq 
	\frac{2}{c_A-1}\frac{1}{\mu},\quad x \leq 0.
	\label{eq:int-5}\\
	\frac{\abs{r(x)}}{q(x)} \int^\infty_x \frac{2}{a(y)} q(y) \mathrm dy 
	&\leq \frac{2}{c_A - 1 }\frac{1}{\mu}, \quad x\geq 0,
	\label{eq:int-6}
	\end{align}
and when 
$\abs{\zeta} \geq 1$,
\begin{align}
\frac{1}{q(x)} \int_x^\infty \frac{2}{a(y)} q(y) \mathrm dy &\leq 
\begin{cases}
\left(1 + \frac{2}{c_A - 1} e^{\frac{1}{c_A - 1}}\right)\frac{1}{\mu},\quad x\in [0,1],\\
\frac{1}{\mu},\quad x \geq 1.
\end{cases}
\label{eq:int-22}
\end{align}
\end{lemma-sec}

\begin{lemma-sec}
\label{lem:ey}
Let the random variable $Y_\infty$ have the stationary distribution of a diffusion process with drift $b(x)$ and \emph{state-dependent} diffusion coefficient $a(x)$. For all $\Lambda>0$, $N\geq 1$ and $\mu\in(0,1)$ satisfying $1\leq R<N$,
\begin{align}
\E \abs{Y_\infty}&\leq 
 \tilde C \left(1 + 1/\abs{\zeta}\right)
 ,
\label{eq:ey}
\end{align}
where $\tilde C$ is specified by~\eqref{eq:tilde_c} in Lemma~\ref{lem:ggradboundsCW}. 
\end{lemma-sec}

\begin{lemma-sec}
\label{lem:a}
Recall the form of $a(x)$ and $r(x)$ in~\eqref{eq:aform-1} and~\eqref{eq:rform-1}. For all $\Lambda>0$, $N\geq 1$ and $\mu\in(0,1)$ satisfying $1\leq R<N$,
\begin{align}
a(x) &\geq (c_A - 1)\mu ,\quad x\in \mathbb{R}; 
\label{eq:abdd_1}\\
\frac{\abs{x a'(x)}}{a(x)} &\leq \left(\frac{1 - \mu }{c_A-1}\vee 2\right) \mathds{1}_{\{x\in (-1/\delta,-\zeta]\}},
\quad x\in\mathbb{R};
\label{eq:abdd_2}\\
\E\abs{Y_\infty}\frac{\abs{a'(x)}}{a(x)}&\leq \tilde C \left(\frac{\mu\delta^{-1}}{c_A-1}
\vee \frac{1}{c_A - 1}\vee 2\right)(1+1/\abs \zeta)
\mathds{1}_{\{x\in (-1/\delta,-\zeta]\}},
\quad x\in \mathbb{R};
\label{eq:abdd_3} \\
\abs{r'(x) a(x) } &\leq 
2 \left(\frac{c_A -\mu}{c_A - 1}\vee 3\right)\mu 
\mathds{1}_{\{x \leq -\zeta\}} 
,\quad x\in \mathbb {R};
\label{eq:abdd_4}
\end{align}
where $a'(x)$ and $r'(x)$ are interpreted as the left derivative at $x = 1/\delta$
and $x = -\zeta$.
\end{lemma-sec}

Now we are going to prove Lemma~\ref{lem:ggradboundsCW} using these three lemmas.

Following from Chapter 3 of the doctoral thesis~\cite{BravDai20166}, the derivatives of $f_h(x)$ have the following forms:
\begin{align}
f'_h(x)&= 
\frac{1}{q(x)} 
\int^x_{-\infty} \frac{2}{a(y)}
\left(\E h(Y_\infty) - h(y) \right)
q(y) \mathrm dy , \quad x\in \mathbb R.
\label{eq:fprime_1} \\
f'_h(x)&= -
\frac{1}{q(x)} 
\int_x^{\infty} \frac{2}{a(y)}
\left(\E h(Y_\infty) - h(y) \right)
q(y) \mathrm dy , \quad x\in \mathbb R.
\label{eq:fprime_2}\\
f''_h(x) &= \frac{1}{q(x)}
\int^x_{-\infty} \left\{
-\frac{2}{a(y)}h'(y)
-\frac{2a'(y)}{a^2(y)} \left[\E h(Y_\infty)
- h(y)\right] 
-r'(y) f'_h(y)
\right\} q(y)\mathrm dy ,
\label{eq:fprimeprime_1}\\
f''_h(x) &= \frac{1}{q(x)}
\int_x^{\infty} \left\{
-\frac{2}{a(y)}h'(y)
-\frac{2a'(y)}{a^2(y)} \left[\E h(Y_\infty)
- h(y)\right] 
-r'(y) f'_h(y)
\right\} q(y)\mathrm dy ,
\label{eq:fprimeprime_2}\\
f_h'''(x)&=
-r'(x) f'_h(x) - r(x) f_h''(x) -\frac{2}{a(x)}h'(x)
-\frac{2a'(x)}{a^2(x)} \left[\E h(Y_\infty)
- h(x)\right],
\label{eq:fprimeprimeprime}
\end{align}
where $a'(x)$ is interpreted as the left derivative at 
the points $x = -1/\delta$ and $x = -\zeta$.

Then, the properties of $h$ implies that 
\begin{align}
\abs{f_h'(x)} &\leq 
\frac{1}{q(x)} \int^x_{-\infty}
\frac{2}{a(y)} \left( \E\abs{Y_\infty}+\abs y \right) q(y)\mathrm dy ,
\label{eq:abs_fprime_1}\\
\abs{f_h'(x)} &\leq 
\frac{1}{q(x)} \int_x^{\infty}
\frac{2}{a(y)} \left( \E\abs{Y_\infty}+\abs y \right) q(y)\mathrm dy .
\label{eq:abs_fprime_2}
\end{align}

For $x\leq 0$, applying~\eqref{eq:int-1}, ~\eqref{eq:int-3},
and~\eqref{eq:ey}
to~\eqref{eq:abs_fprime_1} gives 
\begin{align}
\abs{f_h'(x)} &\leq  
\frac{1}{\mu} \left(1 + 1/\abs \zeta\right)
\left[
1 + \tilde{C} \left(1 + \frac{2}{c_A-1} e^{\frac{1}{c_A-1}}\right)
\right]  .
\label{eq:abs_fprime_case1}
\end{align}
For $x \in [0,-\zeta]$, we need to consider separately the cases
when $\abs \zeta\leq 1$ and $\abs \zeta\geq 1$.

When $\abs \zeta\leq 1$, 
applying~\eqref{eq:int-1}, ~\eqref{eq:int-3}
and~\eqref{eq:ey}
to~\eqref{eq:abs_fprime_1} gives 
\begin{align}
\abs{f_h'(x)} &\leq  
\frac{1}{\mu} \left(1 + 1/\abs \zeta\right)
\tilde{C}\left[
e^{\frac{1}{c_A-1}}
\left(3 +\frac{2}{c_A - 1} +\frac{2}{c_A-1} e^{\frac{1}{c_A-1}}\right)
\right]  , \quad x\in [0,-\zeta] ,
\label{eq:abs_fprime_case3_1}
\end{align}
and when $\abs \zeta\geq 1$, applying~\eqref{eq:int-2}, ~\eqref{eq:int-22}
and~\eqref{eq:ey}
to~\eqref{eq:abs_fprime_2} gives 
\begin{align}
\abs{f_h'(x)} 
&\leq \frac{1}{\mu} (1 + 1/\abs \zeta)
\left[
\tilde C 
\left(1 +\frac{2}{c_A -1} e^{\frac{1}{c_A-1}}\right) 
+ 1 + \frac{c_A}{2}
\right]
  , \quad x\in [0,-\zeta] .
\label{eq:abs_fprime_case3_2}
\end{align}
Using~\eqref{eq:abs_fprime_case1}-\eqref{eq:abs_fprime_case3_2} together, along with the observation that $2\tilde C > 1 + \frac{c_A}{2}$, proves the first half
of~\eqref{eq:grad_bdd_1}.

For $x\geq -\zeta$, 
applying~\eqref{eq:int-2}, ~\eqref{eq:int-4}
and~\eqref{eq:ey}
to~\eqref{eq:abs_fprime_2} gives 
\begin{align}
\abs{f_h'(x)} &\leq  
\frac{1}{2} + \frac{1}{\mu \abs \zeta} \left[
x + \left(\tilde C +\frac{\delta}{2}\right)
+\left(\tilde C +\frac{c_A}{2}\right)\frac{1}{\abs \zeta}
\right] ,
\label{eq:abs_fprime_case2}
\end{align}
which proves the second half of~\eqref{eq:grad_bdd_1}.

Now, we move on to deal with~\eqref{eq:grad_bdd_2} and~\eqref{eq:grad_bdd_3}.

Since the function $h$ satisfies $\abs{h(x)} \leq \abs{x}$ for all $x\in\mathbb R$
and $\norm{h'} \leq 1$,~\eqref{eq:fprimeprime_1} and~\eqref{eq:fprimeprime_2} imply that
\begin{align}
\abs{f''_h(x)} &\leq \frac{1}{q(x)}
\int^x_{-\infty} \left\{
\frac{2}{a(y)}
+\E \abs{Y_\infty}\frac{2\abs{a'(y)}}{a^2(y)} 
+\frac{2\abs{ya'(y)}}{a^2(y)}
+\abs{r'(y) f'_h(y)}
\right\} q(y)\mathrm dy ,
\label{eq:abs_fprimeprime_1}\\
\abs{f''_h(x)} &\leq \frac{1}{q(x)}
\int_x^{\infty} \left\{
\frac{2}{a(y)}
+\E\abs{Y_\infty}\frac{2\abs{a'(y)}}{a^2(y)} 
+\frac{2\abs{ya'(y)}}{a^2(y)}
+\abs{r'(y) f'_h(y)}
\right\} q(y)\mathrm dy .
\label{eq:abs_fprimeprime_2}
\end{align}

Thus, when $x\leq -\zeta$, applying~\eqref{eq:abdd_3},~\eqref{eq:abdd_2},~\eqref{eq:abdd_4} and~\eqref{eq:grad_bdd_1} to~\eqref{eq:abs_fprimeprime_1} implies that 
\begin{align}
\abs{f''_h(x)} &\leq 
\frac{1}{q(x)} \int^x_{-\infty}
\frac{2}{a(y)}q(y)
\bigg\{
1 + \tilde C \left(\frac{\mu\delta^{-1}}{c_A-1}
\vee \frac{1}{c_A - 1}\vee 2\right)(1+1/\abs \zeta)
\mathds{1}_{\{y\in (-1/\delta,-\zeta]\}}
\notag\\
&\hspace{2cm}+ \left(\frac{1 - \mu }{c_A-1}\vee 2\right) \mathds{1}_{\{y\in (-1/\delta,-\zeta]\}}
+ \tilde{C_1} \left(\frac{c_A - \mu}{c_A - 1}\vee 3\right)
(1 + 1/\abs \zeta)\mathds{1}_{\{y \leq -\zeta\}}
\bigg\} \mathrm{dy}
\notag 
\\
&\leq 
\hat C (1 + 1/\abs \zeta) 
\frac{1}{q(x)} \int^x_{-\infty}
\frac{2}{a(y)}q(y) \mathrm{dy},
\label{eq:abs_fprimeprime_3}
\end{align}
where 
\begin{align*}
\hat C &= 1 + (1 +\tilde{C})\left(\frac{\mu\delta^{-1}}{c_A-1}
\vee \frac{1}{c_A - 1}\vee 2\right)
+\tilde{C_1} \left(\frac{c_A - \mu}{c_A - 1}\vee 3\right)
;
\end{align*}
and when $x\geq 0$, we again apply~\eqref{eq:abdd_3},~\eqref{eq:abdd_2},~\eqref{eq:abdd_4} and~\eqref{eq:grad_bdd_1},
but this time to~\eqref{eq:abs_fprimeprime_2},
to see that 
\begin{align}
\abs{f''_h(x)} &\leq 
\frac{1}{q(x)} \int_x^{\infty}
\frac{2}{a(y)}q(y)
\bigg\{
1 + \tilde C \left(\frac{\mu\delta^{-1}}{c_A-1}
\vee \frac{1}{c_A - 1}\vee 2\right)(1+1/\abs \zeta)
\mathds{1}_{\{y\in (-1/\delta,-\zeta]\}}
\notag\\
&\hspace{2cm}+ \left(\frac{1 - \mu }{c_A-1}\vee 2\right) \mathds{1}_{\{y\in (-1/\delta,-\zeta]\}}
+ \tilde{C_1} \left(\frac{c_A - \mu}{c_A - 1}\vee 3\right)
(1 + 1/\abs \zeta)\mathds{1}_{\{y \leq -\zeta\}}
\bigg\} \mathrm{dy}
\notag\\
&\leq 
\left[1+(\hat C - 1)(1 + 1/\abs{\zeta})\mathds{1}_{\{
x< -\zeta
\}}
\right]
\frac{1}{q(x)} \int_x^{\infty}
\frac{2}{a(y)}q(y) \mathrm{dy}
,
\label{eq:abs_fprimeprime_4}
\end{align}
where for the last inequality we used that for all $y> x$,
$\mathds{1}_{\{y \leq -\zeta\}}\leq \mathds{1}_{\{
x< -\zeta
\}}$.

Therefore, when $x\leq 0$, applying~\eqref{eq:int-1}
to~\eqref{eq:abs_fprimeprime_3} implies that
\begin{align}
\abs{f_h''(x)} &\leq 
\frac{\hat C}{\mu} (1 + 1/\abs \zeta)
\left(1 +\frac{2}{c_A-1} e^{\frac{1}{c_A-1}}\right).
\label{eq:abs_fprimeprime_5}
\end{align} 
For $x\in [0,-\zeta]$, we need to consider separately
the cases when $\abs \zeta\leq 1$ and $\abs \zeta\geq 1$.
When $\abs \zeta\leq 1$, applying~\eqref{eq:int-1}
to~\eqref{eq:abs_fprimeprime_3} implies that
\begin{align}
\abs{f_h''(x)} &\leq 
\frac{\hat C}{\mu} (1 + 1/\abs \zeta)\left[e^{\frac{1}{c_A-1}}
\left(1 +\frac{2}{c_A-1}+\frac{2}{c_A-1} e^{\frac{1}{c_A-1}}\right)\right],
\quad x\in [0,-\zeta],
\label{eq:abs_fprimeprime_6}
\end{align} 
and when $\abs \zeta\geq 1$, applying~\eqref{eq:int-22}
to~\eqref{eq:abs_fprimeprime_4} implies that
\begin{align}
\abs{f_h''(x)} &\leq 
\frac{\hat C}{\mu} (1 + 1/\abs \zeta)
\left(1 +\frac{2}{c_A-1} e^{\frac{1}{c_A-1}}\right),
\quad x\in [0,-\zeta] .
\label{eq:abs_fprimeprime_7}
\end{align} 
Combining the bounds in~\eqref{eq:abs_fprimeprime_5}-\eqref{eq:abs_fprimeprime_7} proves the first half
of~\eqref{eq:grad_bdd_2}.

When $x\geq -\zeta$,
applying~\eqref{eq:int-2} to~\eqref{eq:abs_fprimeprime_4}
immediately establishes the remaining half of~\eqref{eq:grad_bdd_2}.

Now we move on to~\eqref{eq:grad_bdd_3}. From the form of
$f_h'''(x)$ in~\eqref{eq:fprimeprimeprime}, along with the properties of the function $h$, we see immediately
that 
\begin{align*}
\abs{f_h'''(x)}&\leq 
\abs{r'(x) f'_h(x)} + \abs{r(x) f_h''(x)} +\frac{2}{a(x)}
+\frac{2\abs{a'(x)}}{a^2(x)} \E \abs{Y_\infty}
+\frac{2\abs{xa'(x)}}{a^2(x)} .
\end{align*}
Applying the bound on $\abs{r'(x)}$ in~\eqref{eq:abdd_4},
the bound on $\abs{f'_h(x)}$ in~\eqref{eq:grad_bdd_1},
and the bounds~\eqref{eq:abdd_1}-\eqref{eq:abdd_3},
we have that
\begin{align}
\abs{f_h'''(x)}&\leq
\frac{2}{c_A-1}\frac{1}{\mu}
\left[
1 + (\hat C - 1)(1+1/\abs \zeta) \mathds{1}_{\{x\leq -\zeta\}}
\right]
+\abs{r(x) f_h''(x)}
.
\label{eq:abs_fprimeprimeprime}
\end{align}
For the last term of the right side above, $\abs{r(x) f_h''(x)}$,
one can simply multiply both sides of~\eqref{eq:abs_fprimeprime_3}
and~\eqref{eq:abs_fprimeprime_4} by $\abs{r(x)}$,
and then apply~\eqref{eq:int-5} and~\eqref{eq:int-6}
to arrive at
\begin{align}
\abs{r(x) f_h''(x)} &\leq \begin{cases}
\frac{2}{c_A-1} \frac{\hat C}{\mu} (1 + 1/\abs \zeta),
\quad x\leq -\zeta,\\
\frac{2}{c_A-1}\frac{1}{\mu},
\quad x\geq -\zeta .
\end{cases}
\label{eq:abs_fprimeprimeprime_1}
\end{align}
Combining~\eqref{eq:abs_fprimeprimeprime} and~\eqref{eq:abs_fprimeprimeprime_1} proves~\eqref{eq:grad_bdd_3}.

Finally, we are going to verify the three lemmas stated at the beginning of this section.
\begin{proof}[Proof of Lemma~\ref{lem:ints_for_grad_bdds}]
First, we claim that
\begin{align}
\frac{1}{q(x)} \int_{-\infty}^{x} \frac{2}{a(y)}q(y) dy \leq \frac{1}{b(x)}, \quad x < 0, 
\label{eq:b1} \\
\frac{1}{q(x)} \int_{x}^{\infty} \frac{2}{a(y)}q(y) dy \leq \frac{1}{\abs{b(x)}}, \quad x > 0. 
\label{eq:b2}
\end{align}
To see why, suppose that $x < 0$. Using the fact that $b(y)/b(x) \geq 1$ for $y \leq x$, there is 
\begin{align}
\frac{1}{q(x)} \int_{-\infty}^{x} \frac{2}{a(y)}q(y) dy \leq&\  \frac{1}{q(x)} \int_{-\infty}^{x} \frac{2 b(y)}{a(y)} \frac{1}{b(x)} q(y) dy \notag \\
=&\  \frac{1}{q(x)} \frac{1}{b(x)} \int_{-\infty}^{x} r(y) e^{\int_{0}^{y} r(u) du } dy \notag  \\
=&\  \frac{1}{q(x)} \frac{1}{b(x)} \big(q(x) - q(-\infty)\big)  \notag \\
\leq&\ \frac{1}{b(x)}. 
\label{eq:negpart}
\end{align}
The proof for~\eqref{eq:b2} is essentially the same and is omitted here.

These two inequalities imply immediately the first part of~\eqref{eq:int-1}
and the second part of~\eqref{eq:int-2}.

It remains to bound the integrals when $x \in [-1, 0]$, and $x \in [0,-\zeta]$. 

When $x \in [-1,0]$,
\begin{align*}
\frac{1}{q(x)} \int_{-\infty}^{x} \frac{2}{a(y)}q(y) dy =&\ \frac{q(-1)}{q(x)} \frac{1}{q(-1)} \int_{-\infty}^{-1} \frac{2}{a(y)}q(y) dy + \frac{1}{q(x)} \int_{-1}^{x} \frac{2}{a(y)}q(y) dy \\
\leq&\ \frac{q(-1)}{q(x)} \frac{1}{\mu} + \frac{1}{q(x)} \int_{-1}^{x} \frac{2}{a(y)}q(y) \mathrm  dy.
\end{align*}
Observe that
\begin{align*}
\frac{q(-1)}{q(x)} = e^{\int_{0}^{-1} r(u) du -\int_{0}^{x} r(u) du}=  e^{-\int_{-1}^{x} r(u) du} \leq 1.
\end{align*}
Furthermore, using the inequality about $r(x)$ from \eqref{eq:r-ineq} and the inequality $a(x) \geq (c_A-1)\mu$ for all $x \in \R$ from~\eqref{eq:abdd_1}, we see that
\begin{align*}
\frac{1}{q(x)} \int_{-1}^{x} \frac{2}{a(y)}q(y) dy &=  e^{\int_{x}^{0} r(u) du}\int_{-1}^{x} \frac{2}{a(y)}e^{-\int_{y}^{0} r(u) du} dy \\
&\leq\ e^{\int_{x}^{0} r(u) du}\int_{-1}^{0} \frac{2}{a(y)} dy \\ 
& \leq e^{\int_x^0 \frac{2}{c_A - 1} (-u) du} \frac{2}{(c_A - 1)\mu}\\
&\leq e^{\frac{1}{c_A-1}} \frac{2}{(c_A-1)\mu}.
\end{align*}
Hence, for $x \in [-1, 0]$,
\begin{align*}
\frac{1}{q(x)} \int_{-\infty}^{x} \frac{2}{a(y)}q(y) dy \leq 
\frac{1}{\mu}(1 + \frac{2}{c_A - 1}e^{\frac{1}{c_A-1 }}).
\end{align*}
This proves the second part of~\eqref{eq:int-1}.

Now fix $\eta > 0$ such that $\eta \leq \abs{\zeta}$. When $x \in [0,\eta]$,
\begin{align*}
\frac{1}{q(x)} \int_{x}^{\infty} \frac{2}{a(y)}q(y) dy =&\ \frac{1}{q(x)} \int_{x}^{\eta} \frac{2}{a(y)}q(y) dy + \frac{p(\eta)}{q(x)} \frac{1}{p(\eta)} \int_{\eta}^{\infty} \frac{2}{a(y)}q(y) dy \\
\leq&\ \frac{1}{q(x)} \int_{x}^{\eta} \frac{2}{a(y)}q(y) dy + \frac{p(\eta)}{q(x)} \frac{1}{\mu \abs{\eta}}.
\end{align*}
To bound the first term above, observe that $r(x) < 0$ and $\abs{r(x)} \leq \frac{2}{c_A-1}x$ for $x \in [0,-\zeta]$. Then
\begin{align*}
\frac{1}{q(x)} \int_{x}^{\eta} \frac{2}{a(y)}q(y) dy =&\ e^{-\int_{0}^{x} r(u) du} \int_{x}^{\eta} \frac{2}{a(y)} e^{\int_{0}^{y} r(u) du} dy \\
\leq&\ e^{\eta^2/(c_A-1)} \int_{x}^{\eta} \frac{2}{a(y)} dy \\
\leq&\ e^{\eta^2/(c_A-1)} \frac{2\eta}{(c_A - 1)\mu}.
\end{align*}
Furthermore,
\begin{align*}
\frac{q(\eta)}{q(x)} = e^{\int_{x}^{\eta} r(u) du } \leq 1.
\end{align*}
Hence when $x \in [0,\eta]$,
\begin{align}
\frac{1}{q(x)} \int_{x}^{\infty} \frac{2}{a(y)}q(y) dy \leq e^{\eta^2/(c_A-1)} \frac{2\eta}{c_A - 1}\frac{1}{\mu} + \frac{1}{\mu \abs{\eta}}
.
\label{eq:int-2-first}
\end{align}
This proves the first part of~\eqref{eq:int-2}.

Lastly, we are going to bound $\frac{1}{q(x)} \int_{-\infty}^{x} \frac{2}{a(y)}q(y) dy$ for $x \in [0,-\zeta]$. As before, observe that
\begin{align*}
\frac{1}{q(x)} \int_{-\infty}^{x} \frac{2}{a(y)}q(y) dy =&\ \frac{q(0)}{q(x)} \frac{1}{q(0)} \int_{-\infty}^{0} \frac{2}{a(y)}q(y) dy +  \frac{1}{q(x)}\int_{0}^{x} \frac{2}{a(y)}q(y) dy \\
\leq&\ \frac{q(0)}{q(x)} \frac{1}{\mu} \left(1 + \frac{2}{c_A-1}e^{\frac{1}{c_A-1}}\right) +  \frac{2}{(c_A-1)\mu} \frac{1}{q(x)}\int_{0}^{x} q(y) dy \\
\leq&\ e^{\frac{1}{c_A-1}\zeta^2} \frac{1}{\mu} \left(1 + \frac{2}{c_A-1}e^{\frac{1}{c_A-1}}\right) +  \frac{2}{(c_A-1)\mu} e^{\frac{1}{c_A-1}\zeta^2}\abs{\zeta},
\end{align*}
where the last inequality comes from 
\begin{align*}
\frac{q(0)}{q(x)} &= \exp\left( \int^x_0 -r(y)\mathrm dy \right)
\leq \exp \left(\int^x_0 \frac{2}{c_A-1}y\mathrm dy\right)
= e^{\frac{1}{c_A-1}x^2} .
\end{align*}
This proves the last part of~\eqref{eq:int-1}.

Therefore,~\eqref{eq:int-1} and~\eqref{eq:int-2} hold true.

Note that when $\abs{\zeta} \geq 1$, taking $\eta = 1$ in~\eqref{eq:int-2-first} gives the first part of~\eqref{eq:int-22}, while~\eqref{eq:b2} gives the second part of it. This proves~\eqref{eq:int-22}.

We move on to~\eqref{eq:int-3}. For $x \leq 0$,
 \begin{align*}
 \frac{1}{q(x)} \int_{-\infty}^{x} \frac{2\abs{y}}{a(y)}q(y) dy = \frac{1}{\mu} \frac{1}{q(x)} \int_{-\infty}^{x} r(y)q(y) dy \leq \frac{1}{\mu},
 \end{align*}
 where the last inequality comes from 
 \begin{align*}
 \int^x_{-\infty} r(y) q(y)\mathrm dy &= \int^x_{-\infty} r(y) e^{\int^y_0 r(u)\mathrm du}
 \mathrm dy 
 = q(x) - q(-\infty).
 \end{align*}
 
 When $x \in [0,-\zeta]$,
\begin{align*}
\frac{1}{q(x)} \int_{-\infty}^{x} \frac{2\abs{y}}{a(y)}q(y) dy =&\ \frac{1}{\mu}\frac{q(0)}{q(x)} \frac{1}{q(0)} \int_{-\infty}^{0} r(y)q(y) dy - \frac{1}{\mu}\frac{1}{q(x)} \int_{0}^{x} r(y)q(y) dy \\
=&\ 
\frac{1}{\mu} \frac{q(0)}{q(x)} \frac{1}{q(0)}(q(0) - q(-\infty))
 - \frac{1}{\mu}\frac{1}{q(x)}(q(x) - q(0))
\\
\leq&\ \frac{2}{\mu}\frac{q(0)}{q(x)} = \frac{2}{\mu} e^{\int^x_0 \abs{r(y)}}\mathrm dy
\leq  \frac{2}{\mu}e^{\frac{1}{c_A-1} \zeta^2}.
\end{align*}
This concludes the proof of~\eqref{eq:int-3}.

We proceed to~\eqref{eq:int-4}. For $x \in [0,-\zeta]$,
\begin{align*}
\frac{1}{q(x)} \int_{x}^{\infty} \frac{2\abs{y}}{a(y)}q(y) dy =&\ -\frac{1}{\mu} \frac{1}{q(x)} \int_{x}^{\abs{\zeta}} r(y)q(y) dy -  \frac{1}{\mu\abs{\zeta}} \frac{1}{q(x)} \int_{\abs{\zeta}}^{\infty} y r(y)q(y) dy \\
=&\ \frac{1}{\mu}\Big(1 - \frac{q(\abs{\zeta})}{q(x)}\Big) -  \frac{1}{\mu\abs{\zeta}} \frac{1}{q(x)} \int_{\abs{\zeta}}^{\infty} y r(y)q(y) dy\\
=&\ \frac{1}{\mu}\Big(1 - \frac{q(\abs{\zeta})}{q(x)}\Big) -  \frac{1}{\mu\abs{\zeta}} \frac{1}{q(x)} \Big[-\abs{\zeta}q(\abs{\zeta}) - \int_{\abs{\zeta}}^{\infty}q(y) dy \Big] \\
=&\ \frac{1}{\mu}  +  \frac{1}{\mu\abs{\zeta}} \frac{1}{q(x)} \int_{\abs{\zeta}}^{\infty}q(y) dy \\
\leq&\ \frac{1}{\mu} +  \frac{1}{\mu\abs{\zeta}} \frac{c_A - \mu + \delta(1 - \mu)\abs{\zeta} + \mu \zeta^2}{2\abs{\zeta}} \frac{q(\abs{\zeta})}{q(x)} \\
\leq&\ \frac{3}{2}\frac{1}{\mu} +  \frac{1}{\mu\abs{\zeta}} \frac{c_A + \delta\abs{\zeta}}{2\abs{\zeta}},
\end{align*}
where in the last inequality we use the fact that for $x \geq -\zeta$,
\begin{align}
\frac{1}{q(x)}\int^\infty_x q(y)\mathrm dy &= \int^\infty_x e^{\int^y_x r(u)\mathrm du}
\mathrm dy  
= \int^\infty_0 e^{\frac{-2\abs{\zeta}}{c_A - \mu + \delta(1-\mu) \abs {\zeta} + \mu\zeta^2} y}
\mathrm dy \notag \\
&= \frac{c_A - \mu + \delta(1-\mu)\abs{\zeta} + \mu\zeta^2}{2\abs{\zeta}} ,
\label{eq:interm-1}
\end{align}
and
\begin{align*}
\frac{q(\abs{\zeta})}{q(x)} = e^{\int_{x}^{\abs{\zeta}} r(u) du} \leq 1.
\end{align*}
For $x \geq -\zeta$,
\begin{align*}
\frac{1}{q(x)} \int_{x}^{\infty} \frac{2\abs{y}}{a(y)}q(y) dy =&\ -  \frac{1}{\mu\abs{\zeta}} \frac{1}{q(x)} \int_{x}^{\infty} y r(y)q(y) dy \\
=&\  -  \frac{1}{\mu\abs{\zeta}} \frac{1}{q(x)} \Big[-x q(x) - \int_{x}^{\infty}q(y) dy \Big] \\
=&\   \frac{x}{\mu\abs{\zeta}} + \frac{1}{\mu\abs{\zeta}} \frac{1}{q(x)} \int_{x}^{\infty}q(y) dy \\
=&\ \frac{x}{\mu\abs{\zeta}} +  \frac{1}{\mu\abs{\zeta}} \frac{c_A - \mu + \delta(1-\mu)\abs{\zeta} + \mu\zeta^2}{2\abs{\zeta}}.
\end{align*}
This proves~\eqref{eq:int-4}.

Finally, we deal with~\eqref{eq:int-5} and~\eqref{eq:int-6}.
For $x < 0$ we use \eqref{eq:b1} to see that
\begin{align*}
\frac{\abs{r(x)}}{q(x)}\int_{-\infty}^{x}\frac{2}{a(y)} q(y)dy \leq \frac{\abs{r(x)}}{b(x)} = \frac{2}{a(x)} \leq \frac{2}{(c_A - 1)\mu}.
\end{align*}
Similarly, we invoke \eqref{eq:b2} to see that when $x \geq 0$,
\begin{align*}
\frac{\abs{r(x)}}{q(x)}\int_{x}^{\infty}\frac{2}{a(y)} q(y)dy \leq \frac{\abs{r(x)}}{\abs{b(x)}} = \frac{2}{a(x)} \leq \frac{2}{(c_A - 1)\mu}.
\end{align*}
This proves~\eqref{eq:int-5} and~\eqref{eq:int-6}, concluding our proof of Lemma~\ref{lem:ints_for_grad_bdds}.
\end{proof}

\begin{proof}[Proof of Lemma~\ref{lem:ey}]
Consider the Lyapunov function $V(x) = x^2$. Using the form of
$G_Y$ in~\eqref{eq:ggy}, we sees immediately that 
\begin{align}
G_Y V(x) &= 2x b(x) + a(x).
\label{eq:gy_vx} 
\end{align}

Recall the form of $b(x)$ and $a(x)$.
When $x\leq -1/\delta$, 
\begin{align*}
G_Y V(x) &= 
2x (-\mu x) + \mu \left[(c_A - 1) + \Lambda\right] 
= - 2 \mu x^2 + \mu(c_A - 1) + \mu^2 \delta^{-2},
\end{align*}
where for the last equality we use $\delta = 1/\sqrt{R}$.

When $ x\in [-1/\delta, -\zeta]$,
\begin{align*}
G_Y V(x) &
= - 2 \mu (1 - \mu/2) x^2 + \delta(1 - \mu ) \mu x 
+ \mu(c_A - \mu) \\
&\leq - 2 \mu (1 - \mu/2) x^2
+ \mu \frac{1-\mu}{2} x^2 + \mu\frac{1-\mu}{2} \delta^2 
+\mu(c_A - \mu)\\
&= -\frac{1}{2}(3-\mu) \mu x^2 + \mu\frac{1-\mu}{2} \delta^2 
+\mu(c_A - \mu)\\
&\leq -\mu x^2 
+ \mu\frac{1-\mu}{2} \delta^2 
+\mu(c_A - \mu),
\end{align*} 
where to get the last inequality we use $\mu < 1$.

When $x\geq -\zeta$, 
\begin{align*}
G_Y V(x) 
&= - 2 \mu \abs \zeta x + \mu(c_A -\mu) +
\mu\left[\delta(1-\mu) \abs{\zeta} + \mu \zeta^2\right].
\end{align*}

Therefore,
\begin{align*}
G_YV(x) &\leq -2\mu x^2 \mathds{1}_{\{x < -1/\delta\}}
-\mu x^2 \mathds{1}_{\{x \in [-1/\delta,-\zeta)\}}
- 2\mu \abs{\zeta} x \mathds{1}_{\{x \geq -\zeta\}}
\\
&+ \mu c_A + 
\mu^2 \delta^{-2} \mathds{1}_{\{x<-1/\delta\}}
+ \frac{1-\mu}{2} \delta^2 \mu \mathds{1}_{\{x\in [-1/\delta,-\zeta)\}} 
\\
&+ 
\delta(1-\mu) \mu \abs{\zeta}\mathds{1}_{\{x \geq -\zeta\}} + \mu^2 \zeta^2 \mathds{1}_{\{x \geq -\zeta\}}.
\end{align*}

According to the standard Foster-Lyapunov criterion~\cite{MeynTwee1993bb}, for any $U,g_1,g_2: \mathbb{R}\rightarrow \mathbb{R_+}$ satisfying
\begin{align*}
G_Y U(x) &\leq -g_1(x) + g_2(x), \quad x\in \mathbb{R},
\end{align*}
there is 
\begin{align*}
\E g_1(Y_\infty) &\leq \E g_2 (Y_\infty) .
\end{align*}

Thus, 
\begin{align}
&2\E \left[Y_\infty^2 \mathds{1}_{\{Y_\infty<-1/\delta\}}\right]
+  \E \left[Y_\infty^2 \mathds{1}_{\{Y_\infty\in [-1/\delta.-\zeta)\}}\right]
+ 2\abs{\zeta} \E \left[Y_\infty \mathds{1}_{\{Y_\infty \geq -\zeta\}}
\right] 
\notag\\
&\hspace{-0.5cm}\leq c_A  + \frac{1-\mu}{2} \delta^2
+
\mu \delta^{-2} \Prob(Y_\infty < -1/\delta) 
+ (1-\mu)\delta\abs{\zeta}\Prob(Y_\infty \geq -\zeta) 
+ \mu \zeta^2 \Prob(Y_\infty \geq -\zeta).
\label{eq:fl}
\end{align}

Note that 
\begin{align*}
\delta^{-1} \Prob(Y_\infty < -1/\delta)&\leq
\E \left[\abs{Y_\infty}\mathds{1}_{\{Y_\infty < -1/\delta\}}\right] 
,\quad
\abs{\zeta}\Prob(Y_\infty \geq -\zeta ) \leq \E\left[ Y_\infty \mathds{1}_{\{Y_\infty \geq -\zeta\}}\right] .
\end{align*}

Applying the above inequalities to~\eqref{eq:fl}, we get that
\begin{align}
&2\E \left[Y_\infty^2 \mathds{1}_{\{Y_\infty<-1/\delta\}}\right] 
+  \E \left[Y_\infty^2 \mathds{1}_{\{Y_\infty\in [-1/\delta.-\zeta)\}}
\right] 
+ 2\abs{\zeta} \E \left[Y_\infty \mathds{1}_{\{Y_\infty \geq -\zeta\}}
\right] 
\notag\\
&\hspace{-0.5cm}\leq c_A  + \frac{1-\mu}{2} \delta^2
+ \mu \delta^{-1} \E\left[ \abs{Y_\infty}\mathds{1}_{\{Y_\infty < -1/\delta\}}\right] 
+ \mu \abs{\zeta} \E \left[Y_\infty \mathds{1}_{\{Y_\infty \geq -\zeta\}}\right] 
\notag\\
&\hspace{-0.5cm}+ 
(1-\mu) \delta^2 \mathds{1}_{\{\abs{\zeta}\leq \delta\}}
+ (1-\mu)\delta \E\left[ \abs{Y_\infty}\mathds{1}_{\{Y_\infty
\geq -\zeta\}}\right] \mathds{1}_{\{\abs \zeta>\delta\}} . 
\label{eq:fl_1}
\end{align}

Since 
\begin{align*}
2\abs \zeta - \mu \abs \zeta - (1-\mu)\delta \mathds{1}_{\{\abs \zeta>\delta\}}
&\geq (2-\mu)\abs \zeta -  (1-\mu) \abs \zeta 
= \abs \zeta ,
\end{align*}
~\eqref{eq:fl_1} implies that 
\begin{align}
 &2\E \left[Y_\infty^2 \mathds{1}_{\{Y_\infty<-1/\delta\}}\right] 
 +  \E \left[Y_\infty^2 \mathds{1}_{\{Y_\infty\in [-1/\delta.-\zeta)\}}
\right] 
+ \abs \zeta \E \left[Y_\infty \mathds{1}_{\{Y_\infty \geq -\zeta\}}
\right] 
\notag\\
&\hspace{-0.5cm}\leq c_A + \frac{3}{2} (1-\mu) \delta^2 
+ \mu\delta^{-1} \E\left[\abs{Y_\infty} \mathds{1}
_{\{Y_\infty < -1/\delta\}}\right] .
\label{eq:fl_2}
\end{align}

From Jensen's inequality and~\eqref{eq:fl_2}, we have
\begin{align*}
\E\left[\abs{Y_\infty} \mathds{1}
_{\{Y_\infty < -1/\delta\}}\right]
&\leq \sqrt{\E \left[Y_\infty^2 \mathds{1}_{\{Y_\infty<-1/\delta\}}\right]}
\\
&\leq \sqrt{\frac{1}{2}\left[
c_A + \frac{3}{2} (1-\mu) \delta^2 
+ \mu\delta^{-1} \E\left[\abs{Y_\infty} \mathds{1}
_{\{Y_\infty < -1/\delta\}}\right] 
\right]},
\end{align*}
which is equivalent to a quadratic inequality in $\E\left[\abs{Y_\infty} \mathds{1}
_{\{Y_\infty < -1/\delta\}}\right]$, 
\begin{align*}
2 \left\{\E\left[\abs{Y_\infty} \mathds{1}
_{\{Y_\infty < -1/\delta\}}\right]\right\}^2 
- \mu \delta^{-1} \E\left[\abs{Y_\infty} \mathds{1}
_{\{Y_\infty < -1/\delta\}}\right] 
- \left[
c_A + \frac{3}{2} (1-\mu) \delta^2 \right]
&\leq 0.
\end{align*}

Solving the above quadratic inequality gives 
\begin{align}
\E\left[\abs{Y_\infty} \mathds{1}
_{\{Y_\infty < -1/\delta\}}\right]
&\leq \frac{1}{4} \mu \delta^{-1}
+\frac{1}{4} \sqrt{\mu^2 \delta^{-2}
+ 8\left[
c_A + \frac{3}{2} (1-\mu) \delta^2 \right]
}
\notag\\
&\leq \frac{1}{2} \mu \delta^{-1}
+\frac{1}{2}\sqrt{2c_A + 3(1-\mu)\delta^2}
\notag\\
&\leq \frac{1}{2} \mu \delta^{-1}
+c_A + \frac{3}{2}(1-\mu)\delta^2 ,
\label{eq:ey_1}
\end{align}
where for the second inequality we use the fact that for any two
non-negative real numbers $y$ and $z$, $\sqrt{y+z}\leq \sqrt{y}
+\sqrt{z}$, and for the last inequality we use $c_A \geq 1$. 

Substituting~\eqref{eq:ey_1} back into~\eqref{eq:fl_2}, we have
\begin{align}
 &2\E \left[Y_\infty^2 \mathds{1}_{\{Y_\infty<-1/\delta\}}\right] 
 +  \E \left[Y_\infty^2 \mathds{1}_{\{Y_\infty\in [-1/\delta.-\zeta)\}}
\right] 
+ \abs \zeta \E \left[Y_\infty \mathds{1}_{\{Y_\infty \geq -\zeta\}}
\right] 
\notag\\
&\hspace{-0.5cm}\leq c_A + \frac{3}{2} (1-\mu) \delta^2 
+\mu \delta^{-1} \left[c_A + \frac{3}{2}(1-\mu)\delta^2\right]
+\frac{1}{2} \mu^2 \delta^{-2} 
.
\label{eq:fl_4}
\end{align}

Then, applying Jensen's inequality to~\eqref{eq:fl_4} implies that
\begin{align}
 \E \left[\abs{Y_\infty} \mathds{1}_{\{Y_\infty\in [-1/\delta.-\zeta)\}}
\right] 
&\leq \sqrt{ \E \left[Y_\infty^2 \mathds{1}_{\{Y_\infty\in [-1/\delta.-\zeta)\}}
\right] }
\notag\\
&\leq \sqrt{1 + \mu\delta^{-1}} \sqrt{c_A + \frac{3}{2}(1-\mu)
\delta^2}
+\sqrt{\frac{1}{2} \mu^2 \delta^{-2}}
\notag\\
&\leq \left( 1 + \mu^{1/2}\delta^{-1/2}\right)
\left[c_A + \frac{3}{2}(1-\mu)
\delta^2\right]
+ \frac{\sqrt{2}}{2} \mu \delta^{-1}
 .
\label{eq:ey_2}
\end{align}

Finally,
\begin{align}
\E\left[Y_\infty \mathds{1}_{\{Y_\infty \geq -\zeta\}}\right]
&\leq \frac{1}{\abs \zeta}
\left\{(1 + \mu\delta^{-1})
\left[c_A + \frac{3}{2}(1-\mu) \delta^2\right]
+\frac{1}{2}\mu^2 \delta^{-2}
\right\} .
\label{eq:ey_3}
\end{align}

Adding up~\eqref{eq:ey_1},~\eqref{eq:ey_2}, and~\eqref{eq:ey_3} proves Lemma~\ref{lem:ey}.  
\end{proof}

\begin{proof}[Proof of Lemma~\ref{lem:a}]

We start with~\eqref{eq:abdd_1}.

From the form of $a(x)$ in~\eqref{eq:aform-1},~\eqref{eq:abdd_1} is obviously true for $x\leq -1/\delta$ and $x\geq -\zeta$. 

When $x\in [-1/\delta,-\zeta]$, 
\begin{align*}
a(x) &= \mu\left(c_A - \mu + \delta(1-\mu)x + \mu x^2\right)
= \mu \left\{
\mu \left(x -x_0\right)^2
+ (c_A -\mu) - \mu x_0^2
\right\},
\end{align*}
where $x_0 = -\frac{\delta(1-\mu)}{2\mu}$.

	If $x_0 \leq -1/\delta$, $a(x)$ is increasing with $x$ over the interval $[-1/\delta,-\zeta]$. Thus 
	\begin{align*}
	\min_{x\in [-1/\delta,-\zeta]}a(x) \geq a(-1/\delta) > (c_A-1)\mu.
	\end{align*}
	
	Otherwise, we have $\delta^2 < \frac{2\mu}{1-\mu}$.
	
	\begin{align*}
	\min_{x\in [-1/\delta,-\zeta]}a(x) &= a\left(x_0\right)
	= \mu(c_A - \mu) - \mu^2 \frac{\delta^2(1-\mu)^2 }{4\mu^2}
	\\
	&> \mu(c_A - \mu) - \frac{1}{4}(1 - \mu)^2\frac{2\mu}{1-\mu}
	\\
	&= \mu(c_A - 1/2 - \mu/2)
	\\
	&>(c_A - 1)\mu ,
	\end{align*}
where to get the last inequality we use $\mu < 1$.

In this way we have established~\eqref{eq:abdd_1}.

Next we proceed to~\eqref{eq:abdd_2}. Note that
\begin{align}
a'(x) &= \mu \left[ 2\mu x + \delta(1-\mu)\right] \mathds{1}_{\{
x\in (-1/\delta,-\zeta]
\}}.
\label{eq:aprime}
\end{align}
Hence 
\begin{align}
\frac{\abs{xa'(x)}}{a(x)}&= \frac{\abs{  2\mu x^2 + \delta(1-\mu)x }}{  \mu x^2+ \delta(1-\mu)x+c_A - \mu }
\mathds{1}_{\{
x\in (-1/\delta,-\zeta]
\}}.
\label{eq:abdd_2_ls}
\end{align}

Consider the function 
\begin{align*}
g_1(x) &= 2 \left[ \mu x^2+ \delta(1-\mu)x+c_A - \mu \right]
- \left[ 2\mu x^2 + \delta(1-\mu)x\right]
= \delta(1-\mu)x + 2(c_A - \mu) 
. 
\end{align*}
When $x\in [-1/\delta,-\zeta]$,
\begin{align}
g_1(x) &\geq - (1 - \mu) + 2(c_A - \mu) 
= (c_A - \mu) + (c_A - 1) > 0
,
\label{eq:abdd_2_g1}
\end{align}
where to get the last inequality we use $\mu < 1$ and $c_A > 1$.

Consider another function
\begin{align*}
g_2(x) &= a \left[ \mu x^2+ \delta(1-\mu)x+c_A - \mu \right]
+ \left[ 2\mu x^2 + \delta(1-\mu)x\right]
\\
&= (a+2) \mu x^2 + (a+1) \delta(1 - \mu) x + a (c_A-\mu)
,
\end{align*}
where $a$ is a positive constant to be determined later.

When $x\in [-1/\delta,-\zeta]$,
\begin{align*}
g_2(x) &\geq -(a+1) (1 - \mu) + a (c_A - \mu)
= (c_A - 1) a - (1 - \mu)
.
\end{align*} 
Taking $a = \frac{1-\mu}{c_A - 1}$, we have 
\begin{align}
g_2(x) &\geq 0,\quad x\in [-1/\delta,-\zeta] .
\label{eq:abdd_2_g2}
\end{align}

Applying~\eqref{eq:abdd_2_g1} and~\eqref{eq:abdd_2_g2} to~\eqref{eq:abdd_2_ls} proves~\eqref{eq:abdd_2}. 

Now we deal with~\eqref{eq:abdd_3}. 

From~\eqref{eq:aform-1} and~\eqref{eq:aprime}, 
\begin{align}
\frac{\abs{a'(x)}}{a(x)}&= \frac{\abs{2\mu x + \delta(1 - \mu)}}{\mu x^2 + \delta(1-\mu)x + (c_A - \mu)}
\mathds{1}_{\{
x\in (-1/\delta,-\zeta]
\}}.
\label{eq:abdd_3_ls}
\end{align}

Consider the function
\begin{align}
g_3(x) &= a\left[\mu x^2 + \delta(1-\mu)x + (c_A - \mu)\right] 
- \left[2\mu x + \delta(1 - \mu)\right]
\notag\\
&= a\mu x^2 + \left[\delta a (1 - \mu) - 2\mu\right] x
+ a (c_A -\mu) -\delta(1 - \mu),
\notag  
\end{align}
where $a$ is a positive constant to be determined later. 

For $\mu \in (0,1)$ such that $\delta a (1 - \mu) - 2\mu \geq 0$,
when $x\in [-1/\delta,-\zeta]$,
\begin{align*}
g_3(x) &\geq -a(1-\mu) + 2 \mu \delta^{-1} + a (c_A -\mu) -\delta(1 - \mu)\\
&\geq -a(1-\mu)  +a(c_A-\mu) - (1 -\mu)
\\
&= (c_A - 1) a - (1 - \mu) ,
\end{align*}
where for the second inequality we used $\delta\leq 1$.

For $\mu \in (0,1)$ such that $\delta a (1 - \mu) - 2\mu < 0$, let $x_0 = \frac{2\mu -
\delta a (1 - \mu)}{2 a \mu} > 0$.
\begin{align*}
g_3(x) &= a\mu (x - x_0)^2 + a(c_A - \mu) 
-\delta(1 - \mu) - a\mu x_0^2 \\
&\geq a(c_A - \mu) 
-\delta(1 - \mu) - a\mu \left(\frac{2\mu}{2 a\mu }\right)^2 \\
&\geq a(c_A - \mu) 
-(1 - \mu) - \mu/a \\
&\geq \left[(c_A - 1)a - 1\right] + \mu (1 - 1/a).
\end{align*}

Taking $a = \frac{1}{c_A - 1}\vee 1$, we have that
\begin{align}
g_3(x) &\geq 0,\quad x\in [-1/\delta,-\zeta] .
\label{eq:abdd_3_g3}
\end{align}

Consider another function
\begin{align*}
g_4(x) &= a\left[\mu x^2 + \delta(1-\mu)x + (c_A - \mu)\right] 
+ \left[2\mu x + \delta(1 - \mu)\right] \\
&= a\mu x^2 +
\left[\delta a (1 - \mu) + 2 \mu\right] x
+ a (c_A - \mu) +\delta(1-\mu) ,
\end{align*}
where $a$ is a positive constant to be determined later.

Denote $x_0 = -\frac{\delta a (1-\mu)+2\mu}{2 a \mu}$. If $x_0 \leq -1/\delta$, $g_4(x)$ is increasing with $x$ over the interval $[-1/\delta,-\zeta
]$. Then, for $x \in [-1/\delta,-\zeta
]$,
\begin{align*}
g_4(x) &\geq g_4(-1/\delta) \\
&= a \mu \delta^{-2} - a (1 -\mu) - 2 \mu \delta^{-1}
+ a (c_A -\mu) + \delta(1 -\mu) \\
&\geq 
(a - 2) \mu \delta^{-1} + (c_A - 1 )a .
\end{align*}
Otherwise, $x_0 > -1/\delta$, that is, 
\begin{align}
\delta a (1-\mu)+2\mu < 2 a \mu \delta^{-1}.
\label{eq:abdd_3_interm}
\end{align}
Then,
\begin{align*}
g_4(x)&= a\mu (x - x_0)^2 + a (c_A -\mu)
+\delta(1 -\mu) - a\mu x_0^2 \\
&\geq  a (c_A -\mu)
+\delta(1 -\mu) - a\mu
 \left(\frac{\delta a (1-\mu)+2\mu}{2 a \mu}\right)^2
 \\
 &> a (c_A -\mu)
+\delta(1 -\mu) - \frac{1}{2}\delta^{-1}
\left[\delta a (1-\mu)+2\mu\right] \\
&\geq a (c_A -\mu)
 - \frac{1}{2}\delta^{-1}
\left[\delta a (1-\mu)+2\mu\right] \\
&= c_A a -\frac{1}{2}(1 + \mu)a -\mu\delta^{-1} \\
&\geq (c_A - 1) a -\mu\delta^{-1} ,
\end{align*}
where to get the third line we use~\eqref{eq:abdd_3_interm}.

Taking $a = \frac{\mu \delta^{-1}}{c_A - 1}\vee
2$, we have that 
\begin{align}
g_4(x)&\geq 0,\quad x\in [-1/\delta,-\zeta] .
\label{eq:abdd_3_g4}
\end{align}

Applying~\eqref{eq:abdd_3_g3} and~\eqref{eq:abdd_3_g4} to~\eqref{eq:abdd_3_ls},
we get that
\begin{align}
\frac{\abs{a'(x)}}{a(x)} &\leq \left(
\frac{\mu \delta^{-1}}{c_A - 1}
\vee
\frac{1}{c_A - 1}
\vee 2 \right)
\mathds{1}_{\{x\in (-1/\delta,-\zeta
]\}}.
\label{eq:abdd_3_0}
\end{align}
Using~\eqref{eq:abdd_3_0} along with~\eqref{eq:ey}
in Lemma~\ref{lem:ey} proves~\eqref{eq:abdd_3}.

Finally we approach~\eqref{eq:abdd_4}. 

Note that when $x \leq -1/\delta$,
\begin{align}
\abs{r'(x) a(x)} = 2\mu ,
\label{eq:abdd_4_1}
\end{align}
and when $x > -\zeta$,
\begin{align}
\abs{r'(x) a(x)} = 0 .
\label{eq:abdd_4_2}
\end{align}

When $x \in (-1/\delta,-\zeta]$, 
\begin{align*}
r'(x) a(x) 
&= \frac{2}{a(x)} 
\left[b'(x) a(x) -  b(x) a'(x)\right] \\
&= -2 \mu + 2\mu \frac{xa'(x)}{a(x)} .
\end{align*}
Hence 
\begin{align*}
\abs{r'(x) a(x)} &\leq
2\mu \left(1 + \frac{\abs{xa'(x)}}{a(x)}\right) 
.
\end{align*}
Applying~\eqref{eq:abdd_2} to the right side of the inequality above, along with~\eqref{eq:abdd_4_1}
and~\eqref{eq:abdd_4_2}, proves~\eqref{eq:abdd_4}.
\end{proof}

\renewcommand{\baselinestretch}{0.9}
\footnotesize
\def\cprime{$'$} \def\cprime{$'$} \def\cprime{$'$} \def\cprime{$'$}
  \def\cprime{$'$} \def\cprime{$'$} \def\cprime{$'$}

\end{document}